\documentclass[11pt]{amsart}

\usepackage[T1]{fontenc}
\usepackage[utf8]{inputenc}
\usepackage{lmodern}
\usepackage{amsmath,amssymb,amsthm}
\usepackage{mathtools}
\usepackage{enumitem}
\usepackage[margin=1.15in]{geometry}
\usepackage[table]{xcolor}
\usepackage{tikz-cd}
\usepackage{longtable}
\usepackage{graphicx}
\usepackage{float}
\usepackage[hidelinks]{hyperref}

\theoremstyle{plain}
\newtheorem{theorem}{Theorem}[section]
\newtheorem{proposition}[theorem]{Proposition}
\newtheorem{lemma}[theorem]{Lemma}
\newtheorem{corollary}[theorem]{Corollary}
\theoremstyle{definition}
\newtheorem{definition}[theorem]{Definition}
\newtheorem{example}[theorem]{Example}
\theoremstyle{remark}
\newtheorem{remark}[theorem]{Remark}

\setlist[enumerate]{label=\textup{(\roman*)}, leftmargin=*, widest=iii}

\newcommand{\ZZ}{\mathbb{Z}}
\newcommand{\NN}{\mathbb{N}}
\newcommand{\A}{\mathrm{A}}
\newcommand{\D}{\mathrm{D}}
\newcommand{\cO}{\mathcal{O}}
\newcommand{\cA}{\mathcal{A}}
\newcommand{\cE}{\mathcal{E}}

\newcommand{\cJ}{\mathcal{J}}
\newcommand{\cN}{\mathcal{N}}
\newcommand{\Spec}{\operatorname{Spec}}
\newcommand{\Sym}{\operatorname{Sym}}
\newcommand{\Hom}{\operatorname{Hom}}
\newcommand{\Lie}{\operatorname{Lie}}
\newcommand{\Bl}{\operatorname{Bl}}
\newcommand{\Sch}{\mathrm{Sch}}
\newcommand{\gp}{\mathrm{gp}}
\newcommand{\fppf}{\mathrm{fppf}}
\newcommand{\id}{\mathrm{id}}
\newcommand{\colim}{\operatorname*{colim}}
\newcommand{\Homeq}[1]{\Hom^{#1}}
\newcommand{\mon}[1]{[#1\rangle}
\newcommand{\units}[1]{#1^{\times}}

\makeatletter
\newcommand{\leankey}[4][]{%
  \expandafter\edef\csname leankey@decl@#2\endcsname{\detokenize{#3}}%
  \expandafter\edef\csname leankey@file@#2\endcsname{\detokenize{#4}}%
  \if\relax\detokenize{#1}\relax
    \expandafter\def\csname leankey@disp@#2\endcsname{\csname leankey@decl@#2\endcsname}%
  \else
    \expandafter\def\csname leankey@disp@#2\endcsname{#1}%
  \fi}
\newcommand{\leansrc}[1]{%
  \ifcsname leankey@decl@#1\endcsname
    \par\nopagebreak
    {\footnotesize\raggedleft\noindent Lean:
      \texttt{\leanbreakable{\csname leankey@disp@#1\endcsname}}\par}
  \else
    \PackageWarning{leansrc}{Undefined Lean key `#1'}%
    \par\noindent\hfill{\footnotesize\textbf{Lean key #1 undefined}}\par
  \fi}
\makeatother
\ExplSyntaxOn
\tl_new:N \l__lean_tl
\cs_new_protected:Npn \leanbreakable #1
  {
    \tl_set:Nx \l__lean_tl { #1 }
    \regex_replace_all:nnN { ([\_\.]) } { \1 \c{allowbreak} } \l__lean_tl
    \tl_use:N \l__lean_tl
  }
\ExplSyntaxOff
\leankey{site}{BigFppfSite, BigFppfOver}{StacksFrontier (dependency)}
\leankey{algspace}{IsAlgebraicSpace, IsSeparatedOver}{StacksFrontier (dependency), Sec3Rest.lean}
\leankey{monoidRing}{X, zRing}{Sec2Monoids.lean, GeneralBase.lean}
\leankey{bialgebra}{counit_X, comul_X, antipode_X}{Sec2Monoids.lean}
\leankey{diagMonoid}{AZ, A, monObjA}{GeneralBase.lean}
\leankey{diagAffine}{flat_str, fact26_i_sheafification, isPullback_AZ_prod, isPullback_A_prod, prodIso}{Sec3Leftover26.lean, FixesFact26.lean, GeneralBase.lean}
\leankey{diagFP}{fact210_i, finitePresentation, instLFP_AZ, instLFP_str, lfp_mul, instAffine_str}{Sec3Remaining.lean, Sec2Monoids.lean, GeneralBase.lean}
\leankey{face}{IsFaceOver}{Sec2Monoids.lean}
\leankey{faceCriterion}{isFace_iff}{Sec2Monoids.lean}
\leankey{unitsFace}{unitsFace, unitsFace_isGroup, faceCond_unitsFace, le_unitsFace, unitsFace_le_face}{Sec2Monoids.lean}
\leankey{facePushout}{faceCond_gluedMagnet, isPushout_glued}{Sec2Monoids.lean, GeneralBasePushout.lean}
\leankey{faceWeight}{exists_weight_of_isFace}{Sec12Topology.lean}
\leankey{action}{AAction, EquivariantHom}{Foundations.lean, Sec3Functorial.lean}
\leankey{actionAN}{coact, aAction, Amap, restrict, faceA, faceExt}{Sec2Monoids.lean, GeneralBase.lean, Sec3Functorial.lean, GeneralBasePushout.lean}
\leankey{fixedPoints}{fixedPoints}{Sec6FPR.lean}
\leankey{attractor}{AttractorPt, attractor}{Foundations.lean}
\leankey{attractorSheaf}{isSheaf_attractor, isSheaf_attractor_of_isAlgebraicSpace}{Foundations.lean, AlgebraicSpace.lean}
\leankey{attractorExtreme}{attractorTopEquiv, attractorBotEquiv}{Foundations.lean, Sec3Functorial.lean}
\leankey{attractorLimits}{fibreProductAction, existsUnique_fibreLift, attractorFibreProductIso, attractorProdIso, attractorOverIso}{Sec3Functorial.lean}
\leankey{magnetMaps}{magnetMap, faceMap, magnetMap_comp_faceMap, attractorMap, mono_attractorMap}{Sec3Functorial.lean}
\leankey{faceClosed}{fact316_spaces}{Sec3LeftoverFace.lean}
\leankey{attractorLfp}{fact319_relativelyLimitPreserving, isLocallyOfFinitePresentation_repMap_spaces}{Final3Fact319.lean, Sec11SpacesFact.lean}
\leankey{attractorActions}{attractorAction, magnetAction, attractorAction_smul_eq, attractorMonoid, attractorGroup}{Sec3Functorial.lean}
\leankey{separatedMono}{ofLE_injective, attractor_toPoint_injective_sep}{Sec5Etale.lean}
\leankey{colimitMagnets}{prop337_exists, prop337_colimit}{Prop337General.lean}
\leankey{subgroupFixed}{attractorFixedEquivOfCover, isFreeAlong_subtype, extSieve_mem, prop335_of_atlas}{Sec8Representability.lean, Prop81General.lean}
\leankey{prescribed}{prescribedLimit}{Sec4Limits.lean}
\leankey[prescribedLimitIso, mono\_prescribedLimit\_$\iota$, prescribedLimit\_toPoint\_injective]{prescribedFibre}{prescribedLimitIso, mono_prescribedLimit_ι, prescribedLimit_toPoint_injective}{Sec4Limits.lean}
\leankey{AMmodule}{ModuleAAction, AModule, aModuleEquiv}{Sec3Remaining.lean, Sec3LeftoverW.lean, Sec3Leftover2W.lean}
\leankey{gradings}{prop221, actionEquivGrading}{Sec3Leftover221Eq.lean, Prop222.lean}
\leankey{affineTheorem}{Zatt, instClosedIAtt, attEquiv, representsAttractor}{GeneralBaseThm320.lean}
\leankey{exampleAL}{ker_iAtt_self}{SelfActionsBase.lean}
\leankey{exampleLinear}{ex12_general, ex12_sub, ex12_mem_iff, ex12_toPoint_injective, ex12_sub_baseChange}{Sec3Leftover12.lean, Ex12Rem313Final.lean}
\leankey{vectorBundle}{prop328}{Sec3Leftover328V.lean}
\leankey{affineConsequences}{instClosedZincl, ker_iAtt_iInf, exists_toPoint_iInf_iff, isPullback_Zmap, instClosedZmap}{GeneralBaseCor.lean}
\leankey{selfAction}{prop42_scheme, lemma41_empty_scheme, prop43_scheme}{SelfActionsMain.lean}
\leankey{changeMonoid}{prop331_i, AlgMag.SD.prop331_space_final, AlgMag.Unrestrict.prop331_ii, prop331_ii_anyS, prop331_iii_gp}{Sec3Final.lean, SDR4Final.lean, Unrestrict1.lean, Prop331Free.lean, Unrestrict2.lean}
\leankey{intersection}{prop332_i, AlgMag.SD.prop332_space_final, prop332_ii, prop332_iii_space_gp}{Sec3Final.lean, SDR4Final.lean, Unrestrict1.lean}
\leankey{etaleFace}{isFaceCartesian_of_fg, isFaceCartesian_of_separated}{Sec5Etale.lean}
\leankey{fpr}{IsFPR, IsZFPR, IsStronglyFPR, IsEquivariantAtlas, IsSAffineAtlas, IsZFPRAtlas, IsStronglyFPRAtlas}{Sec6FPR.lean}
\leankey{fprBasic}{hasFiniteAffineDecomposition_iff, prop63, lemma67_i, lemma67_ii, isZFPR_of_mono, isStronglyFPR_of_mono}{Sec6FPR.lean, AppendixA.lean, Sec6Atlas.lean, Sec7Linearizable.lean}
\leankey{zll}{ZariskiLocallyLinearizable}{Sec7Linearizable.lean}
\leankey{zllAtlas}{fact73}{Sec6Atlas.lean}
\leankey{chart}{ChartLemmaStatement, ChartLemmaStatementFG, chartLemma}{ChartStatement.lean, ChartBridge.lean, ChartFinal.lean}
\leankey{atlas}{thm68_atlas, thm68_family_spaces, thm68_qc_spaces, thm68_i}{FinalNonQCAtlas.lean, ChartFinalSpaces.lean, ChartThm68.lean}
\leankey{strongAtlas}{thm68_atlas, thm68_family_spaces, thm68_qc_spaces}{FinalNonQCAtlas.lean, ChartFinalSpaces.lean}
\leankey{conj69}{conj69_separated}{Conj69Separated.lean}
\leankey{conj69DoubledPlane}{not_conj69_doubledPlane, no_atlas_chart}{Conj69DoubledOrigin.lean}
\leankey{fixedClosed}{prop81FG, prop81_general}{Sec8Representability.lean, Prop81General.lean}
\leankey{representability}{thm83_repr, thm83, thm83_iv, thm83_vi, thm83_vii, thm83_of_lfp}{Sec8Theorem83.lean, Sec8Limits.lean}
\leankey{representabilityFree}{thm83_nonqc, thm83_free}{FinalNonQCAtlas.lean, FinalCorAtlas.lean}
\leankey{finiteMagnet}{prop85_general}{Prop85General.lean}
\leankey{gluedScheme}{prop31, isOpenImmersion_gpMap, prop33_scheme}{XProductsGlue.lean, XProductsAction.lean}
\leankey{monoscheme}{spec, IsMonoscheme}{XProductsF1.lean, XProductsF1Base.lean}
\leankey{attractorFamily}{XN}{XProductsThm.lean}
\leankey{xProduct}{thm41, thm41Equiv, toPoint_injective_gp}{XProductsThm.lean}
\leankey{xProductRep}{prop43_i, prop43_moreover, prop43_ii, prop45}{XProductsRep.lean, XProductsThm.lean}
\leankey{exampleP1}{ex44AttrEquiv, ex44XNEquiv, ex44_iInf, ex44_zero}{XProductsExampleFunctor.lean, XProductsExample.lean, XProductsExampleClass.lean}
\leankey{pi0}{prop121_face, Faces.unitsSub_face_eq, bijective_connectedComponents_repFace, connectedSpace_faceFibre}{Faces/Prop121.lean, Faces/Regrade.lean, Sec12Topology.lean}
\leankey{immersions}{isPullback_repAttractorMap_open, isOpenImmersion_repAttractorMap, isPullback_repAttractorMap, isClosedImm_repAttractorMap}{Sec12Spaces.lean, Sec12Topology.lean}
\leankey{deformation}{deformation_spec, exists_strong_deformation_retract}{Sec12Topology.lean}
\leankey{BB}{Faces.bb_face_zeroSection, Faces.bb_face_attP_zeroSection, Faces.faceChartZero_printed, Faces.bb_face, Faces.bb_face_attP, Faces.exists_chart_cover_face, Faces.attP_smooth_affine_cover, Faces.FaceChart, Faces.faceChart_printed, Faces.face_sym_ring, Faces.mem_face_of_mem_closure, Faces.nbar_sharp, Faces.exists_weight_nbar}{Faces/Spaces.lean, Faces/Cover.lean, Faces/Affine.lean, Faces/Regrade.lean, Faces/ZeroSection.lean}
\leankey{affineFibration}{Faces.bb_face_affineSpace, Faces.AffineSpaceChart, Faces.affineSpaceChart_printed, Faces.exists_cover_mvPolynomial_of_sym, Faces.symBaseChange}{Faces/AffineSpaceFibration.lean}
\leankey{toricSmooth}{cor125, cor125_global}{Sec12PrintedCor.lean, FinalCor125.lean}
\leankey{regular}{IsRegular, SpReg}{Sec13General.lean, Sec13Spaces.lean}
\leankey{dilatation}{IsSpDil, IsSpDil.exists_iso}{Sec13Spaces.lean}
\leankey{dilatationAttractor}{blAction, prop13_1_final, prop13_1_spaces, attractor_isSpDil, exists_spDilatation, attractorOverIso}{Sec13Spaces.lean, Prop131Final.lean, Sec3Functorial.lean}
\leankey{spDilatationExists}{exists_spDilatation}{Sec13Spaces.lean}
\leankey{fact118}{fact118_affine, lemma119, fact118_equivariant_pushout, fact118_overT, fact118_pushout_spaces}{Fact118General.lean, Final3Fact118.lean}
\leankey{dilatationSmooth}{cor13_3', cor13_3_affine, isClosedImmersion_jN}{Cor133NoHyp.lean, FinalCor133.lean}
\leankey{dilatationSmoothAffine}{cor13_3_affine}{FinalCor133.lean}
\leankey{dilatationSmoothSpaces}{cor13_3_free, cor13_3_spaces}{Cor133Spaces.lean, FinalCor133Spaces.lean}
\leankey{ind}{IndPresentation, IsIndAlgebraicSpace, IsIndScheme, IsSeparatedInd}{Sec14Ind.lean, Unrestrict4.lean}
\leankey{indAction}{CatAction, IndActionG, IndActionG.toCatAction}{Sec14Remaining.lean}
\leankey{indAttractor}{prop144_general, prop144_rep_general, prop144_closed_general, prop144_scheme_general, isSeparatedOver_term, prop144_sep, prop144_scheme_sep}{FinalNonQCGeneral.lean, FinalNonQCClauses.lean, Unrestrict4.lean}
\leankey{tangent}{tangent, tangentAction}{Sec16Tangent.lean}
\leankey{tangentBasic}{tangentFibreProductIso, mono_tangentMap, tangentFinalIso}{Sec16Tangent.lean}
\leankey{tangentAttractor}{tangentAttractorIso, tangentAttractorIsoFppf}{Sec16Tangent.lean, Sec16Remaining.lean}
\leankey{lie}{Lie, tangentGrp, lieAction}{Sec16Lie.lean}
\leankey{lieBasic}{lieFibreIso, mono_lieMor, lieTrivialIso}{Sec16Lie.lean}
\leankey{lieAttractor}{lieAttractorIso, lieAttractorLimitIso}{Sec16Lie.lean}
\leankey{semidirect}{attrSemidirect, attrMulIso}{Sec16Lie.lean}
\leankey{rootSmooth}{prop1615_ii_general}{Prop1615General.lean}
\leankey{pureMonotone}{attractorImage, attractorImage_monotone}{Sec15PureMagnets.lean, Sec15General.lean}
\leankey{pure}{IsPure, pureMagnets}{Sec15PureMagnets.lean}
\leankey{fact334}{AttractorOrder, attractorOrder_affine, attractorOrder_of_isAffineHom, fact334_cancel}{Sec15PureMagnets.lean, Sec15General.lean, Thm153Cancel.lean}
\leankey{pureLt}{IsPure.attractorImage_lt, isPure_iff_lt, prop151_affine, prop151_space_final}{Sec15General.lean, Sec15PureMagnets.lean, PrintedProp151.lean, SDR4Final.lean}
\leankey{pureBijection}{pureOf, pureEquivRange, attractorOrder_of_isAffineHom, thm153_cancel'}{Sec15PureMagnets.lean, Sec15General.lean, Thm153Cancel.lean}
\leankey{pureFinite}{thm154_i, thm154_ii, rem156, isLimitPreserving_of_lfp, chartsOfFiniteType_of_lfp, thm154_ii_of_linearizable_qc}{Sec15General.lean, Sec8Limits.lean, Sec15Final.lean}
\leankey{pureReflect}{prop155, prop155_pure, prop155_printed}{Sec15General.lean, Sec15Final.lean}
\leankey{conj157Sep}{conj157_separated, conj157_separated', conj69ii_separated}{Conj157Separated.lean}
\leankey{rank}{monoidRank}{Sec15PureMagnets.lean}
\leankey{rankZn}{monoidRank_bot, monoidRank_closure_singleton, monoidRank_int_pi}{Sec15PureMagnets.lean}
\leankey{pureExamples}{pureMagnets_trivial_action, pureMagnets_eq_of_units, pureMagnets_group_self, pureMagnets_zmod_six, pureMagnets_nat_self, pureMagnets_int_one, pureMagnets_int_two, pureMagnets_p1, pureMagnets_int_prod, pureMagnets_two_three}{Sec15General.lean, Sec15PureMagnets.lean, Sec15Final.lean, SecFinal1316.lean}
\leankey{pureZSix}{magnetIdeal_le_ker_eval, weight_mem_iff_of_magnetIdeal_eq, addSubmonoid_zmod_six, pureMagnets_zmod_six_weights}{Sec15ZModSixWeights.lean}
\leankey{selfPure}{prop31_scheme, prop31EquivScheme, pureMagnets_self_eq, thm46_i_quot_scheme, thm46_ii_scheme, thm46_card_i_scheme, thm46_card_ii_scheme}{Sec15PureMagnets.lean, SelfActionsMain.lean}
\leankey{rooted}{Rooted.trace, Rooted.gen, Rooted.IsStable, Rooted.IsJoinDense, Rooted.isJoinDense_iff}{Rooted/RootedBasic.lean}
\leankey{joinAction}{isPure_sSup_attractorImage, attractorImage_inf_of_le, exists_iter_of_image_le}{JoinAction.lean}
\leankey{joinDenseMinimal}{Conj157.existsUnique_minimal_joinDense}{JoinDenseMinimal.lean}
\leankey{rootedBasic}{Rooted.isJoinDense_iff_eq_span, Rooted.IsJoinDense.pure_of_mem, Rooted.IsJoinDense.pure_gen, Rooted.IsJoinDense.pure_sSup, Rooted.isJoinDense_of_hull}{Rooted/RootedBasic.lean}
\leankey{rootedElements}{Rooted.IsStableElt, Rooted.isJoinDense_mono_iff}{Rooted/RootedBasic.lean}
\leankey{minimalRoot}{Rooted.joinIrreducibles, Rooted.joinIrreducibles_subset, Rooted.isJoinDense_joinIrreducibles, Rooted.existsUnique_minimal_root, Rooted.existsUnique_minimal_root_of_finite}{Rooted/MinimalRoot.lean}
\leankey{rootedPn}{DoubleScalar.PnSite.pnS, DoubleScalar.PnSite.act, DoubleScalar.PnSite.mem_image_iff, DoubleScalar.PnSite.image_eq_iff, DoubleScalar.PnSite.pureMagnets_bijOn, DoubleScalar.PnSite.pn_joinDense, DoubleScalar.pnAction, DoubleScalar.Om, DoubleScalar.isPure_iff, DoubleScalar.pureMagnets_bijOn, DoubleScalar.pn_joinDense}{DoubleScalarLine.lean, DoubleScalarPure.lean, Rooted/DoubleScalarRooted.lean, Rooted/PnSite.lean}
\leankey{rootedP2}{DoubleScalar.PnSite.prop41, DoubleScalar.PnSite.dsAct_joinDense, DoubleScalar.PnSite.ncard_pure, DoubleScalar.dsAct, DoubleScalar.act_formula, DoubleScalar.prop41, DoubleScalar.ncard_pure, DoubleScalar.dsAct_joinDense}{DoubleScalar.lean, Rooted/DoubleScalarRooted.lean, Rooted/PnSite.lean}
\leankey{redSetting}{torusGrading, magnetHopfIdeal, attractorSubgroup, roots, tangentWeights_condition_iff_roots, parabolicDefiningHopfIdeal, leviDefiningHopfIdeal, magnetIdeal_torusGrading_comap_nat, magnetIdeal_torusGrading_comap_bot}{Reductive/TorusGrading.lean, Reductive/AttractorSubgroup.lean, Reductive/TangentRoots.lean, TauCeti/Algebra/AlgebraicGroup/Dynamic/Representability.lean, TauCeti/Algebra/AlgebraicGroup/Dynamic/Levi/Representability.lean, Reductive/PushCocharacter.lean}
\leankey{parabolic}{OpenCell.magnetIdeal_eq_parabolic_of_field}{Reductive/OpenCell/OpenCellFlatness.lean}
\leankey{redEquiv}{reductive_iff_condI_condII_condIII, reductive_of_condI_condII_condIII, condI_condII_condIII_of_reductive}{Reductive/Equivalence.lean}
\leankey{weylGroup}{Weyl.IsNormalizingWith, Weyl.weylGroup, Weyl.mem_weylGroup}{Weyl/Normalizer.lean, Weyl/WeylGroup.lean}
\leankey{weylRoot}{Weyl.joinIrredWeylMagnets, Weyl.weylMagnet, Weyl.fixedQuotient, Weyl.weylMagnet_eq_groupMagnet, Rooted.weylMagnet}{Weyl/Normalizer.lean, Weyl/WeylRootAttractor.lean, Rooted/ReductiveRooted.lean}
\leankey{exRed}{attractor_points_of_no_roots, magnetHopfIdeal_weylMagnet_eq_centralizer, mem_weylMagnet_sPGL2, roots_not_mem_weylMagnet_sPGL2, weylPart_PGL2, cox_eq, cox_cube, mem_weylMagnet_cox, roots_not_mem_weylMagnet_cox, weylPart_PGL3}{Reductive/ExampleRed.lean, Weyl/WeylRootAttractor.lean}
\leankey{weylRootAttractor}{Weyl.orbitHom_normalizing, Weyl.mem_attractorSubgroup_iff_joinIrredWeylMagnets, Weyl.mem_attractorSubgroup_iff_weylMagnet_le, Weyl.groupMagnet, Weyl.magnetHopfIdeal_groupMagnet_eq_centralizer, Weyl.magnetHopfIdeal_weylMagnet_eq_centralizer}{Weyl/Normalizer.lean, Weyl/WeylRootAttractor.lean}
\leankey{rootGroups}{rootGroupIdeal, prescribedUnitHopfIdeal, SplitReductiveData.WithRootSpaces, isRootGroup_of_withRootSpaces, isRootGroup_rootGroupIdeal_of_finrank, exists_generator, bijective_aeval_generator, comul_generator, isHomogeneous_rootGroupIdeal, mem_lieSubalgebra_rootGroupIdeal_iff, limitHom, ker_limitToBot_eq, attractor_semidirect_rootGroup, magnetIdeal_bot_eq_ker, mem_attractorSubgroup_bot_iff, attractor_semidirect_torus, rootGroupIdeal_le_of_isRootGroup, eq_rootGroupIdeal_of_isRootGroup, SplitReductiveData.WithRootSpaces.eq_rootGroupIdeal}{Reductive/RootGroups.lean, Reductive/RootGroupGaLie.lean, Reductive/SplitReductiveRootSpaces.lean, Reductive/GluePrereqs.lean, Reductive/RootGroupUnique.lean}
\leankey{rootGroupDef}{IsRootGroup}{Reductive/RootGroups.lean}
\leankey{weylPart}{weylPart, weylIn_iff_units, weylPart_smul_mem, weylPart_smul_units, weylPart_smul_roots, roots_add_mem, weylIn_of_sub_mem_zmultiples, closure_le_weylPart, closure_normal_in_weylPart, sub_mem_of_mem_closure, WithRootDatum.refl_smul, WithRootDatum.joinIrredWeylMagnets_le_weylPart, WithRootDatum.sub_mem_symLattice, WithRootDatum.joinIrredWeylMagnets_normal}{Reductive/WeylPart.lean, Reductive/WeylPartReflections.lean}
\leankey{redIdentity}{compIdeal, rootComp, compIdeal_eq_of_roots, compIdeal_eq_rootComp, compIdealBar_eq_of_roots, compIdealBar_eq_rootComp}{Reductive/IdentityComponentRoots.lean, Reductive/IdentityComponentBaseChange.lean}
\leankey{redReflects}{mem_of_lieSubalgebra_eq, weylMagnet_le_of_magnetIdeal_eq, attractor_reflects, rep_mem_attractorSubgroup_iff, mul_rep_mem_attractorSubgroup}{Reductive/AttractorReflects.lean}
\leankey{redDescent}{magnetIdeal_eq_iff_baseChange, isPure_baseChange_iff, pureMagnets_baseChange, hull_baseChange, roots_piBar, SplitReductiveData.baseChange}{Reductive/DescentToField.lean, Reductive/ParabolicSubsets.lean}
\leankey{redAttractors}{smooth_quotient_magnetHopfIdeal, compIdeal_eq_rootComp, mem_lieSubalgebra_magnetHopfIdeal_iff, lieSubalgebra_magnetHopfIdeal_eq_of_roots, red_attractors_ii_lattice, componentsMeetNormalizer_lattice, hdet_lattice}{Reductive/AttractorSmooth.lean, Reductive/IdentityComponentRoots.lean, Reductive/LieAttractor.lean, Reductive/RootAttractor.lean, Reductive/ComponentsMeetNormalizerProof.lean}
\leankey{redPure}{Rooted.reductiveHull, red_pure_lattice, reductiveInput_lattice, red_pure_field, red_pure_fieldBar}{Rooted/ReductiveRooted.lean, Reductive/ComponentsMeetNormalizerProof.lean, Reductive/DescentToField.lean}
\leankey{redClassification}{Rooted.pairMagnet, Rooted.IsCompatiblePair, Rooted.bijOn_pairOf, red_classification_i_lattice, red_classification_attractor_lattice, red_classification_i_inj_lattice, finite_pure_lattice, classification_ii_lattice, red_classification_i_field, red_classification_i_inj_field, finite_pure_field, classification_ii_field}{Rooted/ReductiveClassification.lean, Reductive/ComponentsMeetNormalizerProof.lean, Reductive/DescentToField.lean}
\leankey{redMonoPure}{Rooted.pure_mono_iff, red_pure_mono_iff_lattice, Rooted.closure_singleton_inj, red_mono_injOn_lattice, red_mono_bijOn_lattice, red_pure_mono_iff_field, red_mono_bijOn_field}{Rooted/MonogenicPure.lean, Reductive/ComponentsMeetNormalizerProof.lean, Reductive/MonogenicRootBijection.lean, Reductive/DescentToField.lean}
\leankey{redParabolic}{SplitReductiveData.WithRootDatum, Parabolic.inter_roots_eq, red_parabolic_hull_lattice, red_parabolic_lattice, SplitReductiveData.WithRootDatum.parabolic_weyl, Parabolic.weylN_eq_closure_xi, RootSystem.RootDatum.weylN_eq_closure_xi, parabolic_bijOn_lattice, red_parabolic_hull_field, red_parabolic_field, parabolic_bijOn_field, parabolic_std_bijOn_lattice, parabolic_std_bijOn_field}{Reductive/SplitReductiveRootDatum.lean, Reductive/ParabolicBijection.lean, Reductive/WeylParabolic.lean, Reductive/ComponentsMeetNormalizerProof.lean, Reductive/DescentToField.lean}
\leankey{redLevi}{magnetIdeal_closure_symm_eq_comap_bot, magnetIdeal_closure_symm_eq_levi, magnetIdeal_closure_symm_le_ker_iff}{Reductive/LeviAttractor.lean}
\leankey{redRooted}{Rooted.reductiveRoot, Rooted.weylMagnetLower, reductive_joinDense_lattice, joinIrreducibles_isLeast_lattice, joinIrreducibles_eq_lattice, reductive_joinDense_field, joinIrreducibles_isLeast_field, joinIrreducibles_eq_field}{Rooted/ReductiveRooted.lean, Rooted/MinimalRootReductive.lean, Reductive/ComponentsMeetNormalizerProof.lean, Reductive/DescentToField.lean}
\leankey{redSLSp}{SLSp.JoinIrredWeylGenerated, red_SL_Sp_lattice, SLSp.joinIrredWeylGenerated_typeA, SLSp.joinIrredWeylGenerated_typeC, SLSp.joinIrredWeylGenerated_of_image, SLSp.joinIrredWeylGenerated_pi, red_SL_Sp_field}{Reductive/SLSpCriterion.lean, Reductive/ComponentsMeetNormalizerProof.lean, Reductive/DescentToField.lean}
\leankey{homPush}{HomTo, homAtEquivHomTo, HomTo.push, homAtEquivHomTo_pull, prop99_full_tensor}{Sec9Full.lean}
\leankey{comodule}{Coaction, hModuleEquivCoactionFull}{Sec9Full.lean}
\leankey{gradedModule}{Decomp, prop911_full, prop911_full_printed}{Sec9OverB.lean, Sec9Full.lean}
\leankey{regularComp}{regular, regular_comp_ne_zero}{Sec9Hochschild.lean}
\leankey{equivariantLift}{IsThick, Coact, GAct, HochschildH1Vanish, prop101_general}{Prop101General5.lean, Prop101General3.lean, Prop101General4.lean, Prop101General7.lean, Prop101General8.lean}
\leankey{smoothAttractor}{SpacesLift.isFormallySmooth_attractorMap_spaces, SpacesLift.lift103, Sec11Spaces.formallySmoothFor_attractorMap_of_lift}{Sec11SpacesLift.lean, Sec11Spaces.lean, Sec11SpacesCore.lean}
\leankey{smoothAttractorRep}{SpacesLift.isSmooth_repMap_spaces}{Sec11SpacesCor.lean}
\leankey{smoothFace}{SpacesLift.isFormallySmooth_faceMap_spaces}{Sec11SpacesFace.lean}
\leankey{smoothFaceRep}{SpacesLift.isSmooth_repFace_spaces}{Sec11SpacesCor.lean}
\leankey{cochainFormula}{cochainEquivSec, cochainLinearEquiv, delta_phiC}{Sec9OverB.lean, Sec9General.lean}
\leankey{monoidCohomology}{exists_delta_eq, H1_subsingleton, Hn_subsingleton, H1_base, Hn_base, eq_delta_of_cocycle}{Sec9OverB.lean, Sec9General.lean, Sec9Hochschild.lean}
\leankey{invariantLift}{exists_invariant_lift, H1_vanish}{Sec10PrintedMain.lean}
\leankey{equivariantLiftSpaces}{cor103_square}{FinalCor103Square.lean}

\leankey{closedStable}{Torsion.closed_iff_stable, Torsion.closed_iff_stable_rootDatum}{Reductive/TorsionPrimes.lean, Reductive/TorsionPrimesRootDatum.lean}
\leankey{torsionPrimes}{Torsion.torsion_final, Torsion.isTP_iff_candidate, Torsion.torsion_primes, Torsion.torsion_primes_rootDatum, Torsion.torsion_primes_joinIrreducibles}{Torsion/Final.lean, Torsion/Cor34.lean, Reductive/TorsionPrimes.lean, Reductive/TorsionPrimesRootDatum.lean}
\leankey{reflStab}{Torsion.HypR1, Torsion.HypR2, Torsion.A1.hypR1, Torsion.A1.hypR2}{Torsion/Hypotheses.lean, Torsion/A1.lean}
\leankey{thetaLemma}{Torsion.RD.theta_mul, Torsion.RD.theta_of_mem, Torsion.mem_W_of_theta, Torsion.RD.nsmul_theta, Torsion.RD.mem_U_iff}{Torsion/Theta.lean}
\leankey{indexProp}{Torsion.RD.lemma31, Torsion.RD.step33, Torsion.RD.prop33, Torsion.RD.cor34}{Torsion/Lemma31.lean, Torsion/Step33.lean, Torsion/Prop33.lean, Torsion/Cor34.lean}
\leankey{weylLevi}{Torsion.mem_fix_iff_orth, Torsion.fix_eq_hyperplanes, Torsion.fixed_of_Z_le, Torsion.weylMagnetLower_levi, Torsion.index_formula, Torsion.relIndex_eq_det, Torsion.cyclotomic_eval_one, Torsion.cyclotomic_eval_one_not_prime_pow}{Reductive/CyclotomicLevi.lean}

\newcommand{\algmagcommit}{3637c294ac1a}
\DeclareRobustCommand{\stfull}{\textcolor{green!50!black}{$\checkmark$}}
\DeclareRobustCommand{\stplus}{\textcolor{green!50!black}{$\checkmark^{\!+}$}}
\DeclareRobustCommand{\stpart}{\textcolor{orange!90!black}{\tikz[baseline=-0.55ex]{\draw (0,0) circle (0.55ex); \fill (0,0.55ex) arc (90:270:0.55ex) -- cycle;}}}
\DeclareRobustCommand{\stconj}{\textcolor{green!50!black}{$\bigstar$}}
\DeclareRobustCommand{\strepl}{\textcolor{blue!70!black}{$\rightarrow$}}
\DeclareRobustCommand{\stexcl}{\textcolor{black!55}{$\circ$}}
\DeclareRobustCommand{\stnone}{\textcolor{black!55}{--}}

\providecommand{\statuslegend}{\stfull\ formalized; \stplus\ formalized with an added hypothesis;
\stpart\ a special case formalized; \strepl\ corrected and replaced by a statement of this paper;
\stconj\ conjecture proved; \stexcl\ excluded (reductive groups);
\stnone\ no assertion.}

\begin{document}

\title{Algebraic magnetism, formally verified and extended: a survey}
\author{Arnaud Mayeux}
\address{University of Wisconsin--Madison, Madison, WI, USA}
\date{}

\begin{abstract}
We formalize in the Lean proof assistant the three papers on algebraic magnetism published in Semigroup Forum: ``Algebraic magnetism'', ``Algebraic magnetism invariants of self-actions of diagonalizable monoid schemes'' and ``Algebraic magnetism: $X$-products of attractors via $\mathbb{F}_1$-geometry''. We also prove further results and comment on the relation with other work in the literature. The description of this formalization and of these results forms in itself a survey of the theory, from the perspective opened by the formalization. The first main mathematical contribution consists of a Bia\l ynicki-Birula decomposition for faces and of proofs of the existence of strongly fixed-point-reflecting atlases, of the finiteness of pure magnets and of the vanishing of the higher monoid cohomology of $\A(M)$ with quasi-coherent coefficients, three conjectures of ``Algebraic magnetism''. The case of algebraic stacks, where the analogue of the existence of fixed-point-reflecting atlases fails, is treated in joint work with Sylvain Brochard, in which the finiteness of pure magnets was first proved privately, by a comparable method. The second main contribution, for split reductive groups, is a description of the pure magnets and of the join-irreducible ones: besides the classical roots, they are the join-irreducible pure magnets equal to $(1-w)M$ for some $w$ in the Weyl group, which encode the components of attractors; a prime $p$ is a torsion prime of the group, in the sense of Steinberg, if and only if some join-irreducible Weyl magnet $(1-w)M$ has $w$ of $p$-power order.
\end{abstract}

\makeatletter
\def\@settitle{\begin{center}\baselineskip14\p@\relax\bfseries\Large\@title\end{center}}
\makeatother
\maketitle
\setcounter{tocdepth}{1}
\tableofcontents

\section{Introduction}

Let $S$ be a scheme and $M$ a commutative monoid. Algebraic magnetism \cite{AlgMag,SelfActions,XProducts} studies actions of the monoid scheme $\A(M)_S=\Spec(\ZZ[M])\times_{\Spec(\ZZ)}S$ on algebraic spaces $X$ over $S$ through their attractors. If $M$ is an abelian group, $\A(M)_S$ is the underlying monoid scheme of the diagonalizable group scheme $\D(M)_S$ \cite{SGA3}, and an action of $\A(M)_S$ on $X$ is the same as an action of $\D(M)_S$. A \emph{magnet} is a submonoid $N\subseteq M$, and its \emph{attractor} is the fppf sheaf $X^N$ whose $T$-points are the $\A(M)_T$-equivariant morphisms $\A(N)_T\to X_T$. One has $X^M=X$; for $M$ a finitely generated abelian group and $Z\subseteq M$ a subgroup, $X^Z$ is the fixed-point functor of $\D(M/Z)_S$; for $M=\ZZ$ and $N=\NN$, a point of $X^{\NN}$ is a point of $X$ together with an extension of its orbit map to $\mathbb{A}^1$. A face $F$ of $N$ gives morphisms $\iota_{F,N}\colon X^F\to X^N$ and $p_{N,F}\colon X^N\to X^F$ with $p_{N,F}\circ\iota_{F,N}=\id$. A magnet $N$ is \emph{pure} if it is contained in every magnet with the same attractor, and $\mho(a)$ denotes the set of pure magnets of the action $a$.

This text presents algebraic magnetism organized by the mathematics, and every definition and statement outside Section~\ref{sec:reductive} is formalized in a proof assistant. The main statements are the following.

\begin{enumerate}[label=\textup{(\Alph*)},leftmargin=*]
\item (Theorem~\ref{thm:strong-atlas}, Corollary~\ref{cor:representability}.) Let $S$ be separated, $M$ a finitely generated abelian group and $X$ a separated algebraic space over $S$ locally of finite presentation with an action of $\D(M)_S$. Then $X$ has a strongly fixed-point-reflecting atlas, $S$-affine if $|X|$ is quasi-compact, and every attractor $X^N$ is a separated algebraic space locally of finite presentation over $S$. This proves \cite[Conj.~6.9(i)]{AlgMag} for $S$ separated and \cite[Conj.~6.9(ii)]{AlgMag} for $S$ separated and quasi-compact.
\item (Theorem~\ref{thm:BB}, Corollary~\ref{cor:affine-fibration}.) Let $X$ be smooth over $S$ with a Zariski locally linearizable action, $N$ a magnet and $F$ a face of $N$. Then $p_{N,F}\colon X^N\to X^F$ is an affine space fibration: Zariski locally on $X^F$, $X^N\cong X^F\times_S\mathbb{A}^r_S$.
\item (Theorem~\ref{thm:pure}, Corollary~\ref{cor:conj157}.) Let $S$ be separated and quasi-compact, $M$ a finitely generated abelian group and $X$ a separated algebraic space over $S$ of finite presentation. Then $\mho(a)$ is finite, and $N\mapsto X^N$ is a bijection from $\mho(a)$ onto the set of attractors. This proves \cite[Conj.~15.7]{AlgMag} in this case.
\item (Section~\ref{sec:reductive}.) Let $G$ be a split connected reductive group over a field $k$, with split maximal torus $T=\D(M)_k$ acting by conjugation, roots $\Phi\subseteq M$ and Weyl group $W$. Besides the classical roots $\NN\alpha$, $\alpha\in\Phi$, the pure magnets are built from the lattices $(1-w)M=\{m-w(m)\mid m\in M\}$, $w\in W$, the characters vanishing on the fixed points $T^w$; we call them \emph{Weyl magnets}. Magnetism recovers the classical theory: $\alpha\mapsto\NN\alpha$ is a bijection from $\Phi$ onto the set of nonzero monogenerated pure magnets, and the Borel, parabolic and Levi subgroups containing $T$ are attractors of explicit pure magnets (Table~\ref{tab:intro}). It also adds to it: the Weyl magnets refine the Levi subgroups containing $T$ by the torsion of $T^w$, and a prime $p$ is a torsion prime of $G$, in the sense of Steinberg \cite{Steinberg75}, if and only if some join-irreducible Weyl magnet (see below) $(1-w)M$ has $w$ of $p$-power order (Theorem~\ref{thm:torsion-primes}). Precisely, for a magnet $N$ put $R=N\cap\Phi$, $W_N=\{w\in W\mid(1-w)M\subseteq N\}$, $K_R$ the connected component of $G^{\mon R}$ containing the identity point, and $n_w\in N_G(T)(k)$ representing $w$; then, for $k$ algebraically closed, $G^N(k)=\bigcup_{w\in W_N}K_R(k)\,n_w$, and for every $k$, $N\in\mho(a)$ if and only if $N=\mon{R}+\sum_{w\in W_N}(1-w)M$.
\end{enumerate}

A diagonalizable action $a$ is a \emph{join action} if sums of pure magnets are pure; then the pure magnets form a join-semilattice $\mho(a)$, the join being the sum, and we use the language of lattices and semigroups, see e.g.\ \cite{ErneWilke83}. A family $\mathcal F$ of magnets is \emph{join-dense} for $a$ if the pure magnets are exactly the sums of members of $\mathcal F$, that is, if $\mathcal F$ is a join-dense subset of $\mho(a)$ in the sense of lattice theory \cite[p.~354]{ErneWilke83}; we then say that the action $a$ is \emph{join-dense by} $\mathcal F$; the smallest such family, when it exists, is the set $\underline{\mho}(a)$ of \emph{join-irreducible} pure magnets, those that are not the sum of strictly smaller pure magnets. Every action of $\D(M)_S$, $M$ a finitely generated abelian group, on a separated algebraic space of finite presentation over a separated quasi-compact $S$ is a join action (Proposition~\ref{prop:join-action}), and it is \emph{join-irreducible}: join-dense by the family of its join-irreducible pure magnets, which is the unique minimal join-dense family (Corollary~\ref{cor:join-dense-minimal}).

The term \emph{Weyl magnet} is introduced in the present text. It is specific to the conjugation action of a split maximal torus on a split reductive group; it is not a notion of algebraic magnetism for general actions. The Weyl magnet of $w\in W$ is $(1-w)M$, and it does not determine $w$: one has $(1-w^{-1})M=(1-w)M$, since $(1-w^{-1})m=-(1-w)w^{-1}(m)$. So $G$ has at most $|W|$ Weyl magnets, and fewer as soon as some $w\neq w^{-1}$. The attractor of the Weyl magnet $(1-w)M$ is the centralizer $Z_G(T^w)$. Here a connected component of a locally noetherian scheme means a connected component of its underlying topological space, with its structure of open subscheme. Roots and Weyl magnets see different connected components of an attractor $G^N$: by (D), for $k$ algebraically closed, the roots in $N$ give the connected component $K_R$ containing the identity point, and the Weyl magnets in $N$ give the others, the cosets $K_R\,n_w$, $w\in W_N$; for every $k$, $n_w\in G^N$ if and only if $w\in W_N$. For $G=\mathrm{PGL}_2$, with $M=\ZZ\alpha$ and $s$ the reflection, the Weyl magnet $(1-s)M=2\ZZ\alpha$ contains no root. Its attractor is the normalizer $N_G(T)$, whose two connected components are $T$ and $n_sT$, and $2\ZZ\alpha$ is a join-irreducible pure magnet. For the conjugation action of $T$ on $G$, the classical roots $\NN\alpha$ together with the join-irreducible Weyl magnets form the smallest join-dense family. For $\mathrm{GL}_n$, $\mathrm{SL}_n$ and $\mathrm{Sp}_{2n}$, by contrast with $\mathrm{PGL}_2$, every Weyl magnet is generated by the roots it contains (Theorem~\ref{thm:red-SL-Sp}); there the roots in $N$ determine the attractor, $G^N=G^{\mon R}$, and the magnets $\NN\alpha$ alone are join-dense. Table~\ref{tab:intro} lists, for the split reductive groups of rank at most $10$ and the exceptional groups, the classical objects next to the magnetic objects recovering them, and $\underline{\mho}(a)$.

{\tiny\setlength{\tabcolsep}{1.1pt}
\begin{longtable}{ll|rr|rr|rr|rr|cc|rr}
\caption{Split reductive groups $G$: classical objects and the magnetic objects recovering them (equal columns), and what magnetism adds. A reductive group has the row of its derived group. The number of Levi subgroups of $E_8$ containing $T$ is not computed. See Section~\ref{sec:torsion}.}\label{tab:intro}\\
& & \multicolumn{10}{c|}{magnetism recovers} & \multicolumn{2}{c}{magnetism adds}\\
type & group & \rotatebox{90}{roots} & \rotatebox{90}{pure magnets $\NN\alpha$} & \rotatebox{90}{Borel subgroups $\supseteq T$} & \rotatebox{90}{magnets generated by a basis} & \rotatebox{90}{parabolic subgroups $\supseteq T$} & \rotatebox{90}{pure magnets containing a basis} & \rotatebox{90}{Levi subgroups $\supseteq T$} & \rotatebox{90}{group magnets $\ZZ\zeta$} & \rotatebox{90}{torsion primes} & \rotatebox{90}{primes of join-irr.\ Weyl magnets} & \rotatebox{90}{size of $\underline{\mho}(a)$} & \rotatebox{90}{join-irr.\ Weyl magnets}\\ \hline \endfirsthead
type & group & \multicolumn{12}{c}{(continued)}\\ \hline \endhead
$A_{1}$ & $SL_{2}$ & 2 & 2 & 2 & 2 & 3 & 3 & 2 & 2 & none & none & 2 & 0\\
$A_{1}$ & $PGL_{2}$ & 2 & 2 & 2 & 2 & 3 & 3 & 2 & 2 & 2 & 2 & 3 & 1\\
$A_{2}$ & $SL_{3}$ & 6 & 6 & 6 & 6 & 13 & 13 & 5 & 5 & none & none & 6 & 0\\
$A_{2}$ & $PGL_{3}$ & 6 & 6 & 6 & 6 & 13 & 13 & 5 & 5 & 3 & 3 & 7 & 1\\
$A_{3}$ & $SL_{4}$ & 12 & 12 & 24 & 24 & 75 & 75 & 15 & 15 & none & none & 12 & 0\\
$A_{3}$ & $SL_{4}/\mu_{2}$ & 12 & 12 & 24 & 24 & 75 & 75 & 15 & 15 & 2 & 2 & 15 & 3\\
$A_{3}$ & $PGL_{4}$ & 12 & 12 & 24 & 24 & 75 & 75 & 15 & 15 & 2 & 2 & 18 & 6\\
$A_{4}$ & $SL_{5}$ & 20 & 20 & 120 & 120 & 541 & 541 & 52 & 52 & none & none & 20 & 0\\
$A_{4}$ & $PGL_{5}$ & 20 & 20 & 120 & 120 & 541 & 541 & 52 & 52 & 5 & 5 & 26 & 6\\
$A_{5}$ & $SL_{6}$ & 30 & 30 & 720 & 720 & 4683 & 4683 & 203 & 203 & none & none & 30 & 0\\
$A_{5}$ & $SL_{6}/\mu_{2}$ & 30 & 30 & 720 & 720 & 4683 & 4683 & 203 & 203 & 2 & 2 & 45 & 15\\
$A_{5}$ & $SL_{6}/\mu_{3}$ & 30 & 30 & 720 & 720 & 4683 & 4683 & 203 & 203 & 3 & 3 & 50 & 20\\
$A_{5}$ & $PGL_{6}$ & 30 & 30 & 720 & 720 & 4683 & 4683 & 203 & 203 & 2, 3 & 2, 3 & 65 & 35\\
$A_{6}$ & $SL_{7}$ & 42 & 42 & 5040 & 5040 & 47293 & 47293 & 877 & 877 & none & none & 42 & 0\\
$A_{6}$ & $PGL_{7}$ & 42 & 42 & 5040 & 5040 & 47293 & 47293 & 877 & 877 & 7 & 7 & 162 & 120\\
$A_{7}$ & $SL_{8}$ & 56 & 56 & 40320 & 40320 & 545835 & 545835 & 4140 & 4140 & none & none & 56 & 0\\
$A_{7}$ & $SL_{8}/\mu_{2}$ & 56 & 56 & 40320 & 40320 & 545835 & 545835 & 4140 & 4140 & 2 & 2 & 161 & 105\\
$A_{7}$ & $SL_{8}/\mu_{4}$ & 56 & 56 & 40320 & 40320 & 545835 & 545835 & 4140 & 4140 & 2 & 2 & 791 & 735\\
$A_{7}$ & $PGL_{8}$ & 56 & 56 & 40320 & 40320 & 545835 & 545835 & 4140 & 4140 & 2 & 2 & 2051 & 1995\\
$A_{8}$ & $SL_{9}$ & 72 & 72 & 362880 & 362880 & 7087261 & 7087261 & 21147 & 21147 & none & none & 72 & 0\\
$A_{8}$ & $SL_{9}/\mu_{3}$ & 72 & 72 & 362880 & 362880 & 7087261 & 7087261 & 21147 & 21147 & 3 & 3 & 1192 & 1120\\
$A_{8}$ & $PGL_{9}$ & 72 & 72 & 362880 & 362880 & 7087261 & 7087261 & 21147 & 21147 & 3 & 3 & 7912 & 7840\\
$A_{9}$ & $SL_{10}$ & 90 & 90 & 3628800 & 3628800 & 102247563 & 102247563 & 115975 & 115975 & none & none & 90 & 0\\
$A_{9}$ & $SL_{10}/\mu_{2}$ & 90 & 90 & 3628800 & 3628800 & 102247563 & 102247563 & 115975 & 115975 & 2 & 2 & 1035 & 945\\
$A_{9}$ & $SL_{10}/\mu_{5}$ & 90 & 90 & 3628800 & 3628800 & 102247563 & 102247563 & 115975 & 115975 & 5 & 5 & 18234 & 18144\\
$A_{9}$ & $PGL_{10}$ & 90 & 90 & 3628800 & 3628800 & 102247563 & 102247563 & 115975 & 115975 & 2, 5 & 2, 5 & 19179 & 19089\\
$A_{10}$ & $SL_{11}$ & 110 & 110 & 39916800 & 39916800 & 1622632573 & 1622632573 & 678570 & 678570 & none & none & 110 & 0\\
$A_{10}$ & $PGL_{11}$ & 110 & 110 & 39916800 & 39916800 & 1622632573 & 1622632573 & 678570 & 678570 & 11 & 11 & 362990 & 362880\\
\hline
$B_{1}$ & $SO_{3}$ & 2 & 2 & 2 & 2 & 3 & 3 & 2 & 2 & 2 & 2 & 3 & 1\\
$B_{1}$ & $Spin_{3}$ & 2 & 2 & 2 & 2 & 3 & 3 & 2 & 2 & none & none & 2 & 0\\
$B_{2}$ & $SO_{5}$ & 8 & 8 & 8 & 8 & 17 & 17 & 6 & 6 & 2 & 2 & 10 & 2\\
$B_{2}$ & $Spin_{5}$ & 8 & 8 & 8 & 8 & 17 & 17 & 6 & 6 & none & none & 8 & 0\\
$B_{3}$ & $SO_{7}$ & 18 & 18 & 48 & 48 & 147 & 147 & 24 & 24 & 2 & 2 & 21 & 3\\
$B_{3}$ & $Spin_{7}$ & 18 & 18 & 48 & 48 & 147 & 147 & 24 & 24 & 2 & 2 & 19 & 1\\
$B_{4}$ & $SO_{9}$ & 32 & 32 & 384 & 384 & 1697 & 1697 & 116 & 116 & 2 & 2 & 36 & 4\\
$B_{4}$ & $Spin_{9}$ & 32 & 32 & 384 & 384 & 1697 & 1697 & 116 & 116 & 2 & 2 & 37 & 5\\
$B_{5}$ & $SO_{11}$ & 50 & 50 & 3840 & 3840 & 24483 & 24483 & 648 & 648 & 2 & 2 & 55 & 5\\
$B_{5}$ & $Spin_{11}$ & 50 & 50 & 3840 & 3840 & 24483 & 24483 & 648 & 648 & 2 & 2 & 66 & 16\\
$B_{6}$ & $SO_{13}$ & 72 & 72 & 46080 & 46080 & 423857 & 423857 & 4088 & 4088 & 2 & 2 & 78 & 6\\
$B_{6}$ & $Spin_{13}$ & 72 & 72 & 46080 & 46080 & 423857 & 423857 & 4088 & 4088 & 2 & 2 & 114 & 42\\
$B_{7}$ & $SO_{15}$ & 98 & 98 & 645120 & 645120 & 8560947 & 8560947 & 28640 & 28640 & 2 & 2 & 105 & 7\\
$B_{7}$ & $Spin_{15}$ & 98 & 98 & 645120 & 645120 & 8560947 & 8560947 & 28640 & 28640 & 2 & 2 & 197 & 99\\
$B_{8}$ & $SO_{17}$ & 128 & 128 & 10321920 & 10321920 & 197613377 & 197613377 & 219920 & 219920 & 2 & 2 & 136 & 8\\
$B_{8}$ & $Spin_{17}$ & 128 & 128 & 10321920 & 10321920 & 197613377 & 197613377 & 219920 & 219920 & 2 & 2 & 347 & 219\\
$B_{9}$ & $SO_{19}$ & 162 & 162 & 185794560 & 185794560 & 5131725123 & 5131725123 & 1832224 & 1832224 & 2 & 2 & 171 & 9\\
$B_{9}$ & $Spin_{19}$ & 162 & 162 & 185794560 & 185794560 & 5131725123 & 5131725123 & 1832224 & 1832224 & 2 & 2 & 628 & 466\\
$B_{10}$ & $SO_{21}$ & 200 & 200 & 3715891200 & 3715891200 & 148070287697 & 148070287697 & 16430176 & 16430176 & 2 & 2 & 210 & 10\\
$B_{10}$ & $Spin_{21}$ & 200 & 200 & 3715891200 & 3715891200 & 148070287697 & 148070287697 & 16430176 & 16430176 & 2 & 2 & 1168 & 968\\
\hline
$C_{1}$ & $Sp_{2}$ & 2 & 2 & 2 & 2 & 3 & 3 & 2 & 2 & none & none & 2 & 0\\
$C_{1}$ & $PSp_{2}$ & 2 & 2 & 2 & 2 & 3 & 3 & 2 & 2 & 2 & 2 & 3 & 1\\
$C_{2}$ & $Sp_{4}$ & 8 & 8 & 8 & 8 & 17 & 17 & 6 & 6 & none & none & 8 & 0\\
$C_{2}$ & $PSp_{4}$ & 8 & 8 & 8 & 8 & 17 & 17 & 6 & 6 & 2 & 2 & 10 & 2\\
$C_{3}$ & $Sp_{6}$ & 18 & 18 & 48 & 48 & 147 & 147 & 24 & 24 & none & none & 18 & 0\\
$C_{3}$ & $PSp_{6}$ & 18 & 18 & 48 & 48 & 147 & 147 & 24 & 24 & 2 & 2 & 25 & 7\\
$C_{4}$ & $Sp_{8}$ & 32 & 32 & 384 & 384 & 1697 & 1697 & 116 & 116 & none & none & 32 & 0\\
$C_{4}$ & $PSp_{8}$ & 32 & 32 & 384 & 384 & 1697 & 1697 & 116 & 116 & 2 & 2 & 57 & 25\\
$C_{5}$ & $Sp_{10}$ & 50 & 50 & 3840 & 3840 & 24483 & 24483 & 648 & 648 & none & none & 50 & 0\\
$C_{5}$ & $PSp_{10}$ & 50 & 50 & 3840 & 3840 & 24483 & 24483 & 648 & 648 & 2 & 2 & 131 & 81\\
$C_{6}$ & $Sp_{12}$ & 72 & 72 & 46080 & 46080 & 423857 & 423857 & 4088 & 4088 & none & none & 72 & 0\\
$C_{6}$ & $PSp_{12}$ & 72 & 72 & 46080 & 46080 & 423857 & 423857 & 4088 & 4088 & 2 & 2 & 403 & 331\\
$C_{7}$ & $Sp_{14}$ & 98 & 98 & 645120 & 645120 & 8560947 & 8560947 & 28640 & 28640 & none & none & 98 & 0\\
$C_{7}$ & $PSp_{14}$ & 98 & 98 & 645120 & 645120 & 8560947 & 8560947 & 28640 & 28640 & 2 & 2 & 1401 & 1303\\
$C_{8}$ & $Sp_{16}$ & 128 & 128 & 10321920 & 10321920 & 197613377 & 197613377 & 219920 & 219920 & none & none & 128 & 0\\
$C_{8}$ & $PSp_{16}$ & 128 & 128 & 10321920 & 10321920 & 197613377 & 197613377 & 219920 & 219920 & 2 & 2 & 6065 & 5937\\
$C_{9}$ & $Sp_{18}$ & 162 & 162 & 185794560 & 185794560 & 5131725123 & 5131725123 & 1832224 & 1832224 & none & none & 162 & 0\\
$C_{9}$ & $PSp_{18}$ & 162 & 162 & 185794560 & 185794560 & 5131725123 & 5131725123 & 1832224 & 1832224 & 2 & 2 & 26947 & 26785\\
$C_{10}$ & $Sp_{20}$ & 200 & 200 & 3715891200 & 3715891200 & 148070287697 & 148070287697 & 16430176 & 16430176 & none & none & 200 & 0\\
$C_{10}$ & $PSp_{20}$ & 200 & 200 & 3715891200 & 3715891200 & 148070287697 & 148070287697 & 16430176 & 16430176 & 2 & 2 & 133851 & 133651\\
\hline
$D_{2}$ & $SO_{4}$ & 4 & 4 & 4 & 4 & 9 & 9 & 4 & 4 & 2 & 2 & 5 & 1\\
$D_{2}$ & $Spin_{4}$ & 4 & 4 & 4 & 4 & 9 & 9 & 4 & 4 & none & none & 4 & 0\\
$D_{2}$ & $PSO_{4}$ & 4 & 4 & 4 & 4 & 9 & 9 & 4 & 4 & 2 & 2 & 6 & 2\\
$D_{2}$ & $HSpin_{4}$ & 4 & 4 & 4 & 4 & 9 & 9 & 4 & 4 & 2 & 2 & 5 & 1\\
$D_{3}$ & $SO_{6}$ & 12 & 12 & 24 & 24 & 75 & 75 & 15 & 15 & 2 & 2 & 15 & 3\\
$D_{3}$ & $Spin_{6}$ & 12 & 12 & 24 & 24 & 75 & 75 & 15 & 15 & none & none & 12 & 0\\
$D_{3}$ & $PSO_{6}$ & 12 & 12 & 24 & 24 & 75 & 75 & 15 & 15 & 2 & 2 & 18 & 6\\
$D_{4}$ & $SO_{8}$ & 24 & 24 & 192 & 192 & 865 & 865 & 72 & 72 & 2 & 2 & 30 & 6\\
$D_{4}$ & $Spin_{8}$ & 24 & 24 & 192 & 192 & 865 & 865 & 72 & 72 & 2 & 2 & 25 & 1\\
$D_{4}$ & $PSO_{8}$ & 24 & 24 & 192 & 192 & 865 & 865 & 72 & 72 & 2 & 2 & 55 & 31\\
$D_{4}$ & $HSpin_{8}$ & 24 & 24 & 192 & 192 & 865 & 865 & 72 & 72 & 2 & 2 & 30 & 6\\
$D_{5}$ & $SO_{10}$ & 40 & 40 & 1920 & 1920 & 12483 & 12483 & 403 & 403 & 2 & 2 & 50 & 10\\
$D_{5}$ & $Spin_{10}$ & 40 & 40 & 1920 & 1920 & 12483 & 12483 & 403 & 403 & 2 & 2 & 45 & 5\\
$D_{5}$ & $PSO_{10}$ & 40 & 40 & 1920 & 1920 & 12483 & 12483 & 403 & 403 & 2 & 2 & 120 & 80\\
$D_{6}$ & $SO_{12}$ & 60 & 60 & 23040 & 23040 & 216113 & 216113 & 2546 & 2546 & 2 & 2 & 75 & 15\\
$D_{6}$ & $Spin_{12}$ & 60 & 60 & 23040 & 23040 & 216113 & 216113 & 2546 & 2546 & 2 & 2 & 76 & 16\\
$D_{6}$ & $PSO_{12}$ & 60 & 60 & 23040 & 23040 & 216113 & 216113 & 2546 & 2546 & 2 & 2 & 406 & 346\\
$D_{6}$ & $HSpin_{12}$ & 60 & 60 & 23040 & 23040 & 216113 & 216113 & 2546 & 2546 & 2 & 2 & 166 & 106\\
$D_{7}$ & $SO_{14}$ & 84 & 84 & 322560 & 322560 & 4364979 & 4364979 & 17867 & 17867 & 2 & 2 & 105 & 21\\
$D_{7}$ & $Spin_{14}$ & 84 & 84 & 322560 & 322560 & 4364979 & 4364979 & 17867 & 17867 & 2 & 2 & 126 & 42\\
$D_{7}$ & $PSO_{14}$ & 84 & 84 & 322560 & 322560 & 4364979 & 4364979 & 17867 & 17867 & 2 & 2 & 1806 & 1722\\
$D_{8}$ & $SO_{16}$ & 112 & 112 & 5160960 & 5160960 & 100757313 & 100757313 & 137528 & 137528 & 2 & 2 & 140 & 28\\
$D_{8}$ & $Spin_{16}$ & 112 & 112 & 5160960 & 5160960 & 100757313 & 100757313 & 137528 & 137528 & 2 & 2 & 211 & 99\\
$D_{8}$ & $PSO_{16}$ & 112 & 112 & 5160960 & 5160960 & 100757313 & 100757313 & 137528 & 137528 & 2 & 2 & 6077 & 5965\\
$D_{8}$ & $HSpin_{16}$ & 112 & 112 & 5160960 & 5160960 & 100757313 & 100757313 & 137528 & 137528 & 2 & 2 & 7407 & 7295\\
$D_{9}$ & $SO_{18}$ & 144 & 144 & 92897280 & 92897280 & 2616517443 & 2616517443 & 1149079 & 1149079 & 2 & 2 & 180 & 36\\
$D_{9}$ & $Spin_{18}$ & 144 & 144 & 92897280 & 92897280 & 2616517443 & 2616517443 & 1149079 & 1149079 & 2 & 2 & 363 & 219\\
$D_{9}$ & $PSO_{18}$ & 144 & 144 & 92897280 & 92897280 & 2616517443 & 2616517443 & 1149079 & 1149079 & 2 & 2 & 47088 & 46944\\
$D_{10}$ & $SO_{20}$ & 180 & 180 & 1857945600 & 1857945600 & 75496735057 & 75496735057 & 10335766 & 10335766 & 2 & 2 & 225 & 45\\
$D_{10}$ & $Spin_{20}$ & 180 & 180 & 1857945600 & 1857945600 & 75496735057 & 75496735057 & 10335766 & 10335766 & 2 & 2 & 646 & 466\\
$D_{10}$ & $PSO_{20}$ & 180 & 180 & 1857945600 & 1857945600 & 75496735057 & 75496735057 & 10335766 & 10335766 & 2 & 2 & 133876 & 133696\\
$D_{10}$ & $HSpin_{20}$ & 180 & 180 & 1857945600 & 1857945600 & 75496735057 & 75496735057 & 10335766 & 10335766 & 2 & 2 & 318886 & 318706\\
\hline
$G_{2}$ & $G_{2}$ & 12 & 12 & 12 & 12 & 25 & 25 & 8 & 8 & 2 & 2 & 13 & 1\\
\hline
$F_{4}$ & $F_{4}$ & 48 & 48 & 1152 & 1152 & 5089 & 5089 & 268 & 268 & 2, 3 & 2, 3 & 69 & 21\\
\hline
$E_{6}$ & $E_{6}^{\rm sc}$ & 72 & 72 & 51840 & 51840 & 486397 & 486397 & 4598 & 4598 & 2, 3 & 2, 3 & 157 & 85\\
$E_{6}$ & $E_{6}^{\rm ad}$ & 72 & 72 & 51840 & 51840 & 486397 & 486397 & 4598 & 4598 & 2, 3 & 2, 3 & 397 & 325\\
$E_{7}$ & $E_{7}^{\rm sc}$ & 126 & 126 & 2903040 & 2903040 & 39296883 & 39296883 & 90408 & 90408 & 2, 3 & 2, 3 & 3515 & 3389\\
$E_{7}$ & $E_{7}^{\rm ad}$ & 126 & 126 & 2903040 & 2903040 & 39296883 & 39296883 & 90408 & 90408 & 2, 3 & 2, 3 & 4775 & 4649\\
$E_{8}$ & $E_{8}$ & 240 & 240 & 696729600 & 696729600 & 13606376641 & 13606376641 & -- & -- & 2, 3, 5 & 2, 3, 5 & 590335 & 590095\\
\end{longtable}}

\subsection*{Foundations}
Algebraic spaces over $S$ are taken in the sense of the Stacks project \cite{Stacks}: sheaves on the big fppf site of $S$ with representable diagonal and an \'etale atlas by a scheme. Attractors are fppf sheaves on the same site. Schemes are taken in the modern sense, without separatedness condition. In Grothendieck's definition \cite{EGAI}, a scheme is separated and the general notion is a prescheme; the separatedness hypotheses in (A) and (C) are those of the original definition, and the hypothesis on $S$ in (A) cannot be omitted (Remark~\ref{rem:conj69}). Diagonalizable group schemes, maximal tori and split reductive groups are those of \cite{SGA3}, and classical results on reductive groups, such as the decomposition of the Lie algebra into root spaces and the conjugacy of maximal tori, are used with a citation to \cite{SGA3}.

\subsection*{Formalization in Lean}
Every definition and statement of this text outside Section~\ref{sec:reductive} is formalized in the proof assistant Lean~4 \cite{deMoura}; the code is publicly available \cite{AlgMagLean}. All statements of the three Semigroup Forum papers are formalized, except those on reductive groups, which are formalized only in part; a few statements are formalized in special cases, as listed in Section~\ref{sec:formalization}. In this version, the results on reductive groups are not presented as formalized. The reason is that the Lean libraries we build on, Tau Ceti and Mathlib, do not contain, as of 10 November 2026, the full structure theory of reductive groups; we make public the fragments of the code that give a partial formalization. We stress once more that Section~\ref{sec:reductive} is to be considered as not formalized, even though the partial Lean declarations are given there for information. The numerical entries of Table~\ref{tab:intro} for groups of rank at least $4$, and the counts of Remark~\ref{rem:torsion-atoms}, are computed and not formalized; the rows of rank at most $3$ are checked by the Lean kernel. Lean is based on dependent type theory: definitions and proofs are written in a formal language, and each proof is checked by a small kernel. Each definition and statement in the text is followed by the names of the Lean declarations formalizing it, with the same hypotheses. Section~\ref{sec:formalization} gives the pinned commit, the formalization status of every numbered statement of \cite{AlgMag}, and the declarations for the statements of this text. The library consists of $656$ modules, with $193193$ lines of Lean, containing $10864$ declarations, among them $7449$ theorems, and each depends only on the axioms \texttt{propext}, \texttt{Classical.choice} and \texttt{Quot.sound}. These counts do not include the $232$ modules of Tau Ceti imported by the library, with $3248$ declarations (see the paragraph on libraries below).

\subsection*{Libraries}
The formalization is built on Mathlib \cite{Mathlib}, for commutative algebra, linear algebra, Hopf algebras, K\"ahler differentials, root systems and topology, and on Stacks, the Lean library of AxiomMath \cite{AxiomStacks}, for the big fppf site, algebraic spaces and their morphisms following \cite{Stacks}. The library also contains copies of files of Tau Ceti \cite{TauCeti} on affine group schemes, identity components, Lie algebras, cocharacters, and the dynamic description of parabolic and Levi subgroups. Tau Ceti is developed against a more recent version of Mathlib than the one used by AxiomMath Stacks; to combine the two libraries, these files were adapted to the version of Mathlib used here, with changes limited to what this requires, and their mathematical content is that of Tau Ceti. We are grateful to the authors of Tau Ceti for making this work available.

\begin{figure}[H]
\centering
\includegraphics[width=\textwidth,height=0.42\textheight,keepaspectratio]{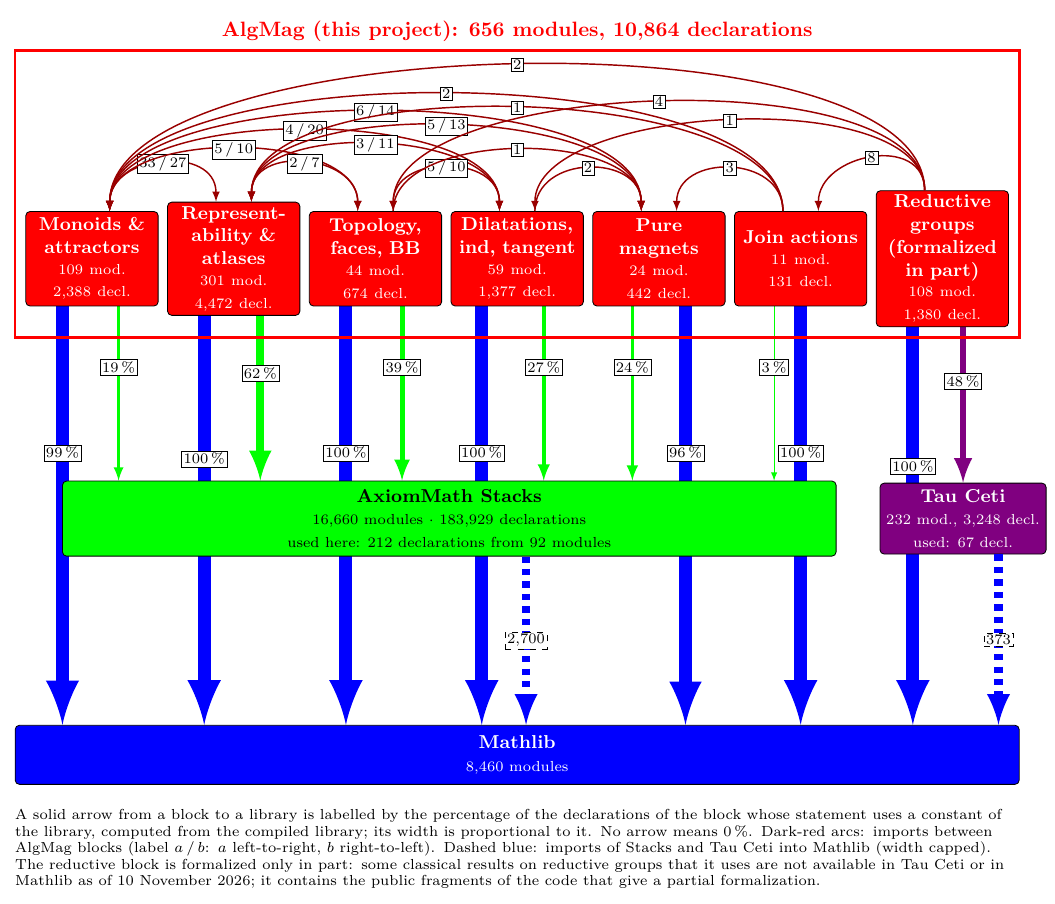}
\caption{Layered overview of the formalization. Red: the main blocks of this library, with their numbers of modules and declarations. Green: AxiomMath Stacks; purple: Tau Ceti; blue: Mathlib. Arcs between blocks show imports between blocks, labelled $a/b$ (import statements from left to right and from right to left). An arrow from a block to a library is labelled by the percentage of the declarations of the block whose statement uses a constant of the library, computed from the compiled library, with width proportional to it; no arrow means $0\,\%$. The boxes of AxiomMath Stacks and Tau Ceti give the number of their declarations used in the statements or proofs of this library. The reductive block is formalized only in part: some classical results on reductive groups that it uses are not available in Tau Ceti or in Mathlib as of 10 November 2026; it contains the public fragments of the code that give a partial formalization.}
\label{fig:layers}
\end{figure}

\begin{figure}[H]
\centering
\includegraphics[width=\textwidth,height=0.36\textheight,keepaspectratio]{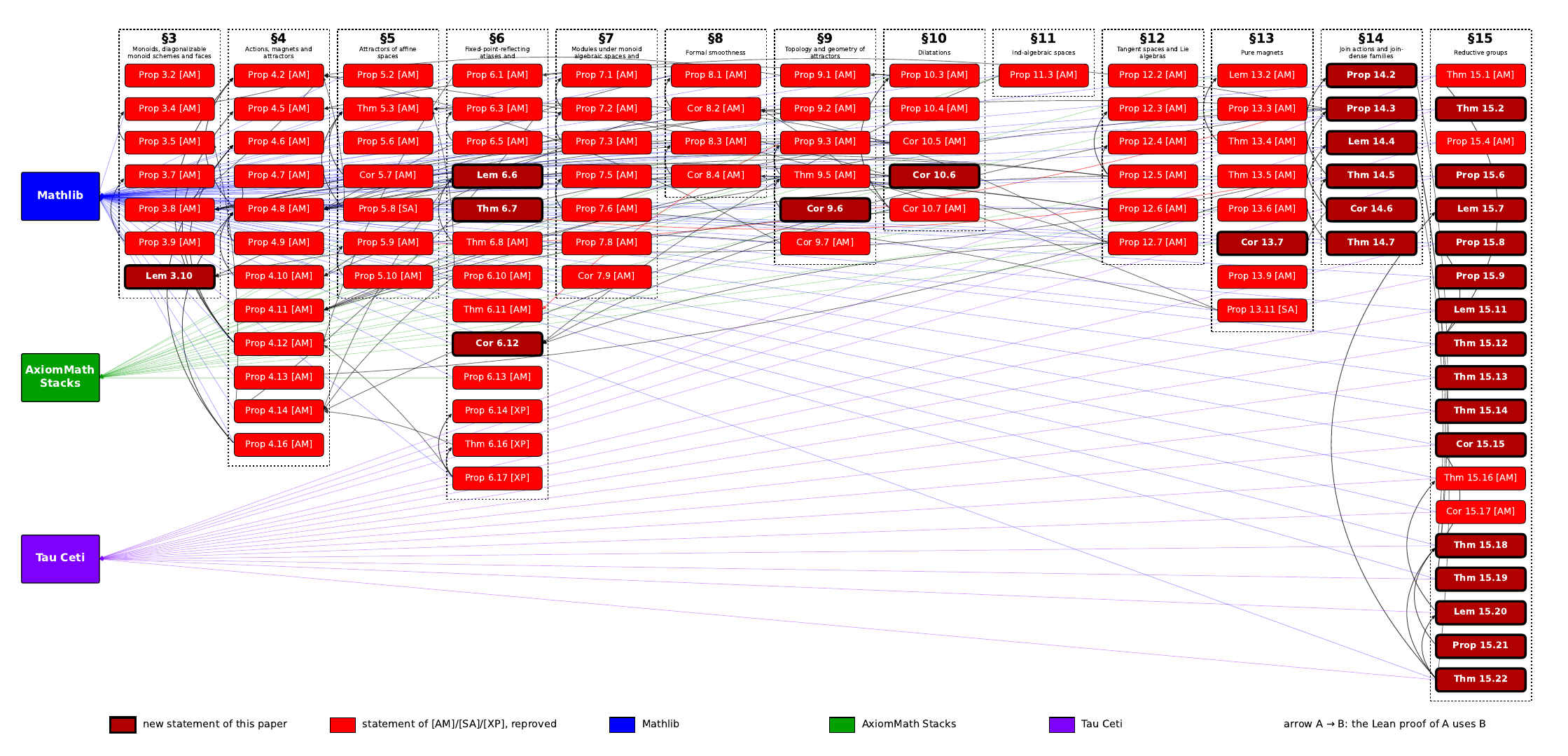}
\caption{Dependency graph of the numbered statements, extracted from the Lean library. An arrow $A\to B$ means that the formal proof of $A$ uses $B$, directly or through declarations that are not statements of the text; arrows implied by others are omitted. Columns are sections. Dark red boxes with thick borders are statements with no counterpart in \cite{AlgMag,SelfActions,XProducts}; red boxes are statements of these papers. Blue, green and purple boxes stand for Mathlib, AxiomMath Stacks and Tau Ceti; a statement is joined to a library when its proof uses that library substantially.}
\label{fig:statements}
\end{figure}

Figures~\ref{fig:layers} and~\ref{fig:statements} show the structure of the formalization at two scales. Figure~\ref{fig:layers} shows the main blocks of the library and their use of the three underlying libraries: for each block and each library, the percentage of the declarations of the block whose statement uses a constant of the library. Figure~\ref{fig:statements} shows the dependencies between the statements of this text.

\subsection*{Algebraic magnetism for stacks}
The natural analogue of Theorem~\ref{thm:strong-atlas} fails for stacks, as was discovered in the joint project with Sylvain Brochard \cite{BrochardMayeux}. Nevertheless, the finiteness of pure magnets still holds for stacks in a natural form, and Theorem~\ref{thm:pure-finite}(ii) was first proved as part of the project \cite{BrochardMayeux}, by a different but comparable method. We refer to \cite{BrochardMayeux} for algebraic magnetism for stacks and for more details on the differences between the cases of algebraic spaces and of stacks. As in the Stacks project \cite{Stacks}, where algebraic spaces and algebraic stacks are developed in separate parts, the two cases of algebraic magnetism are developed in separate texts, because the theory behaves differently in the two cases.

Appendix~A of \cite{AlgMag}, by M.~Romagny after Drinfeld \cite{Drinfeld} and Alper, Hall and Rydh \cite{AHR}, is not discussed in this text. Theorem~6.8 of \cite{AlgMag} is proved here for $S$ and $X$ separated (Theorem~\ref{thm:atlas}), as a corollary of the existence of strongly fixed-point-reflecting atlases (Theorem~\ref{thm:strong-atlas}). The proof uses the existence of equivariant \'etale affine neighbourhoods for diagonalizable groups, a case of the local structure theorem of \cite{AHR}, which is proved here following the architecture of their proof, with algebraic stacks, Tannaka duality and the cotangent complex replaced by graded commutative algebra (Lemma~\ref{lem:chart}).

\subsection*{Errata and generalizations}
Section~\ref{sec:errata} collects a list of errata to \cite{AlgMag}: statements false as printed, statements whose printed proof uses an additional hypothesis, an implicit hypothesis, misprints, and two replaced statements. Section~\ref{sec:generalizations} lists the statements of \cite{AlgMag,SelfActions,XProducts} proved here in greater generality: for every commutative monoid, for every face, or without finiteness or separatedness conditions. The formalization confirms the theory of \cite{AlgMag,SelfActions,XProducts} and strengthens it: the errata concern minor hypotheses and a few statements, while several results hold in greater generality and three conjectures of \cite{AlgMag} are proved.

\section{Conventions}\label{sec:conventions}

Monoids and rings are commutative; monoids are written additively and rings have a unit. For a monoid $P$ and a subset $E\subseteq P$, the submonoid generated by $E$ is denoted $\mon{E}$. For a monoid $N$, the set $\units N=\{x\in N\mid \exists\, y\in N,\ x+y=0\}$ of invertible elements is a subgroup of $N$; if $N$ is a submonoid of an abelian group, $N^{\gp}=\{x-y\mid x,y\in N\}$ is the subgroup it generates. A monoid $N$ is \emph{cancellative} if $x+y=x+z$ implies $y=z$, \emph{sharp} if $\units N=0$, and \emph{fine} if it is cancellative and finitely generated. The set of natural numbers $\NN$ contains $0$.

We fix a base scheme $S$ and work with the big fppf site $(\Sch/S)_{\fppf}$ of \cite[Tag~021R]{Stacks}. Its objects are the schemes over $S$ whose size is bounded by a fixed infinite cardinal $\kappa$ \cite[Tag~000H]{Stacks}. Throughout, $\kappa$ is chosen larger than the cardinality of every monoid under consideration and of every index set of a family of charts under consideration; then for every object $T$ of the site and every such monoid $N$, the scheme $\A(N)_T$ of Section~\ref{sec:monoids} is again an object of the site, and so are the disjoint unions of charts used in Section~\ref{sec:representability}. This is the only size condition, and it is not repeated in the statements.
\leansrc{site}

An \emph{algebraic space over $S$} is a sheaf of sets on $(\Sch/S)_{\fppf}$ which is an algebraic space in the sense of \cite[Tag~025Y]{Stacks}; schemes over $S$ are identified with the sheaves they represent. Properties of algebraic spaces over $S$ such as separated, quasi-separated, quasi-compact, locally of finite presentation and of finite presentation refer to the structure morphism to $S$, and are those of \cite{Stacks}. A morphism $X\to S$ is \emph{$S$-affine} if it is affine; an algebraic space $X$ over $S$ is then called $S$-affine and is identified with $\Spec_S(\cA)$ for the quasi-coherent $\cO_S$-algebra $\cA$ of $X$. For an algebraic space $X$ over $S$ and a morphism of schemes $T\to S$ we write $X_T=X\times_S T$.
\leansrc{algspace}

A statement cited from \cite{AlgMag}, \cite{SelfActions} or \cite{XProducts} keeps the number it has in the published version. Every definition and statement of this text outside Section~\ref{sec:reductive} is accompanied by the name of the Lean declaration that formalizes it.
The Lean library \texttt{algmag-stacks} \cite{AlgMagLean}, at commit \texttt{\algmagcommit}, is written for Lean~4 version \texttt{v4.34.0-rc2} and depends on Mathlib \cite{Mathlib} at commit \texttt{16edf57} and on AxiomMath Stacks \cite{AxiomStacks} at commit \texttt{012108f}.

\section{Monoids, diagonalizable monoid schemes and faces}\label{sec:monoids}

\subsection{Rings of monoids and diagonalizable monoid schemes}

\begin{definition}[{\cite[Def.~2.1, Def.~2.3(i)]{AlgMag}}]\label{def:monoid-ring}
Let $N$ be a monoid. The \emph{ring of $N$} is the ring $\ZZ[N]=\bigoplus_{n\in N}\ZZ X^n$ with multiplication $X^n X^{n'}=X^{n+n'}$ and unit $X^0=1$. For a ring $R$, put $R[N]=\ZZ[N]\otimes_\ZZ R=\bigoplus_{n\in N}RX^n$.
\leansrc{monoidRing}
\end{definition}

\begin{proposition}[{\cite[Fact~2.2]{AlgMag}}]\label{prop:bialgebra}
Let $N$ be a monoid. The ring $\ZZ[N]$ is a bialgebra over $\ZZ$ with counit $X^n\mapsto 1$ and comultiplication $X^n\mapsto X^n\otimes X^n$. If $N$ is a group, $\ZZ[N]$ is a Hopf algebra with antipode $X^n\mapsto X^{-n}$.
\leansrc{bialgebra}
\end{proposition}

\begin{proof}
\cite[Fact~2.2]{AlgMag}
\end{proof}

\begin{definition}[{\cite[Def.~2.4, Def.~2.11(ii)]{AlgMag}}]\label{def:diag-monoid}
Let $N$ be a monoid. The monoid scheme $\A(N)=\Spec(\ZZ[N])$ over $\Spec(\ZZ)$ is the \emph{diagonalizable monoid scheme} of $N$, and $\A(N)_S=\A(N)\times_{\Spec(\ZZ)}S$ is its base change to $S$; the monoid structure is given by the bialgebra structure of Proposition~\ref{prop:bialgebra}. If $N=M$ is an abelian group, $\A(M)$ is a group scheme, denoted $\D(M)$, and $\D(M)_S=\A(M)_S$ is the \emph{diagonalizable group scheme} of $M$ over $S$. The trivial monoid scheme over $S$ is $e_S=\A(0)_S$, which is $S$ as a scheme.
\leansrc{diagMonoid}
\end{definition}

\begin{proposition}[{\cite[Fact~2.6(i), Fact~2.7]{AlgMag}}]\label{prop:diag-affine}
Let $N$ and $L$ be monoids.
\begin{enumerate}
\item The morphism $\A(N)_S\to S$ is affine and flat, and its quasi-coherent $\cO_S$-algebra is the sheaf of bialgebras $\cO_S[N]$, the sheafification of $U\mapsto \cO_S(U)[N]$.
\item There are canonical isomorphisms of monoid schemes $\A(N\times L)=\A(N)\times_{\Spec(\ZZ)}\A(L)$ and $\A(N\times L)_S=\A(N)_S\times_S\A(L)_S$.
\end{enumerate}
\leansrc{diagAffine}
\end{proposition}

\begin{proof}
\cite[proofs of Fact~2.6 and Fact~2.7]{AlgMag}
\end{proof}

\begin{proposition}[{\cite[Fact~2.10]{AlgMag}}]\label{prop:diag-fp}
Let $N$ be a finitely generated monoid. Then $N$ is finitely presented as a monoid, $\ZZ[N]$ is a finitely presented $\ZZ$-algebra, and the morphisms $\A(N)_S\to S$ and $\A(N)_S\times_S\A(N)_S\to\A(N)_S$ (multiplication) are of finite presentation, and $\A(N)_S\to S$ is affine.
\leansrc{diagFP}
\end{proposition}

\begin{proof}
\cite[proof of Fact~2.10]{AlgMag}
\end{proof}

\subsection{Faces}

Let $M$ be a monoid. For submonoids $F\subseteq N$ of $M$, the \emph{projection} $\ZZ[N]\to\ZZ[F]$ is the $\ZZ$-linear map sending $X^n$ to $X^n$ if $n\in F$ and to $0$ if $n\in N\setminus F$.

\begin{definition}[{\cite[Def.~2.14]{AlgMag}}]\label{def:face}
Let $F\subseteq N$ be submonoids of $M$. Then $F$ is a \emph{face} of $N$ if the projection $\ZZ[N]\to\ZZ[F]$ is a morphism of rings.
\leansrc{face}
\end{definition}

\begin{proposition}[{\cite[Prop.~2.16]{AlgMag}}]\label{prop:face-criterion}
Let $F\subseteq N$ be submonoids of $M$. Then $F$ is a face of $N$ if and only if for all $x,y\in N$, $x+y\in F$ implies $x\in F$ and $y\in F$. The same holds with $\ZZ$ replaced by any nonzero commutative ring $R$, for the projection $R[N]\to R[F]$.
\leansrc{faceCriterion}
\end{proposition}

\begin{proof}
\cite[proof of Prop.~2.16]{AlgMag}
\end{proof}

\begin{proposition}[{\cite[Prop.~2.17, Fact~2.18]{AlgMag}}]\label{prop:units-face}
Let $N$ be a monoid. Then $\units N$ is a submonoid of $N$ and a group; it is the largest subgroup of $N$, and it is a face of $N$. Every face of $N$ contains $\units N$.
\leansrc{unitsFace}
\end{proposition}

\begin{proof}
\cite[proofs of Prop.~2.17 and Fact~2.18]{AlgMag}
\end{proof}

The face $\units N$ is the \emph{face of invertible elements} of $N$.

\begin{proposition}[{\cite[Prop.~2.19]{AlgMag}}]\label{prop:face-pushout}
Let $N\subseteq L\subseteq L'$ be submonoids of $M$ such that $L$ is a face of $L'$. Then $N'=L'\setminus(L\setminus N)$ is a submonoid of $M$ and $N$ is a face of $N'$. Moreover the square
\[
\begin{tikzcd}
\A(L)_S \arrow[r] \arrow[d] & \A(L')_S \arrow[d]\\
\A(N)_S \arrow[r] & \A(N')_S
\end{tikzcd}
\]
induced by the inclusions $N\subseteq L$, $N'\subseteq L'$ and by the projections $\ZZ[L']\to\ZZ[L]$, $\ZZ[N']\to\ZZ[N]$ is a push-out in the category of schemes.
\leansrc{facePushout}
\end{proposition}

\begin{proof}
\cite[proof of Prop.~2.19]{AlgMag}
\end{proof}

\begin{lemma}\label{lem:face-weight}
Let $N$ be a finitely generated submonoid of an abelian group $M$ and let $F$ be a face of $N$. Then there is a morphism of monoids $w\colon N\to\NN$ with $w^{-1}(0)=F$.
\leansrc{faceWeight}
\end{lemma}

\begin{proof}
The image $P$ of $N$ in $M/F^{\gp}$ is finitely generated and cancellative. It is sharp: if $x,y\in N$ and $x+y=f_1-f_2$ with $f_1,f_2\in F$, then $x+(y+f_2)=f_1\in F$ with $x,y+f_2\in N$, so $x\in F$ and $y+f_2\in F$ by Proposition~\ref{prop:face-criterion}; applying Proposition~\ref{prop:face-criterion} again to $y+f_2\in F$ gives $y\in F$. By \cite[Prop.~2.2.1, p.~34]{Ogus} there is a morphism $h\colon P\to\NN$ with $h^{-1}(0)=0$. Let $w$ be the composite $N\to P\to\NN$. If $w(x)=0$ then $x\in N\cap F^{\gp}$, say $x=f_1-f_2$, so $x+f_2\in F$ and $x\in F$; conversely $w(F)=0$.
\end{proof}

\section{Actions, magnets and attractors}\label{sec:attractors}

Let $M$ be a monoid.

\subsection{Actions of diagonalizable monoid schemes}

\begin{definition}[{\cite[\S3]{AlgMag}}]\label{def:action}
Let $X$ be a sheaf on $(\Sch/S)_{\fppf}$. An \emph{action} of $\A(M)_S$ on $X$ is an action of the monoid $\A(M)_S(T)$ on $X(T)$ for every object $T$, natural in $T$. A morphism $f\colon X\to Y$ of sheaves with $\A(M)_S$-actions is \emph{equivariant} if $f(g\cdot x)=g\cdot f(x)$ for all $T$, $g\in\A(M)_S(T)$ and $x\in X(T)$. For an object $T$, $\Homeq{\A(M)_T}_T(X_T,Y_T)$ denotes the set of $\A(M)_T$-equivariant $T$-morphisms $X_T\to Y_T$.
\leansrc{action}
\end{definition}

\begin{proposition}[{\cite[Fact~2.12, Fact~2.13, Fact~2.15]{AlgMag}}]\label{prop:action-AN}
Let $N\subseteq N'$ be submonoids of $M$.
\begin{enumerate}
\item The ring map $\ZZ[N]\to\ZZ[N]\otimes_\ZZ\ZZ[M]$, $X^n\mapsto X^n\otimes X^n$, defines an action of $\A(M)$ on $\A(N)$; it is the action through the morphism of monoid schemes $\A(M)\to\A(N)$ induced by $N\subseteq M$. By base change, $\A(M)_S$ acts on $\A(N)_S$.
\item The inclusion $\ZZ[N]\to\ZZ[N']$ induces an $\A(M)_S$-equivariant morphism of monoid schemes $\A(N')_S\to\A(N)_S$.
\item If $N$ is a face of $N'$, the projection $\ZZ[N']\to\ZZ[N]$ induces an $\A(M)_S$-equivariant morphism $\A(N)_S\to\A(N')_S$.
\end{enumerate}
\leansrc{actionAN}
\end{proposition}

\begin{proof}
\cite[Fact~2.12, proofs of Fact~2.13 and Fact~2.15]{AlgMag}
\end{proof}

\begin{definition}[{\cite[Rem.~3.5]{AlgMag}}]\label{def:fixed-points}
Let $X$ be a sheaf with an action of $\A(M)_S$. The \emph{fixed-point functor} $X^{\A(M)_S}$ is the subfunctor of $X$ given by $X^{\A(M)_S}(T)=\{x\in X(T)\mid g\cdot x_{T'}=x_{T'}$ for all $T'\to T$ and $g\in\A(M)_S(T')\}$.
\leansrc{fixedPoints}
\end{definition}

\subsection{Attractors}

From now on, $X$ is an algebraic space over $S$ with an action of $\A(M)_S$. A \emph{magnet} is a submonoid of $M$.

\begin{definition}[{\cite[Def.~3.1]{AlgMag}}]\label{def:attractor}
Let $N$ be a magnet. The \emph{attractor} $X^N$ is the functor on $(\Sch/S)_{\fppf}$ sending $T$ to $\Homeq{\A(M)_T}_T(\A(N)_T,X_T)$, with the action of $\A(M)_T$ on $\A(N)_T$ of Proposition~\ref{prop:action-AN}, and sending $T'\to T$ to the base change of equivariant morphisms. When the monoid has to be recorded, $X^N$ is also written $X^N_M$.
\leansrc{attractor}
\end{definition}

The definition applies to any sheaf on $(\Sch/S)_{\fppf}$ with an action of $\A(M)_S$, and Propositions~\ref{prop:attractor-sheaf}, \ref{prop:attractor-extreme}, \ref{prop:attractor-limits}, \ref{prop:magnet-maps}, \ref{prop:attractor-actions} and~\ref{prop:subgroup-fixed} hold in this generality.

\begin{proposition}[{\cite[Prop.~3.2]{AlgMag}}]\label{prop:attractor-sheaf}
For every magnet $N$, the attractor $X^N$ is a sheaf on $(\Sch/S)_{\fppf}$.
\leansrc{attractorSheaf}
\end{proposition}

\begin{proof}
\cite[proof of Prop.~3.2]{AlgMag}
\end{proof}

\begin{proposition}[{\cite[Rem.~3.5, Fact~3.6]{AlgMag}}]\label{prop:attractor-extreme}
Evaluation at the unit section $e_T\to\A(M)_T$ identifies $X^M$ with $X$, the inverse sending $x\in X(T)$ to its orbit map $h\mapsto h\cdot x$. Evaluation at the unique point of $\A(0)_T=T$ identifies $X^0$ with the fixed-point functor $X^{\A(M)_S}$.
\leansrc{attractorExtreme}
\end{proposition}

\begin{proof}
\cite[Rem.~3.5 and Fact~3.6]{AlgMag}
\end{proof}

\begin{proposition}[{\cite[Prop.~3.7, Prop.~3.8, Prop.~3.9]{AlgMag}}]\label{prop:attractor-limits}
Let $N$ be a magnet.
\begin{enumerate}
\item Let $X_1\to X_2\leftarrow X_3$ be equivariant morphisms of algebraic spaces over $S$ with $\A(M)_S$-actions. Then $\A(M)_S$ acts diagonally on $X_1\times_{X_2}X_3$, the projections are equivariant, the square formed by $X_1\times_{X_2}X_3$ is cartesian in the category of algebraic spaces over $S$ with $\A(M)_S$-action, and $(X_1\times_{X_2}X_3)^N=X_1^N\times_{X_2^N}X_3^N$.
\item For algebraic spaces $Y$, $Z$ over $S$ with $\A(M)_S$-actions, $(Y\times_SZ)^N=Y^N\times_SZ^N$.
\item For every morphism of schemes $T\to S$, $X^N\times_ST=(X_T)^N$ as functors on $T$-schemes, where $X_T$ carries the action of $\A(M)_T$.
\end{enumerate}
\leansrc{attractorLimits}
\end{proposition}

\begin{proof}
\cite[proofs of Prop.~3.7, Prop.~3.8 and Prop.~3.9]{AlgMag}
\end{proof}

\begin{proposition}[{\cite[Fact~3.10, Fact~3.15, Fact~3.17, Fact~3.18]{AlgMag}}]\label{prop:magnet-maps}
Let $N\subseteq L$ be magnets and let $f\colon X\to Y$ be an equivariant morphism of algebraic spaces over $S$ with $\A(M)_S$-actions.
\begin{enumerate}
\item Precomposition with $\A(L)_T\to\A(N)_T$ defines a morphism $\iota_{N,L}\colon X^N\to X^L$.
\item If $N$ is a face of $L$, precomposition with $\A(N)_T\to\A(L)_T$ defines a morphism $p_{L,N}\colon X^L\to X^N$, and $p_{L,N}\circ\iota_{N,L}=\id_{X^N}$.
\item Composition with $f$ defines a morphism $f^N\colon X^N\to Y^N$; if $f$ is a monomorphism, so is $f^N$.
\end{enumerate}
\leansrc{magnetMaps}
\end{proposition}

\begin{proof}
\cite[Fact~3.10, Fact~3.15, Fact~3.17, proof of Fact~3.18]{AlgMag}
\end{proof}

For $L=M$, the composite of $\iota_{N,M}$ with the identification $X^M=X$ of Proposition~\ref{prop:attractor-extreme} is the morphism $X^N\to X$, $\varphi\mapsto\varphi(1)$, evaluation at the unit section of $\A(N)_T$.

\begin{proposition}[{\cite[Fact~3.16]{AlgMag}}]\label{prop:face-closed}
Let $F$ be a face of a magnet $N$. Assume that $X^N$ and $X^F$ are algebraic spaces and that $p_{N,F}\colon X^N\to X^F$ is separated. Then $\iota_{F,N}\colon X^F\to X^N$ is a closed immersion.
\leansrc{faceClosed}
\end{proposition}

\begin{proof}
\cite[proof of Fact~3.16]{AlgMag}
\end{proof}

\begin{proposition}[{\cite[Fact~3.19]{AlgMag}}]\label{prop:attractor-lfp}
Let $f\colon X\to Y$ be an equivariant morphism of algebraic spaces over $S$ with $\A(M)_S$-actions, locally of finite presentation, and let $N$ be a magnet. Then $f^N\colon X^N\to Y^N$ is limit preserving: for every directed system of affine schemes $T_i=\Spec(R_i)$ over $S$ with limit $T=\Spec(\colim R_i)$, the map $\colim_iX^N(T_i)\to X^N(T)\times_{Y^N(T)}\colim_iY^N(T_i)$ is bijective. In particular, if $X^N$ and $Y^N$ are algebraic spaces, $f^N$ is locally of finite presentation.
\leansrc{attractorLfp}
\end{proposition}

\begin{proof}
\cite[proof of Fact~3.19]{AlgMag}, which uses the definition \cite[Tag~049J]{Stacks} of limit preserving transformations and the characterization \cite[Tag~04AK]{Stacks} of morphisms locally of finite presentation; the computation holds for the functors $X^N$, $Y^N$ without assuming that they are algebraic spaces.
\end{proof}

\begin{proposition}[{\cite[Rem.~3.11, Rem.~3.12, Rem.~3.13, Prop.~3.14]{AlgMag}}]\label{prop:attractor-actions}
Let $N$ be a magnet.
\begin{enumerate}
\item The monoid scheme $\A(N)_S$ acts on $X^N$ by $(t\cdot\varphi)(x)=\varphi(tx)$, and $\A(M)_S$ acts on $X^N$ by $(h\cdot\varphi)(x)=h\cdot\varphi(x)=\varphi(hx)$; the second action is the first one through $\A(M)_S\to\A(N)_S$.
\item If $X$ is a monoid (resp.\ group) algebraic space over $S$ and $\A(M)_S$ acts by monoid endomorphisms (resp.\ group endomorphisms), then $X^N$ is a monoid (resp.\ group) functor, with the pointwise law.
\end{enumerate}
\leansrc{attractorActions}
\end{proposition}

\begin{proof}
\cite[Rem.~3.11, Rem.~3.12, Rem.~3.13, proof of Prop.~3.14]{AlgMag}
\end{proof}

\subsection{Separated algebraic spaces}

\begin{proposition}[{\cite[Prop.~3.29, Cor.~3.30]{AlgMag}}]\label{prop:separated-mono}
Assume that $X$ is separated and let $N\subseteq L$ be magnets. Then $\iota_{N,L}\colon X^N\to X^L$ is a monomorphism. In particular $X^N\to X$ is a monomorphism.
\leansrc{separatedMono}
\end{proposition}

\begin{proof}
\cite[proofs of Prop.~3.29 and Cor.~3.30]{AlgMag}
\end{proof}

\begin{proposition}[{\cite[Prop.~3.37]{AlgMag}}]\label{prop:colimit-magnets}
Assume that $X$ is separated and locally of finite presentation. Let $N$ be a magnet, written as a directed union $N=\bigcup_{i\in I}N_i$ of submonoids. Then for every quasi-compact scheme $T$ over $S$,
\[X^N(T)=\colim_{i\in I}X^{N_i}(T).\]
\leansrc{colimitMagnets}
\end{proposition}

\begin{proof}
For $T$ affine this is \cite[proof of Prop.~3.37]{AlgMag}. Let $T$ be quasi-compact and let $T=\bigcup_{k=1}^rV_k$ be a finite affine open cover. The maps $X^{N_i}(T)\to X^N(T)$ are injective by Proposition~\ref{prop:separated-mono}. Given $\varphi\in X^N(T)$, each restriction $\varphi|_{V_k}$ comes from $X^{N_{i_k}}(V_k)$; choose $i\geq i_1,\dots,i_r$. The restrictions of these points to $X^{N_i}(V_k\cap V_l)$ have the same image in $X^N(V_k\cap V_l)$, hence agree, and they glue to a point of $X^{N_i}(T)$ with image $\varphi$, since $X^{N_i}$ is a sheaf (Proposition~\ref{prop:attractor-sheaf}).
\end{proof}

\subsection{Subgroups and fixed points}

\begin{proposition}[{\cite[Prop.~3.35]{AlgMag}}]\label{prop:subgroup-fixed}
Assume that $M$ is a finitely generated abelian group and let $Z\subseteq M$ be a subgroup. Then evaluation at the unit section identifies the attractor $X^Z$ for the action of $\D(M)_S$ with the fixed-point functor $X^{\D(M/Z)_S}$ for the action of $\D(M/Z)_S\subseteq\D(M)_S$. This holds for every sheaf $X$ on $(\Sch/S)_{\fppf}$ with an action of $\D(M)_S$. The same identification holds if $M$ is an arbitrary abelian group and $X$ is an algebraic space over $S$ with a $Z$-FPR atlas (Definition~\ref{def:fpr}).
\leansrc{subgroupFixed}
\end{proposition}

\begin{proof}
\cite[proof of Prop.~3.35]{AlgMag}. The proof uses only that $X$ is an fppf sheaf and that $\D(M)_S\to\D(Z)_S$ is faithfully flat and of finite presentation. In the second case, the $Z$-FPR atlas $U\to X$ provides the extension of values to $\D(M)_S$ which the printed proof obtains from the faithfully flat morphism $\D(M)_S\to\D(Z)_S$.
\end{proof}

\subsection{Prescribed limits}

\begin{definition}[{\cite[Def.~4.1]{AlgMag}}]\label{def:prescribed}
Let $N$ be a magnet, $F$ a face of $N$ and $\mathcal Z\to X^F$ a monomorphism of functors. The \emph{attractor with prescribed limit} $X^N_{F,\mathcal Z}$ is the subfunctor of $X^N$ of points $\varphi\in X^N(T)$ with $p_{N,F}(\varphi)\in\mathcal Z(T)$. For $F=\units N$ it is denoted $X^N_{\mathcal Z}$.
\leansrc{prescribed}
\end{definition}

\begin{proposition}[{\cite[Prop.~4.2, Fact~4.3]{AlgMag}}]\label{prop:prescribed}
In the situation of Definition~\ref{def:prescribed}, $X^N_{F,\mathcal Z}=X^N\times_{X^F}\mathcal Z$. The morphism $X^N_{F,\mathcal Z}\to X^N$ is a monomorphism, as a base change of the monomorphism $\mathcal Z\to X^F$; this holds without separatedness of $X$. For an arbitrary morphism $\mathcal Z\to X^F$ in place of the monomorphism, one defines $X^N_{F,\mathcal Z}$ as $X^N\times_{X^F}\mathcal Z$. If $X$ is separated, $X^N_{F,\mathcal Z}\to X$ is a monomorphism.
\leansrc{prescribedFibre}
\end{proposition}

\begin{proof}
\cite[proofs of Prop.~4.2 and Fact~4.3]{AlgMag}
\end{proof}

\section{Attractors of affine spaces}\label{sec:affine}

Let $M$ be a monoid.

\subsection{Gradings}

\begin{definition}[{\cite[Def.~2.20]{AlgMag}}]\label{def:AM-module}
For an $\cO_S$-module $\cE$, let $W(\cE)$ be the functor $T\mapsto\Gamma(T,p^*\cE)$ on schemes $p\colon T\to S$; it is a module over the ring functor $T\mapsto\Gamma(T,\cO_T)$. An \emph{$\A(M)_S$-$\cO_S$-module} is an $\cO_S$-module $\cE$ together with an action of the monoid functor $\A(M)_S$ on $W(\cE)$ by linear maps: for every $T$, the monoid $\A(M)_S(T)$ acts on $W(\cE)(T)$, with $h\cdot(x+\lambda y)=h\cdot x+\lambda(h\cdot y)$ for $\lambda\in\Gamma(T,\cO_T)$, compatibly with the maps induced by morphisms $T'\to T$.
\leansrc{AMmodule}
\end{definition}

\begin{proposition}[{\cite[Prop.~2.21, Prop.~2.22]{AlgMag}}]\label{prop:gradings}
\begin{enumerate}
\item The category of quasi-coherent $\A(M)_S$-$\cO_S$-modules is equivalent to the category of quasi-coherent $M$-graded $\cO_S$-modules $\cE=\bigoplus_{m\in M}\cE_m$.
\item Let $X=\Spec_S(\cA)$ be an $S$-affine scheme. Actions of $\A(M)_S$ on $X$ correspond to $M$-gradings $\cA=\bigoplus_{m\in M}\cA_m$ of the quasi-coherent $\cO_S$-algebra $\cA$.
\end{enumerate}
\leansrc{gradings}
\end{proposition}

\begin{proof}
\cite[proofs of Prop.~2.21 and Prop.~2.22]{AlgMag}
\end{proof}

\subsection{The affine theorem}

\begin{theorem}[{\cite[Thm.~3.20]{AlgMag}}]\label{thm:affine}
Let $X=\Spec_S(\cA)$ be an $S$-affine scheme with an action of $\A(M)_S$, and let $\cA=\bigoplus_{m\in M}\cA_m$ be the associated grading (Proposition~\ref{prop:gradings}). Let $N$ be a magnet. Then $X^N$ is represented by the closed subscheme of $X$ defined by the quasi-coherent ideal $\cJ_N$ generated by the $\cA_m$, $m\in M\setminus N$. A point of $X^N(T)$ is the morphism $\A(N)_T\to X_T$ given by the $M$-graded homomorphism $p^*\cA\to\cO_T[N]$, $a_m\mapsto F(a_m)X^m$, where $p\colon T\to S$ is the structure morphism and $F\colon p^*\cA/p^*\cJ_N\to\cO_T$ corresponds to the $T$-point of $\Spec_S(\cA/\cJ_N)$.
\leansrc{affineTheorem}
\end{theorem}

\begin{proof}
\cite[proof of Thm.~3.20]{AlgMag}
\end{proof}

\begin{example}[{\cite[Ex.~3.21]{AlgMag}, \cite[Prop.~2.1]{SelfActions}}]\label{ex:AL}
Let $L$ and $N$ be submonoids of $M$ and let $\A(M)_S$ act on $X=\A(L)_S$ as in Proposition~\ref{prop:action-AN}. Then $X^N$ is the closed subscheme $\Spec_S(\cO_S[L]/\cJ_N)$ of $X$, where $\cJ_N$ is the ideal generated by the $X^m$, $m\in L\setminus(N\cap L)$.
\leansrc{exampleAL}
\end{example}

\begin{example}[{\cite[Ex.~1.2]{AlgMag}}]\label{ex:linear}
Let $\cE=\bigoplus_{m\in M}\cE_m$ be an $M$-graded $\cO_S$-module, let $V=\mathbb V(\cE)$ with the induced linear action of $\A(M)_S$, and let $N$ be a submonoid of $M$. For a scheme $p\colon T\to S$, a $T$-point of $V$ is an $\cO_S$-linear map $\xi\colon\cE\to p_*\cO_T$. Then for every $T$ the map $V^N(T)\to V(T)$, evaluation at the unit section, is injective, natural in $T$, and its image is the set of $\xi$ with $\xi|_{\cE_m}=0$ for every $m\in M\setminus N$. In particular, for $V=\mathbb{A}^n_S$ with weights $m_1,\dots,m_n\in M$ on the coordinates, $V=\bigoplus_{m\in M}V_m$ and $V^N=\bigoplus_{m\in N}V_m$.
\leansrc{exampleLinear}
\end{example}

\begin{proposition}[{\cite[Prop.~3.28]{AlgMag}}]\label{prop:vector-bundle}
Let $\cE=\bigoplus_{m\in M}\cE_m$ be a quasi-coherent $M$-graded $\cO_S$-module and let $\mathbb V(\cE)=\Spec_S(\Sym\cE)$ with the induced linear action of $\A(M)_S$. For every magnet $N$, the attractor $\mathbb V(\cE)^N$ is canonically isomorphic to $\mathbb V(\cE_N)$, where $\cE_N=\bigoplus_{n\in N}\cE_n$.
\leansrc{vectorBundle}
\end{proposition}

\begin{proof}
\cite[proof of Prop.~3.28]{AlgMag}
\end{proof}

\begin{corollary}[{\cite[Cor.~3.22, Prop.~3.24, Lem.~3.26]{AlgMag}}]\label{cor:affine-consequences}
Let $X$ be an $S$-affine scheme with an action of $\A(M)_S$.
\begin{enumerate}
\item For magnets $N\subseteq L$, the morphism $\iota_{N,L}\colon X^N\to X^L$ is a closed immersion.
\item For a family $(N_i)_{i\in I}$ of magnets, $\bigcap_{i\in I}X^{N_i}=X^{\bigcap_iN_i}$, the left-hand side being the scheme-theoretic intersection of closed subschemes of $X$; equivalently, for every scheme $T$ over $S$, a $T$-point of $X$ lies in $X^{\bigcap_iN_i}$ if and only if it lies in every $X^{N_i}$.
\item For an equivariant closed immersion $Z\to Y$ of $S$-affine schemes with $\A(M)_S$-actions and a magnet $N$, $Z^N=Z\times_YY^N$; in particular $Z^N\to Y^N$ is a closed immersion.
\end{enumerate}
\leansrc{affineConsequences}
\end{corollary}

\begin{proof}
\cite[proofs of Cor.~3.22, Prop.~3.24 and Lem.~3.26]{AlgMag}
\end{proof}

\subsection{Self-actions}\label{sec:self-actions}

Let $\units M$ be the face of invertible elements of $M$ and let $M/\units M$ be the quotient monoid, with projection $q\colon M\to M/\units M$; it is sharp. The monoid scheme $\A(M)_S$ acts on itself, and $\A(M/\units M)_S$ acts on itself. For a magnet $N$ of $M$, $\A(M)_S^N$ is a closed subscheme of $\A(M)_S$ by Theorem~\ref{thm:affine}.

\begin{proposition}[{\cite[Lem.~4.1, Prop.~4.2, Prop.~4.3]{SelfActions}}]\label{prop:self-action}
Let $N$ and $N'$ be magnets of $M$.
\begin{enumerate}
\item If $\units M\not\subseteq N$, then $\A(M)_S^N=\emptyset$. In particular, if $M$ is a group and $N\neq M$, then $\A(M)_S^N=\emptyset$.
\item If $\units M\subseteq N$ and $\units M\subseteq N'$, then $\A(M)_S^N=\A(M)_S^{N'}$ as closed subschemes of $\A(M)_S$ if and only if $\A(M/\units M)_S^{q(N)}=\A(M/\units M)_S^{q(N')}$ as closed subschemes of $\A(M/\units M)_S$.
\end{enumerate}
\leansrc{selfAction}
\end{proposition}

\begin{proof}
\cite[proofs of Lem.~4.1, Prop.~4.2 and Prop.~4.3]{SelfActions}
\end{proof}

\subsection{Change of monoid and intersections of magnets}

\begin{proposition}[{\cite[Prop.~3.31]{AlgMag}}]\label{prop:change-monoid}
Let $f\colon M\to Z$ be a morphism of monoids, let $Y$ be a magnet of $Z$ and let $N=f^{-1}(Y)$. Let $\A(Z)_S$ act on $X$ through the morphism $\A(Z)_S\to\A(M)_S$ dual to $f$. Assume that one of the following conditions holds:
\begin{enumerate}
\item $X$ is $S$-affine;
\item $Z$ is a finitely generated abelian group and one of the following holds:
\begin{enumerate}[label=(\alph*)]
\item $M$ is countable and $X$ is a scheme;
\item $M$ is a finitely generated abelian group and $X$ is an algebraic space over $S$;
\end{enumerate}
\item $M$ is cancellative, $Z=M^{\gp}$, $f$ is the inclusion, $Y=N$, and $X$ is separated.
\end{enumerate}
Then $X^N_M=X^Y_Z$.
\leansrc{changeMonoid}
\end{proposition}

\begin{proof}
For (i) and (iii), see \cite[proof of Prop.~3.31]{AlgMag}. Case (ii)(a), where $M$ is only a countable monoid, is not covered by the printed proof, written for groups; it is proved in the formal proof. For (ii)(b), the formal proof replaces the fpqc quotient of the printed proof by an fppf cover. Choose a surjection $q\colon\ZZ^r\to Z$. For an affine scheme $T$ over $S$ and a morphism of monoids $P\to Z$, $\Spec\cO(T)[P\times\ZZ^r]\to\Spec\cO(T)[P\times_Z\ZZ^r]$ is free of finite presentation, being the base change of $\cO(T)[\ker q]\to\cO(T)[\ZZ^r]$, and it is an fppf $\D(Z)$-torsor; this reduces the statement to the case $Z=\ZZ^r$. That case is proved for every algebraic space $X$ by a local factorization of equivariant morphisms on étale charts of $X$.
\end{proof}

\begin{proposition}[{\cite[Prop.~3.32]{AlgMag}}]\label{prop:intersection}
Let $N$ and $L$ be magnets, and let $\A(M)_S$ act on $X^N$ as in Proposition~\ref{prop:attractor-actions}. Assume that one of the following conditions holds:
\begin{enumerate}
\item $X$ is $S$-affine;
\item $M$ is a finitely generated abelian group;
\item $M$ is cancellative, $M^{\gp}$ is finitely generated, and $X$ is separated.
\end{enumerate}
Then $(X^N)^L=X^{N\cap L}$.
\leansrc{intersection}
\end{proposition}

\begin{proof}
For (i) and (ii), see \cite[proof of Prop.~3.32]{AlgMag}. For (iii), the proof of \cite[Prop.~3.32(iii)]{AlgMag} applies Proposition~\ref{prop:change-monoid}(iii) to $X$ and to $X^N$, and then (ii). The proof of Proposition~\ref{prop:change-monoid}(iii) uses the separatedness of $X$ only through the following property: if two morphisms $g_1,g_2\colon W\to X_T$ of $T$-schemes agree after composition with a morphism $W'\to W$ whose scheme-theoretic image is $W$, then $g_1=g_2$. This property holds for $X$ separated, since the equalizer of $g_1$ and $g_2$ is then a closed subscheme of $W$ through which $W'\to W$ factors. Since $X^N\to X$ is a monomorphism (Proposition~\ref{prop:separated-mono}), the property passes from $X$ to $X^N$, and the argument applies to $X^N$ without assuming $X^N$ separated.
\end{proof}

\section{Fixed-point-reflecting atlases and representability}\label{sec:representability}

In this section $M$ is an abelian group and $X$ is an algebraic space over $S$ with an action of $\D(M)_S$.

\subsection{\'Etale morphisms and faces}

\begin{proposition}[{\cite[Prop.~5.1]{AlgMag}}]\label{prop:etale-face}
Let $N$ be a magnet, $Q$ a face of $N$, $U$ an $S$-affine scheme with an action of $\D(M)_S$ and $U\to X$ an \'etale equivariant morphism. Assume that $X$ is separated and that one of the following conditions holds:
\begin{enumerate}
\item $N$ is finitely generated;
\item $X$ is locally of finite presentation.
\end{enumerate}
Then the natural morphism $U^N\to X^N\times_{X^Q}U^Q$ is an isomorphism.
\leansrc{etaleFace}
\end{proposition}

\begin{proof}
(i) We follow \cite[proof of Prop.~5.1(i)]{AlgMag}. Let $T$ be a scheme over $S$ and let $\varphi\colon\A(N)_T\to X_T$ be an equivariant morphism whose restriction to $\A(Q)_T$ lifts to an equivariant morphism $\psi_0\colon\A(Q)_T\to U_T$. Let $w\colon N\to\NN$ be a morphism of monoids with $w^{-1}(0)=Q$ (Lemma~\ref{lem:face-weight}). The ideal $I$ of $\A(Q)_T$ in $\A(N)_T$ is spanned by the $X^n$, $n\in N\setminus Q$, so $I^k$ is spanned by the $X^n$ with $n$ a sum of $k$ elements of $N\setminus Q$ and of an element of $N$, and such $n$ satisfy $w(n)\geq k$. The proof constructs, by formal smoothness of $U\to X$, equivariant lifts $\psi_k$ of $\varphi$ on the infinitesimal neighbourhoods $V(I^k)$ of $\A(Q)_T$ in $\A(N)_T$, and an equivariant morphism $\psi\colon\A(N)_T\to U_T$ restricting to $\psi_k$ on $V(I^k)$ for every $k$. Write $U=\Spec_S(\cA)$ with the grading of Proposition~\ref{prop:gradings}(ii) and let $p\colon T\to S$ be the structure morphism; the $\psi_k$ are compatible graded homomorphisms $p^*\cA\to\cO_T[N]/I^k$. For $a$ homogeneous of degree $m\in N$ and $k>w(m)$, the degree-$m$ part of $\cO_T[N]/I^k$ is $\cO_T\cdot X^m$, and for $m\in M\setminus N$ it is $0$; so these homomorphisms define a graded homomorphism $p^*\cA\to\cO_T[N]$, which is $\psi$. The uniqueness of the lifts gives injectivity. It remains to show that $\psi$ lies over $\varphi$. The morphisms $\varphi$ and $\psi$ composed with $U_T\to X_T$ agree on every $V(I^k)$. Since $X$ is separated, their equalizer is a closed subscheme of $\A(N)_T$ containing every $V(I^k)$; locally on $T=\Spec(R)$ its ideal $J\subseteq R[N]$ is therefore contained in $\bigcap_kI^k$. Since $X^n\in I^k$ implies $w(n)\geq k$, $\bigcap_kI^k=0$, so $J=0$ and the two morphisms agree.

(ii) Write $N=\bigcup_iN_i$ as the directed union of its finitely generated submonoids. For every $i$, $Q\cap N_i$ is a face of $N_i$ (Proposition~\ref{prop:face-criterion}), and by (i), $U^{N_i}=X^{N_i}\times_{X^{Q\cap N_i}}U^{Q\cap N_i}$. The scheme $U$ is separated, and it is locally of finite presentation over $S$, being \'etale over $X$. Let $T$ be an affine scheme over $S$. Proposition~\ref{prop:colimit-magnets}, applied to $X$ and $U$, to $N=\bigcup_iN_i$ and to $Q=\bigcup_i(Q\cap N_i)$, gives $X^N(T)=\colim_iX^{N_i}(T)$, $X^Q(T)=\colim_iX^{Q\cap N_i}(T)$, and the same for $U$. Since filtered colimits of sets commute with fibre products,
\[U^N(T)=\colim_i\bigl(X^{N_i}(T)\times_{X^{Q\cap N_i}(T)}U^{Q\cap N_i}(T)\bigr)=X^N(T)\times_{X^Q(T)}U^Q(T).\]
Both sides are Zariski sheaves, so this holds for every $T$.
\end{proof}

\subsection{Fixed-point-reflecting atlases}

\begin{definition}[{\cite[Def.~6.1]{AlgMag}}]\label{def:fpr}
Let $f\colon U\to X$ be a $\D(M)_S$-equivariant morphism of algebraic spaces over $S$ and let $Z\subseteq M$ be a subgroup.
\begin{enumerate}
\item $f$ is \emph{fixed-point reflecting} (FPR) if the canonical morphism $U^{\D(M)_S}\to U\times_XX^{\D(M)_S}$ is an isomorphism.
\item $f$ is \emph{$Z$-FPR} if it is fixed-point reflecting for the induced action of $\D(M/Z)_S$, that is, if $U^{\D(M/Z)_S}\to U\times_XX^{\D(M/Z)_S}$ is an isomorphism.
\item $f$ is \emph{strongly-FPR} if it is $Z$-FPR for every subgroup $Z\subseteq M$.
\item $f$ is an \emph{equivariant atlas} if $f$ is \'etale and surjective and $U=\coprod_{\alpha\in A}U_\alpha$ is a disjoint union of $\D(M)_S$-stable $S$-affine schemes; the equivariant atlas $f$ is \emph{$S$-affine} if $A$ may be chosen finite.
\item $f$ is a \emph{$Z$-FPR atlas} if it is a $Z$-FPR equivariant atlas.
\item $f$ is a \emph{strongly-FPR atlas} if it is a strongly-FPR equivariant atlas.
\end{enumerate}
\leansrc{fpr}
\end{definition}

By Proposition~\ref{prop:subgroup-fixed}, for $M$ finitely generated, $f$ is $Z$-FPR if and only if $U^Z\to U\times_XX^Z$ is an isomorphism.

\begin{proposition}[{\cite[Rem.~6.2, Prop.~6.3, Lem.~6.7, Fact~7.1]{AlgMag}}]\label{prop:fpr-basic}
Let $Z\subseteq M$ be a subgroup.
\begin{enumerate}
\item An equivariant atlas is $S$-affine if and only if it can be written with a single $S$-affine scheme $U$.
\item If $M=\ZZ^r$, every \'etale equivariant morphism $U\to X$ is FPR.
\item Let $X'\to X$ be an equivariant affine morphism of algebraic spaces over $S$ and $U\to X$ a $Z$-FPR atlas (resp.\ an $S$-affine $Z$-FPR atlas). Then $U\times_XX'\to X'$ is a $Z$-FPR atlas (resp.\ an $S$-affine $Z$-FPR atlas).
\item Every equivariant monomorphism $U\to X$ is strongly-FPR.
\end{enumerate}
\leansrc{fprBasic}
\end{proposition}

\begin{proof}
(i), (iv) \cite[Rem.~6.2, proof of Fact~7.1]{AlgMag}; the proof of \cite[Fact~7.1]{AlgMag} is written for every equivariant monomorphism.

(iii) Put $K=\D(M/Z)_S$ and $U'=U\times_XX'$. The morphism $U'\to X'$ is \'etale and surjective, as a base change of $U\to X$, and $U'=\coprod_\alpha U_\alpha\times_XX'$, where each $U_\alpha\times_XX'$ is $\D(M)_S$-stable and $S$-affine, being affine over the $S$-affine scheme $U_\alpha$; the index set is unchanged. Fixed-point functors commute with fibre products, so $U'^K=U^K\times_{X^K}X'^K$. As $U\to X$ is $Z$-FPR, $U^K=U\times_XX^K$, hence $U'^K=(U\times_XX^K)\times_{X^K}X'^K=U\times_XX'^K=U'\times_{X'}X'^K$, and this identification is the canonical morphism. So $U'\to X'$ is a $Z$-FPR atlas, $S$-affine if $U\to X$ is.

(ii) Let $f\colon U\to X$ be \'etale and equivariant. Both $U^{\D(M)_S}$ and $U\times_XX^{\D(M)_S}$ are subfunctors of $U$, so it suffices to show that a point $u\in U(T)$ with $f(u)\in X^{\D(M)_S}(T)$ lies in $U^{\D(M)_S}(T)$. Let $G=\D(M)_T$ with its universal point $g\in\D(M)(G)$, and let $a=g\cdot u_G$ and $b=u_G$ in $U(G)$. Then $f(a)=g\cdot f(u)_G=f(u)_G=f(b)$. The morphism $f$ is \'etale, hence unramified, so its diagonal is an open immersion \cite[Tag~05W1]{Stacks}, and the preimage $W$ of the diagonal under $(a,b)\colon G\to U\times_XU$ is an open subscheme of $G$; a morphism $r\colon Y\to G$ factors through $W$ if and only if $a\circ r=b\circ r$. Let $x$ be a point of $G$, with residue field $K$, over $\tau\colon\Spec(K)\to T$, and let $g_0\in\D(M)(K)$ be the corresponding element. Let $q_1,q_2\colon\D(M)_K\to G$ be the morphisms over $\tau$ given by $h\mapsto h$ and $h\mapsto g_0h^{-1}$. The unit section $e_T\colon T\to G$ factors through $W$, since $a\circ e_T=u=b\circ e_T$, and $q_1(1)=e_T\circ\tau=q_2(g_0)$; so the open subschemes $q_1^{-1}(W)$ and $q_2^{-1}(W)$ of $\D(M)_K$ contain the unit and $g_0$ respectively. Since $M\cong\ZZ^r$, $\D(M)_K=\Spec(K[t_1^{\pm1},\dots,t_r^{\pm1}])$ is irreducible, so they meet in a point with residue field $L\supseteq K$, corresponding to $h\in\D(M)(L)$ with $h\cdot u_L=u_L$ and $g_0h^{-1}\cdot u_L=u_L$. Hence $g_0\cdot u_L=u_L$, so $x\in W$, since the morphism $\Spec(L)\to G$ given by $g_0$ over $\tau$, with image $x$, factors through the open $W$. Therefore $W=G$, $a=b$, that is, $g\cdot u_G=u_G$, and $u\in U^{\D(M)_S}(T)$.
\end{proof}

\begin{definition}[{\cite[Def.~7.2]{AlgMag}}]\label{def:zll}
The action of $\D(M)_S$ on $X$ is \emph{Zariski locally linearizable} if $X$ is covered by $\D(M)_S$-stable open subspaces that are $S$-affine schemes.
\leansrc{zll}
\end{definition}

\begin{proposition}[{\cite[Fact~7.3]{AlgMag}}]\label{prop:zll-atlas}
If the action of $\D(M)_S$ on $X$ is Zariski locally linearizable, then $X$ has a strongly-FPR atlas.
\leansrc{zllAtlas}
\end{proposition}

\begin{proof}
\cite[proof of Fact~7.3]{AlgMag}
\end{proof}

\begin{lemma}[{\cite[Cor.~10.2]{AHR}}]\label{lem:chart}
Assume that $S$ is affine and of finite type over $\ZZ$ and that $M$ is finitely generated. Let $X$ be a separated algebraic space over $S$ of finite type with an action of $\D(M)_S$, let $k$ be a field and let $x\colon\Spec(k)\to X$ be a closed immersion. Then there exist an affine scheme $U$ with an action of $\D(M)_S$, an equivariant \'etale morphism $f\colon U\to X$ and a morphism $u\colon\Spec(k)\to U$ with $f\circ u=x$ such that, for every $k$-scheme $T$ and every $g\in\D(M)(T)$, $g\cdot f(u_T)=f(u_T)$ implies $g\cdot u_T=u_T$.
\leansrc{chart}
\end{lemma}

So $u$ has residue field $\kappa(x)$, and $u$ and $x$ have the same stabilizer, the reverse implication being a consequence of equivariance.

\begin{proof}
The lemma is the case of the diagonalizable group $\D(M)_S$ of \cite[Cor.~10.2]{AHR}. The proof below follows the architecture of the proof of \cite[Thm.~1.1]{AHR}, with algebraic stacks, Tannaka duality and the cotangent complex replaced by graded commutative algebra: by Proposition~\ref{prop:gradings}(ii), an action of $\D(M)_S$ on an $S$-affine scheme $\Spec_S(\cA)$ is an $M$-grading of $\cA$, equivariant morphisms of such schemes are graded homomorphisms, and the $\D(M)_S$-stable closed subschemes are those with homogeneous ideals. Rings of $S$-affine schemes with an action are graded in this way, and their stable ideals are homogeneous.

\emph{Step 1: the point and its orbit.} The morphism $\Spec(k)\to\Spec(\ZZ)$ is of finite type, so $k$ is finite \cite[Tags~02J6, 01TB]{Stacks}. Let $q=|k|$ and let $c$ generate $k^\times$. The $S$-scheme $Z'=\Spec\bigl(\cO(S)[t]/(t^q-t)\bigr)_{qt^{q-1}-1}$ is \'etale over $S$, and $t\mapsto c$ defines a monomorphism $\Spec(k)\to Z'$ over $S$. The projection $X\times_SZ'\to X$, with the action on the first factor, is equivariant and \'etale, and $(x,c)$ is a closed point of $X\times_SZ'$ with residue field $k$ over $x$; a chart for $(x,c)$ composed with this projection is a chart for $x$ with the required property. So we may assume that $X$ has an invariant morphism to an $S$-scheme $Z'$ through which $x$ becomes a monomorphism $\Spec(k)\to Z'$. The stabilizer $H\subseteq\D(M)_k$ of $x$, the fibre product of $g\mapsto(g\cdot x,x)$ with the diagonal of $X$, is a closed subgroup scheme since $X$ is separated, so $H=\D(M/Y)_k$ for a subgroup $Y\subseteq M$ \cite[Exp.~IX, Prop.~8.1]{SGA3}. The orbit morphism $o\colon O=\D(Y)_k\to X$, $g\mapsto g\cdot x$, is then a monomorphism, locally of finite presentation. It is flat on a nonempty open subscheme of $O$ (generic flatness of a finite type epimorphism from a field onto a noetherian ring), hence flat everywhere, since a $\D(M)_k$-stable closed subset of $O$ is empty or $O$. A flat monomorphism locally of finite presentation is an open immersion \cite[Tag~025G]{Stacks}, applied \'etale locally on $X$; so $o$ is an open immersion into its scheme-theoretic image, which is $\D(M)_S$-stable. Removing the closed $\D(M)_S$-stable complement of $|O|$ in this image, we may assume that $o$ is a closed immersion.

\emph{Step 2: thickenings and graded completion.} Let $I$ be the ideal of $O$ in $X$, which is stable since $\D(M)_S$ is flat, and $X_n=V(I^{n+1})$. As $O$ is affine, $X_n$ is an affine scheme \cite[Tags~05ZR, 06AD]{Stacks}; let $O_n=\Gamma(X_n)$. Put $C=\bigoplus_{d\in M}\lim_n(O_n)_d$ and $\mathfrak m=\ker(C\to k[Y])$. The ring $\bigoplus_iI^i/I^{i+1}$ is a finitely generated $k$-algebra graded by $\NN\times M$; by Gordan's lemma its homogeneous component of degree $0$ for $M$ is a finitely generated algebra, over which every homogeneous component is finite. Hence $C$ is noetherian, $C/\mathfrak m^{n+1}=O_n$, and $(C,\mathfrak m)$ is a \emph{complete pair}: $C_d\to\lim_n(C/\mathfrak m^n)_d$ is bijective for every $d$. For a complete pair $(C,\mathfrak m)$ with $C$ noetherian, $\bigcap_n\mathfrak m^nN=0$ for every finitely generated graded $C$-module $N$, $(C/\mathfrak b,\mathfrak m)$ is a complete pair for every homogeneous ideal $\mathfrak b$, and every compatible system of finitely generated graded $C/\mathfrak m^n$-modules is the system of reductions of a finitely generated graded $C$-module, unique up to unique isomorphism.

\emph{Step 3: equivariant effectivity.} Let $X'$ be a separated algebraic space over $S$ of finite type with an action of $\D(M)_S$ and $(C,\mathfrak m)$ a complete pair with $C$ noetherian. We show that every compatible system of equivariant morphisms $\varphi_n\colon\Spec(C/\mathfrak m^n)\to X'$ is induced by a unique equivariant morphism $\varphi\colon T=\Spec(C)\to X'$. For uniqueness, the equalizer of two such morphisms is a closed subscheme of $T$, as $X'$ is separated, whose ideal is homogeneous and contained in $\bigcap_n\mathfrak m^n=0$. Existence is proved by noetherian induction on the $\D(M)_S$-stable closed subspaces of $X'$, with the following three ingredients.

(a) \emph{Affine stratum.} If $X'\neq\emptyset$, it has a nonempty $\D(M)_S$-stable open subspace which is an affine scheme. Let $n\geq1$ kill the torsion of $M$. If $n$ is invertible at some point of $X'$, put $N=n$ and $K=0$; otherwise every residue characteristic of $X'$ is a prime dividing $n$; choose one such prime $\ell$, and put $N$ the prime-to-$\ell$ part of $n$ and $K$ the $\ell$-primary torsion subgroup of $M$. Then the open subspace $X'[1/N]$ is nonempty and its residue fields have characteristic $\ell$, $K$ acts trivially on its points with values in fields, and $\D(M/K)_S$ preserves the reduced closed subspace $X''$ of $X'[1/N]$, since $R[M/K]$ is reduced for every reduced ring $R$ in which $N$ is invertible. Let $V\to X''$ be an \'etale morphism with $V$ an integral affine noetherian scheme over an affine open subscheme of $S$. By the linear independence of characters over the fraction field of $V$, there are a nonempty open $D(s)\subseteq V$ and a subgroup $Z\subseteq M/K$ such that, over $D(s)$, a point of $\D(M/K)$ fixes the point $V\to X''$ if and only if it is trivial on $Z$. The image $W$ of the action morphism $\D(M/K)_S\times_SD(s)\to X''$ is a stable open subspace, on which the action induces a free action of $\D(Z)_S$; the quotient $W/\D(Z)_S$ is an algebraic space \cite[Tag~04S6]{Stacks}, it is quasi-separated, hence has a dense open subspace which is a scheme \cite[Tag~086U]{Stacks}, and the morphism $W\to W/\D(Z)_S$ is affine, being fppf locally the projection $\D(Z)\times V'\to V'$ \cite[Tag~041Z]{Stacks}. The preimage in $W$ of a nonempty affine open subscheme of the quotient is a nonempty stable affine open subscheme of $X''$. The open subspace of $X'$ with the same points is then affine \cite[Tag~07VT]{Stacks}, and it is $\D(M)_S$-stable since $K$ acts trivially on points with values in fields.

(b) \emph{The two strata.} Let $U_0$ be as in (a), $J$ the stable ideal of a closed subspace structure on $|X'|\setminus|U_0|$, $\mathfrak a\subseteq C$ the homogeneous ideal $\bigoplus_d\lim_n(J\cdot C/\mathfrak m^n)_d$, and $T_U=T\setminus V(\mathfrak a)$. For $h\in\Gamma(U_0)$ of degree $d$ and $g\in\mathfrak a$ homogeneous of degree $\delta$, there are $N\geq0$ and $t\in C_{d+N\delta}$ with $t\equiv g^Nh$ modulo every $\mathfrak m^n$, computed on an affine \'etale chart of $X'$ from the sections of $J^N$ which restrict to multiples of $h$ on $U_0$; the elements $t/g^N\in\Gamma(D(g))$ glue to a section of $\cO(T_U)$ of degree $d$, and $h\mapsto t/g^N$ is a graded homomorphism, which defines an equivariant morphism $T_U\to U_0$. By the induction hypothesis applied to the closed subspaces $V(J^k)$ and the complete pairs $(C/\mathfrak a^k,\mathfrak m)$, there are compatible equivariant morphisms $\Spec(C/\mathfrak a^k)\to V(J^k)$.

(c) \emph{Gluing.} Let $\hat C$ be the $\mathfrak a$-adic completion of $C$. The morphisms $\Spec(C/\mathfrak a^k)\to X'$ are induced by a morphism $\theta\colon\Spec(\hat C)\to X'$ by Grothendieck's algebraization theorem for the separated target $X'$ \cite[Tag~0GHK]{Stacks}, and $\theta$ agrees with $T_U\to U_0$ on $\Spec(\hat C)\setminus V(\mathfrak a\hat C)$, both being given there by the formula of (b). As $T_U\sqcup\Spec(\hat C)\to T$ is an fpqc covering and $X'$ is an fpqc sheaf \cite[Tag~0APL]{Stacks}, it remains to see that the two pullbacks of $\theta$ to $\Spec(\hat C\otimes_C\hat C)$ agree. Their equalizer is closed, as $X'$ is separated, and contains the part over $T_U$ and every $V(\mathfrak a^k)$. The following separating lemma shows that it is everything: if $R$ is a flat $C$-algebra and $r\in R$ satisfies $\mathfrak a^sr=0$ and $r\in\bigcap_k\mathfrak a^kR$, then $r=0$; indeed flatness gives $(0:\mathfrak a^s)_R=(0:\mathfrak a^s)_CR$ and $\mathfrak a^kR\cap(0:\mathfrak a^s)_CR=(\mathfrak a^k\cap(0:\mathfrak a^s)_C)R$, and $\mathfrak a^k\cap(0:\mathfrak a^s)_C=0$ for $k$ large by the Artin--Rees lemma. This gives $\varphi\colon T\to X'$. The same lemma with the flat $C$-algebra $R=\ZZ[M]\otimes C$ shows that $\varphi$ is equivariant, and $\varphi$ restricts to $\varphi_n$ by construction of $\theta$.

\emph{Step 4: approximation.} Step~3 applied to $X$ and to the closed immersions $X_n\to X$ gives an equivariant morphism $\varphi\colon T=\Spec(C)\to X$ restricting to $X_n\to X$ on $\Spec(C/\mathfrak m^{n+1})$. Choose homogeneous $x_1,\dots,x_r\in C$ such that $\ZZ[x_1,\dots,x_r]\to C/\mathfrak m^n$ is surjective for every $n$, and let $P=\ZZ[x_1,\dots,x_r]$, graded by the degrees of the $x_i$. By Gordan's lemma, $P_0$ is of finite type over $\ZZ$ and every $P_d$ is finite over $P_0$; let $\mathfrak p_0=\ker(P_0\to k)$ and $Q=P\otimes_{P_0}(P_0)^\wedge_{\mathfrak p_0}$. Then $Q\to C$ is surjective, with kernel generated by homogeneous $f_1,\dots,f_s$; let $(r_j)$ be homogeneous generators of the module of relations among the $f_i$. Since $X$ is locally of finite presentation, the data $(f,r,\varphi)$ is a solution over $(P_0)^\wedge_{\mathfrak p_0}$ of a system of polynomial equations over $P_0$. The local rings of $P_0$ are G-rings \cite[Tag~07PX]{Stacks}, so by Artin approximation \cite[Tag~0CAR]{Stacks}, for every $\nu\geq1$ there are an \'etale ring map $P_0\to P_0'$, a prime $\mathfrak p'$ over $\mathfrak p_0$ with $\kappa(\mathfrak p')=k$ and a solution $(\tilde f,\tilde r,\tilde\varphi)$ over $P_0'$ congruent to $(f,r,\varphi)$ modulo $\mathfrak p'^{\nu}$. Put $W=\Spec\bigl((P\otimes_{P_0}P_0')/(\tilde f)\bigr)$, an affine scheme with an action of $\D(M)_S$ and an equivariant morphism $\tilde\varphi\colon W\to X$, and let $w\in W$ be the point given by the unit of $O$; it lies on a closed stable orbit $O_W\subseteq W$ mapped isomorphically onto $O$. For $\nu$ larger than the Artin--Rees constants of the ideal $(f)$ and of the module of relations, the initial ideal of $(\tilde f)$ is contained in that of $(f)$ in the associated graded ring of the completion of $Q$ at $w$, so $\cO_{W,w}^\wedge$ has length at least that of $\cO_{T,w}^\wedge$ modulo every power of the maximal ideal. As $\cO_{X,x}^\wedge\to\cO_{T,w}^\wedge$ is an isomorphism (the $X_n$ are the thickenings of $O$) and $\cO_{X,x}^\wedge\to\cO_{W,w}^\wedge$ is surjective (it is so modulo $\mathfrak m_w^2$), the latter is an isomorphism; here $\cO_{X,x}$ is the local ring at $x$ of an affine \'etale neighbourhood of $x$ with residue field $k$ \cite[Tag~0BBP]{Stacks}. Hence $\tilde\varphi$ is \'etale at $w$ \cite[Tag~039M]{Stacks}.

\emph{Step 5: shrinking.} The \'etale locus $\Omega$ of $\tilde\varphi$ is open and $\D(M)_S$-stable and contains $O_W$. Let $\mathfrak J$ be the radical ideal of $W\setminus\Omega$. For $a\in\mathfrak J$, the image of $a$ under the coaction is nilpotent modulo $\mathfrak J\otimes\ZZ[M]$, so all homogeneous components of a power of $a$ lie in $\mathfrak J$; these components, for $a$ running over $\mathfrak J$, generate a homogeneous ideal $\mathfrak J'$ with $V(\mathfrak J')=W\setminus\Omega$. As $O_W\cap V(\mathfrak J')=\emptyset$, the ideal of $O_W$ and $\mathfrak J'$ generate $\Gamma(W)$, and taking components of degree $0$, $1=i+h$ with $i$ in the ideal of $O_W$ and $h\in\mathfrak J'_0$. Then $U=D(h)$ is a $\D(M)_S$-stable affine open subscheme of $W$ containing $O_W$ and contained in $\Omega$. Let $f=\tilde\varphi|_U$ and let $u\colon\Spec(k)\to U$ be the point $w$; then $f\circ u=x$. If $T$ is a $k$-scheme and $g\in\D(M)(T)$ fixes $f(u_T)=x_T$, then $g\in H(T)$, and $g\cdot u_T=u_T$ since $O_W\to O$ is an isomorphism and $H$ is the stabilizer of the unit of $O=\D(M)_k/H$.
\end{proof}

\begin{theorem}\label{thm:strong-atlas}
Assume that $S$ is separated and that $M$ is finitely generated. Let $X$ be a separated algebraic space over $S$ locally of finite presentation with an action of $\D(M)_S$. Then $X$ has a strongly-FPR atlas; if $|X|$ is quasi-compact, it has an $S$-affine strongly-FPR atlas.
\leansrc{strongAtlas}
\end{theorem}

\begin{proof}
For a subgroup $Z\subseteq M$ put $K_Z=\D(M/Z)_S$, so that $X^Z=X^{K_Z}$ (Proposition~\ref{prop:subgroup-fixed}).

\emph{Step 1: the bad locus.} Let $f\colon U\to X$ be an equivariant \'etale morphism with $U$ an affine scheme. Let $E\subseteq\D(M)_U$ be the fibre product of $(g,v)\mapsto(g\cdot f(v),f(v))$, $\D(M)_U\to X\times_SX$, with the diagonal of $X$; it is a closed subscheme since $X$ is separated. The morphisms $E\to U$, $(g,v)\mapsto g\cdot v$ and $(g,v)\mapsto v$, have the same composite with $f$, so their equalizer $O\subseteq E$ is the preimage of the diagonal of $U$ over $X$, which is an open immersion as $f$ is unramified \cite[Tag~05W1]{Stacks}; $O$ is also closed in $E$, as $U$ is separated. Let $B\subseteq U$ be the image of $E\setminus O$ under the projection. If $T$ is a scheme, $v\in U(T)$ has image in $U\setminus B$ and $g\in\D(M)(T)$ satisfies $g\cdot f(v)=f(v)$, then $(g,v)\colon T\to E$ has image in $O$, so $g\cdot v=v$. Consequently, if $W\subseteq U\setminus B$ is a $\D(M)_S$-stable open subscheme, then $W\to X$ is strongly-FPR: for a subgroup $Z$, $W^{K_Z}$ and $W\times_XX^{K_Z}$ are subfunctors of $W$, and a point $v\in W(T)$ with $f(v)$ fixed by $K_Z$ is fixed by $K_Z$, since every $g\in K_Z(T')\subseteq\D(M)(T')$, $T'\to T$, fixes $f(v_{T'})$. By equivariance, $B$ is $\D(M)_S$-stable.

In Steps 2 to 4, $S$ is affine of finite type over $\ZZ$ and $X$ is of finite type over $S$, so that every affine scheme \'etale over $X$ is noetherian.

\emph{Step 2: $B$ is closed.} Let $U=\Spec(A)$ in Step~1. Then $E=V(I)$ for an ideal $I\subseteq A[M]$ stable under the endomorphisms induced by $g\mapsto g^n$, $n\geq0$ ($E$ is a subgroup scheme of $\D(M)_U$), and $E\setminus O=V(I)\cap D(e)$ with $e$ idempotent modulo $I$. The set $B$ is constructible by Chevalley's theorem \cite[Tag~054K]{Stacks}. For a subgroup $Y\subseteq M$, let $\mathfrak c_Y\subseteq A$ be the ideal generated by the coefficients of the images in $A[M/Y]$ of the elements of $I$, so that $\D(M/Y)_{V(\mathfrak c_Y)}\subseteq E$, and let $C_Y\subseteq U$ be the image of $(E\setminus O)\cap\D(M/Y)_{V(\mathfrak c_Y)}$. Then $C_Y$ is closed: it is the image in $\Spec(A/\mathfrak c_Y)$ of $D(\bar e)\subseteq\Spec((A/\mathfrak c_Y)[M/Y])$, where $\bar e$ is the image of $e$, an idempotent; this image is the complement of $V(\mathfrak c)$, for $\mathfrak c$ the ideal generated by the coefficients of $\bar e$, and $\mathfrak c=\mathfrak c^2$ as $\bar e=\bar e^2$, so $\mathfrak c$ is generated by an idempotent and $V(\mathfrak c)$ is open and closed \cite[Tag~00EH]{Stacks}. Every point of $B$ lies in some $C_Y$: let $(g,v)$ be a point of $E\setminus O$ over $v\in B$, with values in a field $L$, and $Y=\{m\in M:\chi_m(g)=1\}$. For $a\in I$ with image $\sum_{c\in M/Y}a_c\chi_c$ in $A[M/Y]$, the points $(g^n,v)$, $n\geq0$, lie in $E$, so $\sum_ca_c(v)\chi_c(g)^n=0$ for all $n$; the $\chi_c(g)$, $c\in M/Y$, are pairwise distinct, so $a_c(v)=0$ for all $c$ by the invertibility of Vandermonde matrices. Hence $v\in V(\mathfrak c_Y)$, $(g,v)\in\D(M/Y)_{V(\mathfrak c_Y)}$, and $v\in C_Y$. Since $C_Y\subseteq B$ for all $Y$, $B=\bigcup_YC_Y$ is stable under specialization, hence closed \cite[Tag~0542]{Stacks}.

\emph{Step 3: the good locus.} Let $v\in U\setminus B$. Let $\mathfrak J\subseteq A$ be the radical ideal of $B$ and $b\in\mathfrak J$ with $b(v)\neq0$. As $B$ is $\D(M)_S$-stable, the image of $b$ under the coaction is nilpotent modulo $\mathfrak J\otimes\ZZ[M]$, so the homogeneous components of $b^e$ lie in $\mathfrak J$ for some $e\geq1$. Their sum $b^e$ does not vanish at $v$, so some homogeneous component $a$ of $b^e$ satisfies $v\in D(a)\subseteq U\setminus B$. For $a$ homogeneous of degree $m$ the coaction maps $a$ to $a\otimes X^m$, so $D(a)$ is a $\D(M)_S$-stable affine open subscheme of $U$, and $D(a)\to X$ is strongly-FPR by Step~1.

\emph{Step 4: $S$ of finite type over $\ZZ$.} Let $x\in|X|$ be a closed point, represented by a closed immersion $\Spec(k)\to X$. Lemma~\ref{lem:chart} gives an affine scheme $U$, an equivariant \'etale morphism $f\colon U\to X$ and $u\colon\Spec(k)\to U$ over $x$ with the stabilizer property; let $v$ be the image of $u$. Then $v\notin B$: a point of $E$ over $v$, with values in a field $L$, is $(g,u_L)$ with $g\in\D(M)(L)$ fixing $x_L$, since $\kappa(v)=k$ as $f\circ u=x$; so $g\cdot u_L=u_L$ and the point lies in $O$. By Step~3 there is a homogeneous $a\in\Gamma(U)$ with $v\in D(a)$ and $D(a)\to X$ strongly-FPR; its image is an open neighbourhood of $x$. The space $|X|$ is Jacobson \cite[Tags~02J6, 0BA3, 0BA6]{Stacks}, $X$ being decent as it is separated, so an open subset of $|X|$ containing all closed points is $|X|$. As $X$ is quasi-compact, finitely many $D(a_1),\dots,D(a_n)$ cover $X$, and $\coprod_iD(a_i)\to X$ is a strongly-FPR atlas by an affine scheme.

\emph{Step 5: $S$ affine and $X$ quasi-compact.} Then $X$ is of finite presentation over $S$. Write $S$ as a limit of affine schemes $S_\lambda$ of finite type over $\ZZ$; $X$ with its action descends to a separated algebraic space $X_\lambda$ of finite presentation over some $S_\lambda$ with an action of $\D(M)_{S_\lambda}$: the algebraic space $X$, the action morphism and the identities expressing the action axioms descend by \cite[Tag~07SK]{Stacks}, since $\D(M)_S\times_SX$ is of finite presentation over $S$, and separatedness descends by \cite[Tag~0851]{Stacks}. The base change $U\to X$ of the atlas $U_\lambda\to X_\lambda$ of Step~4 is an equivariant \'etale surjection from an affine scheme, and it is strongly-FPR since fixed-point functors commute with base change, so that the isomorphisms $U_\lambda^Z\to U_\lambda\times_{X_\lambda}X_\lambda^Z$ base change to $U^Z\to U\times_XX^Z$.

\emph{Step 6: general case.} Cover $S$ by affine open subschemes $S_i$ and each $X_{S_i}$ by quasi-compact open subspaces $V$. The image $\D(M)\cdot V$ of $\D(M)_S\times_SV$ under the action is open, as the image of the flat projection composed with the automorphism $(g,y)\mapsto(g,g\cdot y)$; it is quasi-compact and $\D(M)_S$-stable. Step~5 applied to each $\D(M)\cdot V$ over $S_i$ gives strongly-FPR equivariant \'etale surjections from affine schemes. The open immersions $\D(M)\cdot V\to X$ are strongly-FPR (Proposition~\ref{prop:fpr-basic}(iv)), so the disjoint union of these morphisms is a strongly-FPR equivariant \'etale surjection onto $X$. Its pieces are affine schemes, hence $S$-affine since $S$ is separated \cite[Tag~01SG]{Stacks}, so it is a strongly-FPR atlas. If $|X|$ is quasi-compact, finitely many pieces cover $X$, which gives an $S$-affine strongly-FPR atlas.
\end{proof}

\begin{theorem}[{\cite[Thm.~6.8]{AlgMag}}]\label{thm:atlas}
Assume that $S$ is separated and that $M$ is finitely generated. Let $X$ be a separated algebraic space over $S$ locally of finite presentation with an action of $\D(M)_S$, and let $Z\subseteq M$ be a subgroup. Then there exists a $Z$-FPR atlas $U\to X$; it may be chosen $S$-affine if the topological space $|X|$ \cite[Tag~03BY]{Stacks} is quasi-compact.
\leansrc{atlas}
\end{theorem}

\begin{proof}
By Theorem~\ref{thm:strong-atlas}, $X$ has a strongly-FPR atlas, $S$-affine if $|X|$ is quasi-compact. A strongly-FPR atlas is a $Z$-FPR atlas for every subgroup $Z\subseteq M$ (Definition~\ref{def:fpr}).
\end{proof}

\begin{remark}\label{rem:conj69}
Theorem~\ref{thm:strong-atlas} proves \cite[Conj.~6.9(i)]{AlgMag} when $S$ is separated, and \cite[Conj.~6.9(ii)]{AlgMag} when $S$ is separated and quasi-compact: if $S$ and $X$ are quasi-compact, then $|X|$ is quasi-compact.
\leansrc{conj69}
The hypothesis that $S$ is separated cannot be dropped. Let $S$ be the plane with doubled origin over a field, with affine charts $A_1,A_2$, let $X=A_1$ and $M=0$. Then $X$ is separated, quasi-compact and locally of finite presentation over $S$. If $U\to X$ is \'etale with $U$ affine over $S$ and with nonempty fibre $E$ over the origin, then $U=U\times_SA_1$ is affine and regular, and $U\setminus E=U\times_SA_2$ is affine; this is impossible since $E$ has codimension $2$ in $U$ \cite[Tag~0BCU]{Stacks}. Hence $X$ has no strongly-FPR atlas.
\leansrc{conj69DoubledPlane}
\end{remark}

\subsection{Representability}

\begin{proposition}[{\cite[Prop.~8.1]{AlgMag}}]\label{prop:fixed-closed}
Let $Z\subseteq M$ be a subgroup. If $X$ has a $Z$-FPR atlas, then $X^Z=X^{\D(M/Z)_S}$ is represented by a closed subspace of $X$.
\leansrc{fixedClosed}
\end{proposition}

\begin{proof}
Let $f\colon U=\coprod_aU_a\to X$ be a $Z$-FPR atlas, let $K=\D(M/Z)_S$ and let $\pi\colon\D(M)_S\to\D(Z)_S$ be the restriction, with kernel $K$. As $\ZZ[M]$ is free over $\ZZ[Z]$, $\pi$ is affine and faithfully flat, so its base changes $\pi_T$ are fpqc coverings. For $\varphi\in X^Z(T)$ and $k\in K$, $k\cdot\varphi(1)=\varphi(\pi(k))=\varphi(1)$, so evaluation at the unit is a map $X^Z\to X^K$.

\emph{Injectivity.} Let $\varphi,\varphi'\in X^Z(T)$ with $\varphi(1)=\varphi'(1)=x$. By equivariance, $\varphi\circ\pi_T$ and $\varphi'\circ\pi_T$ are both $h\mapsto h\cdot x$, so $\varphi$ and $\varphi'$ agree after the fpqc covering $\pi_T$. An algebraic space satisfies the sheaf property for the fpqc topology \cite[Tag~0APL]{Stacks}, so $\varphi=\varphi'$.

\emph{Surjectivity.} Let $x\in X^K(T)$. For each $a$, let $T_a=T\times_{x,X}U_a$ and let $u_a\in U_a(T_a)$ be the tautological lift of $x_{T_a}$; the $T_a\to T$ form an \'etale covering. As $f$ is $Z$-FPR, $U^K=U\times_XX^K$, so $u_a$ is fixed by $K$. The orbit morphism $o_a\colon\D(M)_{T_a}\to U_a$, $h\mapsto h\cdot u_a$, then satisfies $o_a(hk)=o_a(h)$ for $k\in K$, that is, it coequalizes the two projections of $\D(M)_{T_a}\times_{\D(Z)_{T_a}}\D(M)_{T_a}$. Since $U_a$ is a scheme, $o_a$ descends along the fpqc covering $\pi_{T_a}$ to a morphism $\psi_a\colon\D(Z)_{T_a}\to U_a$ \cite[Tag~023Q]{Stacks}, which is equivariant since $o_a$ is, with $\psi_a(1)=u_a$. So $f\circ\psi_a\in X^Z(T_a)$ has value $x_{T_a}$ at the unit. By injectivity over $T_a\times_TT_b$, these points agree on overlaps, and since $X^Z$ is an fppf sheaf (Proposition~\ref{prop:attractor-sheaf}) they glue to $\varphi\in X^Z(T)$ with $\varphi(1)=x$.

\emph{Closedness.} For an $S$-affine piece $U_a=\Spec_S(\cA_a)$, $U_a^K$ is the closed subscheme cut out by the homogeneous components of degree $m\notin Z$. As $f$ is $Z$-FPR, the base change of $X^K\to X$ along $f$ is the closed immersion $U^K\to U$. For a scheme $T$ with a morphism $T\to X$, the fppf sheaf $T\times_XX^K$ therefore becomes the closed subscheme $(T\times_XU)\times_UU^K$ of $T\times_XU$ after the \'etale covering $T\times_XU\to T$. This closed subscheme carries the descent datum of a subsheaf, which is effective \cite[Tag~03I0]{Stacks}, so $T\times_XX^K$ is a closed subscheme of $T$ \cite[Tag~02W5]{Stacks}. Hence $X^Z=X^K$ is represented by a closed subspace of $X$.
\end{proof}

\begin{theorem}[{\cite[Thm.~8.3]{AlgMag}}]\label{thm:representability}
Assume that $M$ is finitely generated. Let $N$ be a magnet and assume that $X$ has a $\units N$-FPR atlas. Assume that $X$ is separated and that one of the following conditions holds:
\begin{enumerate}[label=\textup{(\alph*)}]
\item $N$ is finitely generated;
\item $X$ is locally of finite presentation.
\end{enumerate}
Then $X^N$ is represented by an algebraic space over $S$, $X^{\units N}$ is a closed subspace of $X$, and:
\begin{enumerate}
\item the morphism $p_{N,\units N}\colon X^N\to X^{\units N}$ is affine;
\item $X^N$ is separated, in particular quasi-separated;
\item if $X$ is locally of finite presentation, so is $X^N$;
\item if $|X|$ is quasi-compact, so is $|X^N|$;
\item if $X$ is of finite presentation, so is $X^N$;
\item if $X$ is a scheme, so is $X^N$.
\end{enumerate}
\leansrc{representability}
\end{theorem}

\begin{proof}
\cite[proof of Thm.~8.3]{AlgMag}, with Proposition~\ref{prop:etale-face} in place of \cite[Prop.~5.1]{AlgMag}; the assertions (ii) and (iii) of \cite[Thm.~8.3]{AlgMag} are (ii) here.
\end{proof}

\begin{corollary}\label{cor:representability}
Assume that $S$ is separated and $M$ is finitely generated. Let $X$ be a separated algebraic space over $S$ locally of finite presentation with an action of $\D(M)_S$. Then for every magnet $N$, the attractor $X^N$ is a separated algebraic space locally of finite presentation over $S$, and $|X^N|$ is quasi-compact if $|X|$ is quasi-compact.
\leansrc{representabilityFree}
\end{corollary}

\begin{proof}
Combine Theorem~\ref{thm:atlas}, applied to $Z=\units N$, with Theorem~\ref{thm:representability}(b).
\end{proof}

\begin{proposition}[{\cite[Prop.~8.5]{AlgMag}}]\label{prop:finite-magnet}
Assume that $M$ is finitely generated and that $S$ is quasi-compact. Let $X$ be a separated algebraic space over $S$, quasi-compact and locally of finite presentation, with an action of $\D(M)_S$, and let $N$ be a magnet. Then there exists a finitely generated submonoid $N_c\subseteq N$ with $\units{N_c}=\units N$ and $X^{N_c}=X^N$.
\leansrc{finiteMagnet}
\end{proposition}

\begin{proof}
Let $S=V_1\cup\dots\cup V_r$ be a finite cover by affine open subschemes. For each $j$, \cite[proof of Prop.~8.5]{AlgMag} applies to $X_{V_j}$ over $V_j$, where the $V_j$-affine $\units N$-FPR atlas is given by Theorem~\ref{thm:atlas}, $V_j$ being separated and $|X_{V_j}|$ quasi-compact, and Proposition~\ref{prop:colimit-magnets} is applied to affine schemes $T$, which suffices since the attractors are sheaves; it gives a finitely generated submonoid $N_j\subseteq N$ with $(X_{V_j})^{N_j}=(X_{V_j})^N$. Let $N_c$ be the submonoid generated by $N_1,\dots,N_r$ and $\units N$; it is finitely generated since $M$ is, and $\units{N_c}=\units N$. For every scheme $T$ over $S$, the map $X^{N_c}(T)\to X^N(T)$ is injective (Proposition~\ref{prop:separated-mono}), and it is surjective: a point of $X^N(T)$ comes, over each $T\times_SV_j$, from $X^{N_j}\subseteq X^{N_c}$ (Proposition~\ref{prop:attractor-limits}(iii)), and these local points glue, by injectivity and since $X^{N_c}$ is a sheaf.
\end{proof}

\subsection{Products of attractors}

Let $M$ be a finitely generated abelian group and let $\cN=(N_\alpha)_{\alpha\in A}$ be a nonempty family of fine submonoids of $M$ with a common group $N_\alpha^{\gp}=Z$. For every $\alpha$, the morphism $\D(Z)\to\A(N_\alpha)$ is an open immersion.

\begin{proposition}[{\cite[Prop.~3.1, Prop.~3.3]{XProducts}}]\label{prop:glued-scheme}
The diagram of open immersions $\D(Z)\to\A(N_\alpha)$, $\alpha\in A$, has a colimit $\A(\cN)$ in the category of schemes, obtained by gluing the $\A(N_\alpha)$ along $\D(Z)$; each $\A(N_\alpha)\to\A(\cN)$ is an open immersion. The actions of $\D(M)$ on the $\A(N_\alpha)$ glue to an action on $\A(\cN)$, and $\A(\cN)$ is the colimit of the same diagram in the category of schemes with a $\D(M)$-action.
\leansrc{gluedScheme}
\end{proposition}

\begin{proof}
\cite[proofs of Prop.~3.1 and Prop.~3.3]{XProducts}
\end{proof}

The monoschemes $\operatorname{spec}(N_\alpha)$ of \cite{Deitmar} glue along $\operatorname{spec}(Z)$ to a monoscheme $\cN_{\mathbb F_1}$ with $\cN_{\mathbb F_1}\otimes_{\mathbb F_1}\ZZ=\A(\cN)$ \cite[Prop.~3.2]{XProducts}.
\leansrc{monoscheme}

\begin{definition}[{\cite[Def.~4.2]{XProducts}}]\label{def:attractor-family}
Let $X$ be an algebraic space over $S$ with an action of $\D(M)_S$. The \emph{attractor} $X^{\cN}$ is the functor $T\mapsto\Homeq{\D(M)_T}_T(\A(\cN)_T,X_T)$.
\leansrc{attractorFamily}
\end{definition}

\begin{theorem}[{\cite[Thm.~4.1]{XProducts}}]\label{thm:x-product}
Let $X$ be an algebraic space over $S$ with an action of $\D(M)_S$. Restriction along the $\A(N_\alpha)\to\A(\cN)$ gives an isomorphism, compatible with the morphisms to $X$,
\[X^{\cN}\xrightarrow{\ \sim\ }\textstyle\prod_{X,\,\alpha\in A}X^{N_\alpha},\]
the fibre product over $X$ of the $X^{N_\alpha}\to X$.
\leansrc{xProduct}
\end{theorem}

\begin{proof}
\cite[proof of Thm.~4.1]{XProducts}. The hypotheses 1 and 2 of \cite[Thm.~4.1]{XProducts} are used there only to show that $X^Z\to X$ is a monomorphism. By Proposition~\ref{prop:subgroup-fixed}, $X^Z$ is identified with the fixed-point functor $X^{\D(M/Z)_S}$, a subfunctor of $X$; so $X^Z\to X$ is a monomorphism for every algebraic space $X$ over $S$.
\end{proof}

\begin{proposition}[{\cite[Prop.~4.3, Prop.~4.5]{XProducts}}]\label{prop:x-product-rep}
Let $X$ be an algebraic space over $S$ with an action of $\D(M)_S$.
\begin{enumerate}
\item If $X$ is $S$-affine, then $X^{\cN}=X^{\bigcap_\alpha N_\alpha}$; in particular $X^{\cN}$ is a closed subscheme of $X$.
\item If $A$ is finite and every $X^{N_\alpha}$ is an algebraic space, then $X^{\cN}$ is an algebraic space.
\item If $X$ is separated, then $X^{\cN}\to X$ is a monomorphism.
\end{enumerate}
\leansrc{xProductRep}
\end{proposition}

\begin{proof}
(i) \cite[proof of Prop.~4.3]{XProducts}. (ii) By Theorem~\ref{thm:x-product}, $X^{\cN}$ is a finite fibre product of algebraic spaces. (iii) By Theorem~\ref{thm:x-product}, $X^{\cN}\to X$ is the fibre product over $X$ of the morphisms $X^{N_\alpha}\to X$, which are monomorphisms by Proposition~\ref{prop:separated-mono}; a fibre product of monomorphisms is a monomorphism.
\end{proof}

\begin{example}[{\cite[Ex.~4.4]{XProducts}}]\label{ex:P1}
Let $M=Z=\ZZ$ and $\cN=(\NN,-\NN)$, so that $\A(\cN)=\mathbb P^1$ with $t\cdot[x:y]=[tx:y]$, and let $X=\mathbb P^1$ over $S=\Spec(\ZZ)$. Over a connected base, $X^{\NN}\to X$ is $(\mathbb P^1\setminus\{\infty\})\sqcup\{\infty\}\to\mathbb P^1$, $X^{-\NN}\to X$ is $(\mathbb P^1\setminus\{0\})\sqcup\{0\}\to\mathbb P^1$, and $X^{\cN}=\mathbb G_m\sqcup\{0\}\sqcup\{\infty\}$, the points being the maps $z\mapsto az$ with $a$ a unit and the two constant maps, while $X^{\NN\cap-\NN}=X^0=\{0,\infty\}$. Hence the formula of Proposition~\ref{prop:x-product-rep}(i) does not extend to non-affine $X$.
\leansrc{exampleP1}
\end{example}

\section{Modules under monoid algebraic spaces and equivariant deformations}\label{sec:modules}

\subsection{Modules over a base}

In this subsection $B$ is a presheaf of sets on $(\Sch/S)_{\fppf}$. Let $(\Sch/B)$ be the category of pairs $(T,b)$ with $T$ an object of $(\Sch/S)_{\fppf}$ and $b\in B(T)$, and let $\cO_B$ be the presheaf of rings $(T,b)\mapsto\cO(T)$ on it; we write $T$ for $(T,b)$. An \emph{$\cO_B$-module} is a presheaf of $\cO_B$-modules on $(\Sch/B)$. An $\cO_B$-module $F$ is \emph{quasi-coherent} if it is a sheaf for the families given by open covers and, for every affine $T$, its restriction to $(\Sch/T)_{\fppf}$ is a quasi-coherent module. Let $F$, $G$ be $\cO_B$-modules and let $X$ be a presheaf of sets on $(\Sch/B)$, with structure morphism $f\colon X\to B$.
\begin{itemize}
\item $\Hom_{\cO_B}(F,G)(X)$ is the set of families $(\varphi_x)$, indexed by the objects $T$ and the $x\in X(T)$, of $\cO(T)$-linear maps $\varphi_x\colon F(T)\to G(T)$ such that $\varphi_{x|_{T'}}(v|_{T'})=\varphi_x(v)|_{T'}$ for every morphism $T'\to T$. A morphism $m\colon X'\to X$ induces $m^*\colon\Hom_{\cO_B}(F,G)(X)\to\Hom_{\cO_B}(F,G)(X')$, $(m^*\varphi)_{x'}=\varphi_{m(x')}$.
\item $f_*f^*G$ is the $\cO_B$-module whose value on $T$ is the set of families $\sigma=(\sigma(t,x))$, indexed by the morphisms $t\colon T'\to T$ and the $x\in X(T')$, with $\sigma(t,x)\in G(T')$ and $\sigma(ts,x|_{T''})=\sigma(t,x)|_{T''}$ for every $s\colon T''\to T'$; that is, $(f_*f^*G)(T)$ is the set of morphisms of presheaves $T\times_BX\to G$. The ring $\cO(T)$ acts through restriction, and a morphism $m\colon X'\to X$ induces $m^\sharp\colon f_*f^*G\to f'_*f'^*G$, $m^\sharp(\sigma)(t,x')=\sigma(t,m(x'))$.
\end{itemize}

\begin{proposition}[{\cite[Prop.~9.9]{AlgMag}}]\label{prop:hom-push}
Let $F$, $G$ be $\cO_B$-modules and $f\colon X\to B$ a presheaf of sets on $(\Sch/B)$.
\begin{enumerate}
\item The map $\varphi\mapsto\psi$, $\psi(v)(t,x)=\varphi_x(v|_{T'})$ for $v\in F(T)$, $t\colon T'\to T$ and $x\in X(T')$, is a bijection
\[\Hom_{\cO_B}(F,G)(X)\xrightarrow{\ \sim\ }\Hom_{\cO_B}(F,f_*f^*G).\]
\item For a morphism $m\colon X'\to X$, this bijection carries $m^*$ to composition with $m^\sharp$.
\item Assume that $f$ is affine, with $\cO_B$-algebra $\cA$, $\cA(T)=\cO(T\times_BX)$, and that $G$ is quasi-coherent. Let $J$ be a Grothendieck topology on $(\Sch/B)$ for which the families of morphisms given by open covers are covering, and assume that $G$ is a sheaf for $J$. Then $f_*f^*G$ is a sheaf for $J$, and the canonical morphism from the presheaf $T\mapsto G(T)\otimes_{\cO(T)}\cA(T)$ to $f_*f^*G$ is a sheafification for $J$.
\end{enumerate}
\leansrc{homPush}
\end{proposition}

\begin{proof}
(i) The inverse sends $\psi$ to $\varphi_x(v)=\psi(v)(\id_T,x)$ for $x\in X(T)$; the two composites are the identity by the compatibility of $\varphi$ and $\psi$ with restrictions. (ii) Both sides send $\varphi$ to $v\mapsto\bigl((t,x')\mapsto\varphi_{m(x')}(v|_{T'})\bigr)$. (iii) An element of $(f_*f^*G)(T)$ is a family of elements of $G(T')$ indexed by the pairs $(t,x)$, $t\colon T'\to T$, $x\in X(T')$, compatible with restrictions. A $J$-covering sieve of $T$ pulls back along every $t\colon T'\to T$ to a $J$-covering sieve of $T'$, so a $J$-covering sieve of $T$ pulls back to a $J$-covering family of $T\times_BX$; since $G$ is a sheaf for $J$, $f_*f^*G$ is a sheaf for $J$. For $T$ affine, $T\times_BX$ is the affine scheme $\Spec(\cA(T))$, $f$ being affine, so $(f_*f^*G)(T)=G(T\times_BX)$, and the canonical map $G(T)\otimes_{\cO(T)}\cA(T)\to G(T\times_BX)$ is bijective by quasi-coherence of $G$ and \cite[Tag~01I8]{Stacks}. The families given by open covers are $J$-covering and every object has an open cover by affine objects, so the canonical morphism is locally injective and locally surjective for $J$; its target being a sheaf for $J$, it is a sheafification for $J$. This is the proof of \cite[Prop.~9.9]{AlgMag}, with these two steps written out.
\end{proof}

For a presheaf of monoids $H$ on $(\Sch/B)$, an \emph{$H$-$\cO_B$-module structure} on an $\cO_B$-module $F$ is a family of actions $H(T)\times F(T)\to F(T)$ of the monoids $H(T)$ by $\cO(T)$-linear maps, compatible with restrictions. With $f\colon H\to B$ the underlying presheaf of sets, a \emph{comodule structure} on $F$ is a morphism of $\cO_B$-modules $\rho\colon F\to f_*f^*F$ such that $\rho(\rho(v)(t,h))(\id_{T'},g)=\rho(v)(t,gh)$ and $\rho(v)(\id_T,1)=v$ for all $v\in F(T)$, $t\colon T'\to T$ and $g,h\in H(T')$. For $f$ affine with $\cO_B$-algebra $\cA$ and $F$ quasi-coherent, these are the conditions $(\rho\otimes\id_\cA)\rho=(\id_F\otimes\Delta_\cA)\rho$ and $(\id_F\otimes\varepsilon_\cA)\rho=\id_F$ on $F\to F\otimes_{\cO_B}\cA$ of \cite[Prop.~9.10]{AlgMag}, by Proposition~\ref{prop:hom-push}(iii).

\begin{proposition}[{\cite[Prop.~9.10]{AlgMag}}]\label{prop:comodule}
Let $H$ be a presheaf of monoids on $(\Sch/B)$ and $F$ an $\cO_B$-module. The map sending an $H$-$\cO_B$-module structure to $\rho(v)(t,g)=g\cdot v|_{T'}$ is a bijection from the set of $H$-$\cO_B$-module structures on $F$ onto the set of comodule structures on $F$.
\leansrc{comodule}
\end{proposition}

\begin{proof}
\cite[proof of Prop.~9.10]{AlgMag}, with Proposition~\ref{prop:hom-push}(i) in place of \cite[Prop.~9.9]{AlgMag}: the inverse sends $\rho$ to $g\cdot v=\rho(v)(\id_T,g)$; associativity and unitality of the action are the two conditions on $\rho$.
\end{proof}

For a monoid $M$, $\A(M)_B$ is the presheaf of monoids $T\mapsto\Hom(\ZZ[M],\cO(T))$ on $(\Sch/B)$.

\begin{proposition}[{\cite[Prop.~9.11]{AlgMag}}]\label{prop:graded-module}
Let $M$ be a countable monoid and $F$ a quasi-coherent $\cO_B$-module. The $\A(M)_B$-$\cO_B$-module structures on $F$, equivalently by Proposition~\ref{prop:comodule} the comodule structures on $F$, correspond bijectively to the families $(\rho_m)_{m\in M}$ of $\cO_B$-linear endomorphisms of $F$ such that
\begin{enumerate}[label=\textup{(\alph*)}]
\item $\rho_m\circ\rho_m=\rho_m$, and $\rho_m\circ\rho_k=0$ for $k\neq m$;
\item for every affine $T$ and every $v\in F(T)$, $\rho_m(v)=0$ for all but finitely many $m$, and $\sum_m\rho_m(v)=v$.
\end{enumerate}
\leansrc{gradedModule}
\end{proposition}

\begin{proof}
\cite[proof of Prop.~9.11]{AlgMag}, with Proposition~\ref{prop:comodule} in place of \cite[Prop.~9.10]{AlgMag}: on an affine $T$, $\rho(v)$ is an element $\sum_m\rho_m(v)X^m$ of $F(T)\otimes_\ZZ\ZZ[M]$, and coassociativity and counitality are (a) and $\sum_m\rho_m(v)=v$.
\end{proof}

\begin{remark}\label{rem:graded-module}
In \cite[Prop.~9.11]{AlgMag}, condition (b) is replaced by the condition that only finitely many of the $\rho_m$ are nonzero locally. For a nonzero ring $R$, the comodule $R[M]$ over $R[M]$, with $\rho(X^m)=X^m\otimes X^m$, has all its components $\rho_m$ nonzero; for $M$ infinite it does not satisfy this condition, and it satisfies the analogue of (b).
\leansrc{regularComp}
\end{remark}

Let $H$ be a presheaf of monoids on $(\Sch/B)$ and $F$ an $\cO_B$-module with an $H$-$\cO_B$-module structure. For an object $V$ of $(\Sch/B)$ and $n\geq0$, the $\cO(V)$-module $C^n(H,F)(V)$ of \cite[Def.~9.4, Def.~9.5]{AlgMag} is the set of families $c=(c_t)$, indexed by the morphisms $t\colon T\to V$, of maps $c_t\colon H(T)^n\to F(T)$ such that $c_{ts}(g_1|_{T'},\dots,g_n|_{T'})=c_t(g_1,\dots,g_n)|_{T'}$ for every $s\colon T'\to T$. The differential $\delta\colon C^n(H,F)(V)\to C^{n+1}(H,F)(V)$ is
\begin{multline*}(\delta c)_t(g_0,\dots,g_n)=g_0\cdot c_t(g_1,\dots,g_n)\\+\sum_{i=1}^n(-1)^ic_t(g_0,\dots,g_{i-1}g_i,\dots,g_n)+(-1)^{n+1}c_t(g_0,\dots,g_{n-1}),\end{multline*}
and $H^n_{\mathrm{monoid}}(H|_V,F|_V)$ is the quotient of the kernel of $\delta$ on $C^n(H,F)(V)$ by the image of $\delta$ from $C^{n-1}(H,F)(V)$ ($0$ for $n=0$) \cite[Def.~9.7]{AlgMag}. For $B$ itself, $C^n(H,F)(B)$ is the set of families $(c_T)$ indexed by the objects $T$, with the same compatibility, and $\delta$ is given by the same formula.

\begin{proposition}[{\cite[Prop.~9.12]{AlgMag}}]\label{prop:cochain-formula}
Let $M$ be a countable commutative monoid, $F$ a quasi-coherent $\cO_B$-module with an $\A(M)_B$-$\cO_B$-module structure, $(\rho_m)_{m\in M}$ the family of Proposition~\ref{prop:graded-module}, and $n\geq0$. For an object $V$ of $(\Sch/B)$, let
\begin{gather*}\phi_V\colon\bigoplus_{\mathbf m\in M^n}F(V)\longrightarrow C^n(\A(M)_B,F)(V),\\ fX^{m_1}\cdots X^{m_n}\longmapsto\bigl((g_1,\dots,g_n)\mapsto g_1(m_1)\cdots g_n(m_n)\,f|_T\bigr)_{t\colon T\to V}.\end{gather*}
\begin{enumerate}
\item For every object $V$, the map sending $c$ to the family $(\phi_E^{-1}(c|_E))_{a\colon E\to V}$, indexed by the morphisms $a\colon E\to V$ with $E$ affine, is a bijection from $C^n(\A(M)_B,F)(V)$ onto the set of families $(x_a)$ with $x_a\in\bigoplus_{\mathbf m\in M^n}F(E)$ and $x_{as}=s^*x_a$ for every $s\colon E'\to E$ with $E'$ affine.
\item If $V$ is quasi-compact, $\phi_V$ is an isomorphism of $\cO(V)$-modules.
\item If $V$ is quasi-compact, $\delta\circ\phi_V=\phi_V\circ\delta'$, where
\begin{multline*}\delta'\bigl(fX^{m_1}\cdots X^{m_n}\bigr)=\sum_{k\in M}\rho_k(f)X^kX^{m_1}\cdots X^{m_n}\\+\sum_{i=1}^n(-1)^if\,X^{m_1}\cdots X^{m_i}X^{m_i}\cdots X^{m_n}+(-1)^{n+1}f\,X^{m_1}\cdots X^{m_n}X^0.\end{multline*}
\end{enumerate}
If $B$ is a scheme, (i)--(iii) hold for every commutative monoid $M$.
\leansrc{cochainFormula}
\end{proposition}

\begin{proof}
A cochain is a morphism of presheaves $\A(M)_V^n=\A(M^n)_V\to F$, an element of $\Hom_{\cO_B}(\cO_B,F)(\A(M^n)_V)$. On an affine $E\to V$, $\A(M^n)_E=\Spec(\cO(E)[M^n])$ is an object of $(\Sch/B)$, $M$ being countable, and evaluation at its universal point gives $F(E)\otimes_{\cO(E)}\cO(E)[M^n]=\bigoplus_{\mathbf m\in M^n}F(E)$, by quasi-coherence of $F$. This is (i), and (ii) follows since a section over a quasi-compact $V$ has finitely many nonzero coefficients, by a finite affine open cover. (iii) is the formula for $\delta$ on the monomial cochain $fX^{\mathbf m}$, with $g_0\cdot f=\sum_kg_0(k)\rho_k(f)$. This is \cite[proof of Prop.~9.12]{AlgMag}.
\end{proof}

\begin{proposition}[{\cite[Prop.~9.13, Rem.~9.14]{AlgMag}}]\label{prop:monoid-cohomology}
Let $M$ be a countable commutative monoid and $F$ a quasi-coherent $\cO_B$-module with an $\A(M)_B$-$\cO_B$-module structure. Let $n\geq0$. For every object $V$ of $(\Sch/B)$, and for $V=B$, every $c\in C^{n+1}(\A(M)_B,F)(V)$ with $\delta c=0$ is of the form $c=\delta e$ with $e\in C^n(\A(M)_B,F)(V)$. In particular $H^n_{\mathrm{monoid}}(\A(M)_V,F|_V)=0$ and $H^n_{\mathrm{monoid}}(\A(M)_B,F)=0$ for all $n\geq1$. If $B$ is a scheme, this holds for every commutative monoid $M$.
\leansrc{monoidCohomology}
\end{proposition}

\begin{proof}
For a cochain $x=\sum_{\mathbf w\in M^{n+1}}x_{\mathbf w}X^{\mathbf w}$ on an affine $E\to V$, in the form of Proposition~\ref{prop:cochain-formula}, let $\epsilon(x)=\sum_{\mathbf m\in M^n}x_{(\mathbf m,0)}X^{\mathbf m}$ be the coefficient of $X^0$ in the last variable. The formula of Proposition~\ref{prop:cochain-formula}(iii) gives $\epsilon(\delta x)=\delta(\epsilon(x))+(-1)^nx$, with no condition on the $\rho_k$. If $\delta x=0$, then $x=\delta\bigl((-1)^{n+1}\epsilon(x)\bigr)$. For $n=0$ this is the element $e_U=-f_0$ of \cite[proof of Prop.~9.13]{AlgMag}. The elements $(-1)^{n+1}\epsilon(x)$ are given by a formula compatible with restriction along the morphisms of affine objects, so by Proposition~\ref{prop:cochain-formula}(i) they form a cochain $e$ over $V$, and over $B$, with $\delta e=c$, as in \cite[proof of Prop.~9.13]{AlgMag}.
\end{proof}

\begin{remark}\label{rem:monoid-cohomology}
This proves the statement conjectured in \cite[Rem.~9.14]{AlgMag}.
\end{remark}

\subsection{Equivariant infinitesimal deformations}

In this subsection $S=\Spec(R)$, $H$ is a commutative Hopf $R$-algebra and $G=\Spec(H)$. For an $R$-algebra $E$, $G(E)$ is the group of $R$-algebra homomorphisms $H\to E$, and $u=\id_H\in G(H)$. An action of $G$ on an affine $S$-scheme $\Spec(A)$ is given by a coassociative and counital $R$-algebra homomorphism $\rho\colon A\to H\otimes_RA$. An $S$-action of $G$ on an $S$-scheme $Y$ is given by actions $G(E)\times Y(E)\to Y(E)$ for the $R$-algebras $E$, natural in $E$.

\begin{proposition}[{\cite[Prop.~10.1]{AlgMag}}]\label{prop:equivariant-lift}
Let $X'=\Spec(A_1)$ and $X=\Spec(A)$ be affine $S$-schemes with actions $\rho_1$, $\rho_0$ of $G$, and let $\phi\colon A_1\to A$ be a surjective homomorphism of $R$-algebras whose kernel has square zero and with $(\id_H\otimes\phi)\rho_1=\rho_0\phi$; put $i=\Spec(\phi)\colon X\to X'$. Let $Y$ and $Z$ be $S$-schemes with $S$-actions of $G$, and let $f\colon Y\to Z$, $a\colon X\to Y$ and $b\colon X'\to Z$ be equivariant $S$-morphisms with $f\circ a=b\circ i$. For an $R$-algebra $E$, let $C_E$ be the kernel of $E\otimes_RA_1\to E\otimes_RA$ and $\mathcal F(E)=\Hom(a_E^*\Omega_{Y/Z},C_E)$, with the action of $G(E)$ and the maps $\mathcal F(E)\to\mathcal F(E')$ of \cite[proof of Prop.~10.1]{AlgMag}. Assume:
\begin{enumerate}[label=\textup{(\alph*)}]
\item every $c\in\mathcal F(H)$ with $\Delta^*c=i_1^*c+i_1\cdot i_2^*c$ in $\mathcal F(H\otimes_RH)$, where $\Delta$ is the comultiplication and $i_1,i_2\colon H\to H\otimes_RH$ the two inclusions, is of the form $c=u\cdot v_H-v_H$ for some $v\in\mathcal F(R)$;
\item there is an $S$-morphism $\alpha\colon X'\to Y$ with $\alpha\circ i=a$ and $f\circ\alpha=b$.
\end{enumerate}
Then there is an $S$-morphism $l\colon X'\to Y$ with $l\circ i=a$ and $f\circ l=b$ such that $l\circ g_{X'_E}=g\cdot l_E$ as morphisms $X'_E\to Y$ for every $R$-algebra $E$ and every $g\in G(E)$, where $l_E\colon X'_E\to X'\to Y$ and $g_{X'_E}$ is the action of $g$ on $X'_E$.
\leansrc{equivariantLift}
\end{proposition}

If $H$ is flat over $R$, condition (a) is the vanishing of $H^1_{\mathrm{monoid}}(G,\mathcal F)$ of \cite[Prop.~10.1]{AlgMag}, written with the universal point $u$. If $Y\to Z$ is formally smooth, condition (b) holds by the definition of formal smoothness, $X'$ being affine.

\begin{proof}
\cite[proof of Prop.~10.1]{AlgMag}: the $1$-cocycle $g\mapsto g\cdot\alpha-\alpha$ of $G$ with values in $\mathcal F$ is determined by its value $c=u\cdot\alpha_H-\alpha_H\in\mathcal F(H)$, which satisfies $\Delta^*c=i_1^*c+i_1\cdot i_2^*c$. By (a), $c=u\cdot v_H-v_H$ with $v\in\mathcal F(R)$, and $l=\alpha-v$ is invariant under $G(E)$ for every $R$-algebra $E$.
\end{proof}

For $R=\ZZ$ and $G=\D(M)$ with $M$ a commutative group, the conclusion of Proposition~\ref{prop:equivariant-lift} holds for $S$-schemes $Y$, $Z$ without condition (a): the vanishing of $H^1_{\mathrm{monoid}}$ needed in the proof is proved for the coefficients at hand, as in Proposition~\ref{prop:monoid-cohomology}, and $l=\alpha-v$ as in the proof.
\leansrc{invariantLift}

\begin{corollary}[{\cite[Cor.~10.3]{AlgMag}}]\label{cor:equivariant-lift-spaces}
Let $S$ be a scheme and $M$ a commutative group of cardinality less than the cutoff of the site $(\Sch/S)_{\fppf}$. Let $A_1$, $A$ be rings, $\phi\colon A_1\to A$ a surjective ring homomorphism whose kernel has square zero, $X'=\Spec(A_1)$ with a morphism $X'\to S$, $X=\Spec(A)$ over $S$ through $i=\Spec(\phi)\colon X\to X'$, both objects of the site. Let $Y$, $Z$ be algebraic spaces over $S$. Let $X'$, $X$, $Y$, $Z$ carry actions of $\D(M)_S$, and let $i$, $a\colon X\to Y$, $b\colon X'\to Z$ and $f\colon Y\to Z$ be equivariant, with $f$ formally smooth and $f\circ a=b\circ i$. Then there is an equivariant morphism $l\colon X'\to Y$ with $l\circ i=a$ and $f\circ l=b$.
\leansrc{equivariantLiftSpaces}
\end{corollary}

\begin{proof}
The action of $\D(M)_S$ on $X'$ at the universal point of $\D(M)$ over $\Spec(A_1[M])$ is a ring homomorphism $\rho\colon A_1\to A_1[M]$, coassociative and counital, and similarly $\rho_0\colon A\to A[M]$ with $\phi[M]\circ\rho=\rho_0\circ\phi$. Formal smoothness of $f$, $X'$ being affine, gives $\alpha\colon X'\to Y$ with $\alpha\circ i=a$ and $f\circ\alpha=b$. The integral of $g\mapsto g^{-1}\cdot\alpha(g\,\cdot)$ over $\D(M)$ is an equivariant morphism $l$ with the same two properties; this is \cite[proof of Cor.~10.3]{AlgMag}, Proposition~\ref{prop:monoid-cohomology} for the coefficients at hand replacing \cite[Prop.~9.13]{AlgMag}. The equivariance of $l$ on every object of the site follows from the equivariance at the universal point.
\end{proof}

\begin{remark}\label{rem:equivariant-lift-spaces}
In \cite[Cor.~10.3]{AlgMag}, $M$ is finitely generated, $X'\to S$ is quasi-compact and quasi-separated, and $Y\to Z$ is locally of finite type; these assumptions are not used in Corollary~\ref{cor:equivariant-lift-spaces}, where $X'$ is affine.
\end{remark}

\section{Formal smoothness}\label{sec:smooth}

In this section $M$ is an abelian group. Following \cite[Tag~049S]{Stacks}, a morphism $g\colon F\to G$ of presheaves on $(\Sch/S)_{\fppf}$ is \emph{formally smooth as a transformation of functors} if for every first order thickening $T_0\to T$ of affine $S$-schemes, every $x_0\in F(T_0)$ and every $y\in G(T)$ with $g(x_0)=y|_{T_0}$, there is $x\in F(T)$ with $x|_{T_0}=x_0$ and $g(x)=y$.

\begin{proposition}[{\cite[Prop.~11.1]{AlgMag}}]\label{prop:attractor-fsmooth}
Let $X$ and $Y$ be algebraic spaces over $S$ with actions of $\D(M)_S$, let $f\colon X\to Y$ be an equivariant smooth morphism and let $N$ be a magnet. Then $f^N\colon X^N\to Y^N$ is formally smooth as a transformation of functors.
\leansrc{smoothAttractor}
\end{proposition}

\begin{proof}
The morphism $f$ is formally smooth \cite[Tag~04AM]{Stacks}. Let $T_0\to T$ be a first order thickening of affine $S$-schemes; then $\A(N)_{T_0}\to\A(N)_T$ is a first order thickening of affine schemes. Given a commutative square formed by an equivariant $\varphi_0\colon\A(N)_{T_0}\to X$ and an equivariant $\phi\colon\A(N)_T\to Y$, formal smoothness of $f$ gives a morphism $\alpha\colon\A(N)_T\to X$ with $\alpha|_{\A(N)_{T_0}}=\varphi_0$ and $f\circ\alpha=\phi$. The integral over $\D(M)$ of $g\mapsto g^{-1}\cdot\alpha(g\,\cdot)$ is an equivariant morphism $\varphi\colon\A(N)_T\to X$ with the same two properties: this is Corollary~\ref{cor:equivariant-lift-spaces} for the affine scheme $\A(N)_T$. Then $\varphi$ is the required point of $X^N(T)$, as in \cite[proof of Prop.~11.1]{AlgMag}.
\end{proof}

\begin{corollary}[{\cite[Cor.~11.2]{AlgMag}}]\label{cor:attractor-smooth}
In the situation of Proposition~\ref{prop:attractor-fsmooth}, assume that $X^N$ and $Y^N$ are algebraic spaces. Then $f^N\colon X^N\to Y^N$ is smooth.
\leansrc{smoothAttractorRep}
\end{corollary}

\begin{proof}
The morphism $f$ is locally of finite presentation, so $f^N$ is locally of finite presentation (Proposition~\ref{prop:attractor-lfp}), and it is formally smooth by Proposition~\ref{prop:attractor-fsmooth}; hence $f^N$ is smooth \cite[Tag~04AM]{Stacks}, as in \cite[proof of Cor.~11.2]{AlgMag}.
\end{proof}

\begin{proposition}[{\cite[Prop.~11.7]{AlgMag}}]\label{prop:face-fsmooth}
Let $X$ be an algebraic space over $S$ with an action of $\D(M)_S$ such that $X\to S$ is formally smooth, let $N$ be a magnet and $F$ a face of $N$. Then $p_{N,F}\colon X^N\to X^F$ is formally smooth as a transformation of functors.
\leansrc{smoothFace}
\end{proposition}

\begin{proof}
\cite[proof of Prop.~11.7]{AlgMag}. Let $T_0=\Spec(D_1)\to T=\Spec(D)$ be a first order thickening of affine $S$-schemes, $\varphi_0\in X^N(T_0)$ and $\phi\in X^F(T)$ with $p_{N,F}(\varphi_0)=\phi|_{T_0}$. With $R$, $\nu_1$, $\nu_2$ as in Proposition~\ref{prop:fact118}, the equivariant morphisms $\phi\colon\A(F)_T\to X$ and $\varphi_0\colon\A(N)_{T_0}\to X$ glue to a morphism $k\colon\Spec(R)\to X$ (Proposition~\ref{prop:fact118}(iii)), and $\Spec(R)\to\A(N)_T$ is a first order thickening of affine schemes (Proposition~\ref{prop:fact118}(ii)). Since $X\to S$ is formally smooth, $k$ extends to a morphism $\alpha\colon\A(N)_T\to X$ over $S$. The integral over $\D(M)$ of $g\mapsto g^{-1}\cdot\alpha(g\,\cdot)$, as in the proof of Proposition~\ref{prop:attractor-fsmooth}, is an equivariant extension $\varphi\colon\A(N)_T\to X$ of $k$; then $\varphi\in X^N(T)$, $\varphi|_{T_0}=\varphi_0$ and $p_{N,F}(\varphi)=\phi$.
\end{proof}

\begin{corollary}[{\cite[Cor.~11.10]{AlgMag}}]\label{cor:face-smooth}
Let $X$ be an algebraic space over $S$ with an action of $\D(M)_S$, smooth over $S$, let $N$ be a magnet and $F$ a face of $N$, and assume that $X^N$ and $X^F$ are algebraic spaces. Then $p_{N,F}\colon X^N\to X^F$ is smooth.
\leansrc{smoothFaceRep}
\end{corollary}

\begin{proof}
Let $S$ carry the trivial action, so that $S^N=S^F=S$. By Proposition~\ref{prop:attractor-lfp} applied to $X\to S$, the morphisms $X^N\to S$ and $X^F\to S$ are locally of finite presentation, so $p_{N,F}$ is locally of finite presentation \cite[Tag~05WT]{Stacks}. It is formally smooth by Proposition~\ref{prop:face-fsmooth}, hence smooth \cite[Tag~04AM]{Stacks}, as in \cite[proof of Cor.~11.10]{AlgMag}.
\end{proof}

\section{Topology and geometry of attractors}\label{sec:topology}

In this section $M$ is a finitely generated abelian group and $X$ is an algebraic space over $S$ with an action of $\D(M)_S$. For an algebraic space $Y$, $|Y|$ denotes its underlying topological space \cite[Tag~03BY]{Stacks} and $\pi_0(|Y|)$ its set of connected components.

\begin{proposition}[{\cite[Prop.~12.2]{AlgMag}}]\label{prop:immersions}
Let $f\colon X\to Y$ be an equivariant morphism of algebraic spaces over $S$ with $\D(M)_S$-actions and let $N$ be a magnet.
\begin{enumerate}
\item Let $F$ be a face of $N$ that is a subgroup of $M$, for instance the group of units $\units N$, and let $e_F$ be the vertex of $\A(N)$ at $F$, the point with value $1$ on $F$ and $0$ on $N\smallsetminus F$. If $f$ is an open immersion, then $X^N=Y^N\times_YX$, the morphism $X^N\to Y^N\times_YX$ being given by $f^N$ and evaluation at $e_F$. In particular, if $X^N$ and $Y^N$ are representable, $f^N$ is an open immersion.
\item If $f$ is a closed immersion, then $X^N=Y^N\times_YX$, the morphism $X^N\to Y^N\times_YX$ being given by $f^N$ and evaluation at the unit section. In particular, if $X^N$ and $Y^N$ are representable, $f^N$ is a closed immersion.
\end{enumerate}
\leansrc{immersions}
\end{proposition}

\begin{proof}
(i) For a scheme $T$ over $S$, an equivariant $\psi\colon\A(N)_T\to Y_T$ and $x\in X(T)$ with $f(x)=\psi(e_F)$, there is a unique equivariant $\varphi\colon\A(N)_T\to X_T$ with $f\circ\varphi=\psi$ and $\varphi(e_F)=x$; this is the inverse of $X^N(T)\to Y^N(T)\times_{Y(T)}X(T)$. The last assertion follows with $F=\units N$, open immersions being stable under base change; it is \cite[Prop.~12.2(i)]{AlgMag}, proved there by smoothness and the monomorphism property of $f^N$. (ii) Let $T$ be a scheme over $S$, let $\psi\colon\A(N)_T\to Y_T$ be equivariant and let $x\in X(T)$ with $f(x)=\psi(1)$. Consider the orbit morphism $o\colon\A(M)_T\to\A(N)_T$, $h\mapsto h\cdot1$; locally on $T=\Spec(R)$ it corresponds to the inclusion $R[N]\to R[M]$, which is injective, so the scheme-theoretic image of $o$ is $\A(N)_T$ \cite[Tag~01R8]{Stacks}. For every point $h$ of $\A(M)_T$, $\psi(o(h))=h\cdot\psi(1)=f(h\cdot x)$, so $\psi\circ o$ factors through $X_T$. Hence $o$ factors through the closed subscheme $\psi^{-1}(X_T)$ of $\A(N)_T$, which is therefore $\A(N)_T$. So $\psi=f\circ\varphi$ for a unique $\varphi\colon\A(N)_T\to X_T$, since $f$ is a monomorphism; $\varphi$ is equivariant and $\varphi(1)=x$, again because $f$ is a monomorphism. This gives the inverse of $X^N(T)\to Y^N(T)\times_{Y(T)}X(T)$. The last assertion follows, closed immersions being stable under base change.
\end{proof}

\begin{proposition}[{\cite[Prop.~12.3]{AlgMag}}]\label{prop:deformation}
Let $N$ be a finitely generated magnet and $F$ a face of $N$, and let $X$ be any sheaf on $(\Sch/S)_{\fppf}$ with an action of $\D(M)_S$. Then $(\iota_{F,N},p_{N,F})$ is a strong deformation retract: there is a morphism $\Psi\colon X^N\times_S\mathbb{A}^1_S\to X^N$ such that
\begin{enumerate}
\item $\Psi\circ j_0=\iota_{F,N}\circ p_{N,F}$, where $j_0\colon X^N\to X^N\times_S\mathbb{A}^1_S$ is given by the section $0$;
\item $\Psi\circ j_1=\id_{X^N}$, where $j_1$ is given by the section $1$;
\item $\Psi\circ(\iota_{F,N}\times\id_{\mathbb{A}^1_S})=\iota_{F,N}\circ\mathrm{pr}_1$;
\item $p_{N,F}\circ\Psi=p_{N,F}\circ\mathrm{pr}_1$.
\end{enumerate}
\leansrc{deformation}
\end{proposition}

\begin{proof}
Let $w\colon N\to\NN$ be a morphism of monoids with $w^{-1}(0)=F$ (Lemma~\ref{lem:face-weight}). For a scheme $T$ and $c\in\Gamma(T,\cO_T)$, let $\lambda_c\colon\A(N)_T\to\A(N)_T$ be the morphism given by $X^n\mapsto c^{w(n)}X^n$, with $0^0=1$; it is $\D(M)_T$-equivariant, since it preserves degrees. Put $\Psi(\varphi,c)=\varphi\circ\lambda_c$. Then $\lambda_1=\id$, which is (ii). The morphism $\lambda_0$ is given by $X^n\mapsto X^n$ for $n\in F$ and $X^n\mapsto0$ otherwise, which is the composite $\A(N)_T\to\A(F)_T\to\A(N)_T$ of the morphisms of Proposition~\ref{prop:action-AN}(ii), (iii); this gives (i). Finally the composite of $\lambda_c$ with $\A(N)_T\to\A(F)_T$ is $\A(N)_T\to\A(F)_T$, since $w$ vanishes on $F$, which gives (iii) and, composing on the other side, $p_{N,F}(\varphi\circ\lambda_c)=p_{N,F}(\varphi)$, which is (iv). This is \cite[proof of Prop.~12.3]{AlgMag}, with $w$ in place of \cite[Prop.~3.3.1(iii), p.~75]{Ogus}.
\end{proof}

\begin{proposition}[{\cite[Prop.~12.1]{AlgMag}}]\label{prop:pi0}
Let $N$ be a finitely generated magnet and $F$ a face of $N$. Assume that $X$ is a separated algebraic space and has a $\units N$-FPR atlas. Then $X^N$ and $X^F$ are algebraic spaces, the morphism $p_{N,F}\colon X^N\to X^F$ has connected fibres over all field-valued points of $X^F$, and it induces a bijection $\pi_0(|X^N|)\to\pi_0(|X^F|)$.
\leansrc{pi0}
\end{proposition}

\begin{remark}\label{rem:pi0-representable}
The proof below uses only Proposition~\ref{prop:deformation} and the representability of $X^N$ and $X^F$: for any algebraic space $X$ over $S$ with an action of $\D(M)_S$ such that $X^N$ and $X^F$ are representable, $p_{N,F}$ has connected fibres over all field-valued points of $X^F$ and induces a bijection $\pi_0(|X^N|)\to\pi_0(|X^F|)$.
\end{remark}

\begin{proof}
The face $F$ is finitely generated and $\units F=\units N$, so $X^N$ and $X^F$ are algebraic spaces by Theorem~\ref{thm:representability}(a). Let $x$ be a point of $X^F$ with values in a field $K$, $L\supseteq K$ a field and $y$ an $L$-point of the fibre $p_{N,F}^{-1}(x)$. With $\Psi$ as in Proposition~\ref{prop:deformation}, the morphism $\mathbb{A}^1_L\to X^N$, $t\mapsto\Psi(y,t)$, lands in the fibre, by Proposition~\ref{prop:deformation}(iv), and joins $y$ ($t=1$) to $\iota_{F,N}(x_L)$ ($t=0$). So every point of the fibre lies in the connected component of $\iota_{F,N}(x)$, and the fibre is connected. The map $|p_{N,F}|$ is continuous and has the continuous section $|\iota_{F,N}|$, and every point $y$ of $|X^N|$ lies in the connected component of $\iota_{F,N}(p_{N,F}(y))$; hence the bijection on $\pi_0$, as in \cite[proof of Prop.~12.1]{AlgMag}.
\end{proof}

\begin{theorem}[{\cite[Thm.~12.4]{AlgMag}}]\label{thm:BB}
Assume that $X$ is smooth over $S$ and that the action of $\D(M)_S$ on $X$ is Zariski locally linearizable, with a cover $X=\bigcup_aU_a$ by $\D(M)_S$-stable open subspaces that are $S$-affine schemes. Let $N$ be a magnet and $F$ a face of $N$.
\begin{enumerate}
\item The morphisms $U_a^N\to X^N$ and $U_a^F\to X^F$ are open immersions, $U_a^N=X^N\times_{X^F}U_a^F$, the $U_a^N$ cover $X^N$ and the $U_a^F$ cover $X^F$; the $U_a^N$ and $U_a^F$ are affine over $S$ and smooth over $S$.
\item Let $V\subseteq S$ be an affine open subscheme and $U_a\times_SV=\Spec(A)$. Put $A'=\cO(U_a^N\times_SV)$ and $A'_F=\cO(U_a^F\times_SV)$, let $\pi\colon A'\to A'_F$ be the restriction of functions and $I=\ker(\pi)$. Then $p_{N,F}$ corresponds to a section $\rho\colon A'_F\to A'$ of $\pi$, the $A'_F$-module $I/I^2$ is finite projective, and there is an isomorphism of $A'_F$-algebras $\Sym_{A'_F}(I/I^2)\cong A'$, under which $p_{N,F}$ is the projection $\Spec(\Sym_{A'_F}(I/I^2))\to\Spec(A'_F)$ and $\iota_{F,N}$ is the zero section.
\end{enumerate}
\leansrc{BB}
\end{theorem}

For $F=\units N$, Theorem~\ref{thm:BB} is the affine space fibration $X^N\to X^{\units N}$ of \cite[Thm.~12.4]{AlgMag}.

\begin{proof}
Let $q\colon M\to M/F^{\gp}$ be the quotient map and $\bar N=q(N)$. As in the proof of Lemma~\ref{lem:face-weight}, $N\cap F^{\gp}=F$ and $\bar N$ is sharp.

(i) Let $U\subseteq X$ be a $\D(M)_S$-stable open subspace and $L$ a face of $N$. Both $U^N$ and $X^N\times_{X^L}U^L$ are subfunctors of $X^N$, and they are equal: a point $\psi\colon\A(N)_T\to X_T$ of $X^N$ whose restriction to $\A(L)_T$ lands in $U$ lands in $U$. The open subscheme $W=\psi^{-1}(U_T)$ is $\D(M)_T$-stable and contains $\A(L)_T$. Let $z$ be a point of $\A(N)_T$ with values in a field $K$, and let $e_L\in\A(N)(K)$ be the point with $e_L(n)=1$ for $n\in L$ and $e_L(n)=0$ otherwise; then $e_L\cdot z\in\A(L)_K\subseteq W$. The set of $v\in\A(N)_K$ with $v\cdot z\in W$ is open, $\D(M)_K$-stable and contains $e_L$. If its closed complement contained $1$, it would contain the closure of the orbit $\D(M)_K\cdot1$, which is dense in $\A(N)_K$ since $K[N]\subseteq K[M]$; so it would contain $e_L$. Hence $z=1\cdot z\in W$, and $W=\A(N)_T$. This gives $U_a^N=X^N\times_{X^F}U_a^F$, and, applied to the monoid $F$ and its face $\units F=\units N$, $U_a^F=X^F\times_{X^{\units N}}U_a^{\units N}$. The open immersion $U_a\to X$ is strongly-FPR (Proposition~\ref{prop:fpr-basic}(iv)), so $U_a^{\units N}=X^{\units N}\times_XU_a$; these are open subfunctors of $X^{\units N}$ covering it, and by the two cartesian squares the $U_a^F$ and the $U_a^N$ are open subfunctors covering $X^F$ and $X^N$. They are $S$-affine schemes by Theorem~\ref{thm:affine}, smooth over $S$ by \cite[Cor.~11.2]{AlgMag}.

(ii) Write $A'=\bigoplus_{n\in N}A'_n$. By Theorem~\ref{thm:affine} and Proposition~\ref{prop:intersection}(i), $U_a^F\times_SV=(U_a^N\times_SV)^F$ is defined in $\Spec(A')$ by the $A'_n$, $n\in N\setminus F$. Since $N\setminus F$ is an ideal of $N$, $I=\bigoplus_{n\in N\setminus F}A'_n$, and $\rho$ is the inclusion of $A'_F=\bigoplus_{n\in F}A'_n$. The algebra $A$ is generated by finitely many homogeneous elements; their degrees lying in $N$ generate a finitely generated submonoid $N_0\subseteq N$ containing the support of $A'$, and $F\cap N_0$ is a face of $N_0$; so we may assume $N$ finitely generated. Grade $A'$ by $\bar N$, $A'_\nu=\bigoplus_{q(n)=\nu}A'_n$; then $A'_0=A'_F$ and $I=\bigoplus_{\nu\neq0}A'_\nu$. Let $h\colon\bar N\to\NN$ be a morphism with $h^{-1}(0)=0$ (proof of Lemma~\ref{lem:face-weight}). The closed immersion $U_a^F\times_SV\to U_a^N\times_SV$ of smooth $V$-schemes \cite[Cor.~11.2]{AlgMag} is Koszul-regular \cite[Tag~067U]{Stacks}, hence quasi-regular \cite[Tag~063K]{Stacks}, and $I/I^2$ is finite projective over $A'_F$ \cite[Tag~063M]{Stacks}. From here \cite[proof of Thm.~12.4]{AlgMag} applies to the $\bar N$-graded algebra $A'$, with $A'_F$ in degree $0$ and the weight $h$: the graded pieces $(I/I^2)_\nu$ are direct summands of $I/I^2$, hence projective over $A'_F$; graded sections $(I/I^2)_\nu\to I_\nu$ define $\phi\colon\Sym_{A'_F}(I/I^2)\to A'$, surjective by induction on $h$ and injective since $\operatorname{Gr}_I(A')\cong\Sym_{A'_F}(I/I^2)$ \cite[Tags~062Z, 063M]{Stacks}. The isomorphism $\phi$ is the identity on $A'_F$ and maps $I/I^2$ into $I$, which gives the projection and the zero section.
\end{proof}

\begin{corollary}\label{cor:affine-fibration}
In the situation of Theorem~\ref{thm:BB}, the morphism $p_{N,F}\colon X^N\to X^F$ is an affine space fibration: Zariski-locally on $X^F$, $X^N\cong X^F\times_S\mathbb{A}^r_S$. More precisely, in the notation of Theorem~\ref{thm:BB}(ii), there are $f_1,\dots,f_k\in A'_F$ generating the unit ideal and integers $r_j\geq0$ such that $(I/I^2)_{f_j}$ is free of rank $r_j$ over $(A'_F)_{f_j}$ and $A'_{f_j}\cong(A'_F)_{f_j}[x_1,\dots,x_{r_j}]$ as $(A'_F)_{f_j}$-algebras.
\leansrc{affineFibration}
\end{corollary}

\begin{proof}
The finite projective $A'_F$-module $I/I^2$ is finitely presented, so it is free of finite rank on the members of a finite cover of $\Spec(A'_F)$ by basic opens $D(f_j)$ \cite[Tag~00NX]{Stacks}. The symmetric algebra commutes with base change, $(A'_F)_{f_j}\otimes_{A'_F}\Sym_{A'_F}(I/I^2)\cong\Sym_{(A'_F)_{f_j}}((I/I^2)_{f_j})$, and the symmetric algebra of a free module of rank $r_j$ is a polynomial algebra in $r_j$ variables. With Theorem~\ref{thm:BB}(ii) this gives $A'_{f_j}\cong(A'_F)_{f_j}\otimes_{A'_F}A'\cong(A'_F)_{f_j}[x_1,\dots,x_{r_j}]$, compatibly with $\rho$.
\end{proof}

An element $e$ of a sharp monoid $N$ is \emph{irreducible} if $e\neq0$ and $e=x+y$ with $x,y\in N$ implies $x=0$ or $y=0$; $\mathrm{Irr}(N)$ denotes the set of irreducible elements of $N$.

\begin{corollary}[{\cite[Cor.~12.5]{AlgMag}}]\label{cor:toric-smooth}
Let $N$ be a fine sharp monoid. If $\A(N)_S\to S$ is smooth, then $\A(N)_S\to S$ is a vector bundle. If moreover $S$ is nonempty, $N$ is the free monoid on its irreducible elements and $\A(N)_S\cong\mathbb{A}^r_S$, where $r$ is the number of irreducible elements of $N$.
\leansrc{toricSmooth}
\end{corollary}

\begin{proof}
The first assertion is \cite[proof of Cor.~12.5]{AlgMag}. For the second, let $s\in S$ and $k=\kappa(s)$. Since $N$ is sharp, the map $X^n\mapsto0$ for $n\neq0$ is a $k$-point $v$ of $\A(N)_k$ with maximal ideal $\mathfrak m$ spanned by the $X^n$, $n\neq0$; $\mathfrak m^2$ is spanned by the $X^n$ with $n$ a sum of two nonzero elements, so the $X^e$, $e$ irreducible, form a basis of $\mathfrak m/\mathfrak m^2$. The irreducible elements generate $N$, by induction on $w(n)$ for $w$ as in Lemma~\ref{lem:face-weight} with $F=0$. Since $\A(N)_k\to\Spec(k)$ is smooth, the proof of Theorem~\ref{thm:BB}(ii) \cite[proof of Thm.~12.4]{AlgMag}, applied to $X=\A(N)_k$ with the action of $\D(N^{\gp})_k$ by translation and $F=0$, with the section of $\mathfrak m\to\mathfrak m/\mathfrak m^2$ sending the class of $X^e$ to $X^e$ for $e$ irreducible, gives an isomorphism of $k$-algebras $\Sym_k(\mathfrak m/\mathfrak m^2)\cong k[N]$ sending the class of $X^e$ to $X^e$. Hence the monomials in the $X^e$, $e\in\mathrm{Irr}(N)$, are linearly independent, so the morphism $\NN^{(\mathrm{Irr}(N))}\to N$ is injective, and $N=\mon{\mathrm{Irr}(N)}\cong\NN^r$. The corresponding morphism $\A(N)_S\to\A(\NN^{(\mathrm{Irr}(N))})_S=\mathbb{A}^r_S$ is an isomorphism over $S$.
\end{proof}

\section{Dilatations}\label{sec:dilatations}

Let $S'\subseteq S$ be a closed subscheme whose ideal is locally generated by one element. A scheme $T$ over $S$ is \emph{$S'$-regular} if every local generator of the ideal of $S'$ pulls back to a nonzerodivisor of $\cO_T$, that is, if $T\times_SS'\to T$ is an effective Cartier divisor. A sheaf $F$ on $(\Sch/S)_{\fppf}$ is \emph{$S'$-regular} if there is a surjective \'etale morphism $U\to F$, representable by schemes, from an $S'$-regular scheme $U$.
\leansrc{regular}

\begin{definition}[{\cite[\S13]{AlgMag}, \cite{MultiDil}}]\label{def:dilatation}
Let $X$ be a sheaf on $(\Sch/S)_{\fppf}$ and let $Y\to X_{S'}$ be a closed immersion. A \emph{dilatation} of $X$ with center $Y$ along $S'$ is an $S'$-regular sheaf $\Bl^{X_{S'}}_YX$ with a morphism $b\colon\Bl^{X_{S'}}_YX\to X$ such that the restriction of $b$ to $S'$ factors through $Y$, and which is universal with this property: for every $S'$-regular sheaf $T$ and every morphism $x\colon T\to X$ whose restriction to $S'$ factors through $Y$, there is a unique $\tilde x\colon T\to\Bl^{X_{S'}}_YX$ with $b\circ\tilde x=x$.
\leansrc{dilatation}
\end{definition}

\begin{remark}\label{rem:dilatation-presheaves}
In the Lean formalization the universal property quantifies over all $S'$-regular presheaves $T$ on $(\Sch/S)_{\fppf}$ (with $S'$-regularity defined by a surjective \'etale morphism, representable by schemes, from an $S'$-regular scheme), not only over sheaves.
\end{remark}

A dilatation is unique up to unique isomorphism over $X$. For an algebraic space $X$ over $S$ and a closed subspace $Y$ of $X_{S'}$, the dilatation exists, is an algebraic space and is affine over $X$ by \cite{MultiDil}, and by \cite{MRR} if $X$ is a scheme.
\leansrc{spDilatationExists}

\begin{proposition}[{\cite[Prop.~13.1]{AlgMag}}]\label{prop:dilatation}
Let $M$ be a monoid with $\#M$ below the cutoff of the site (for instance $M$ countable), let $X$ be an algebraic space over $S$ with an action of $\A(M)_S$, let $Y$ be a closed subspace of $X_{S'}$ stable under $\A(M)_{S'}$, and let $b\colon\Bl^{X_{S'}}_YX\to X$ be the dilatation. Then $\A(M)_S$ acts on $\Bl^{X_{S'}}_YX$ and $b$ is equivariant. Let $N$ be a magnet. Assume that $X^N$ is an algebraic space, that $Y^N$ and $(\Bl^{X_{S'}}_YX)^N$ are representable, and that $(\Bl^{X_{S'}}_YX)^N$ is flat over $S$, that is, there is a surjective \'etale morphism, representable by schemes, from a scheme flat over $S$ to $(\Bl^{X_{S'}}_YX)^N$; when $(\Bl^{X_{S'}}_YX)^N$ is an algebraic space, this is flatness over $S$ in the sense of \cite[Tag~03MM]{Stacks}. Then $(X_{S'})^N=(X^N)_{S'}$, the dilatation $\Bl^{(X^N)_{S'}}_{Y^N}X^N$ exists and is an algebraic space affine over $X^N$, $b^N$ is a dilatation of $X^N$ with center $Y^N$ along $S'$, and $b^N$ induces a canonical isomorphism over $X^N$
\[(\Bl^{X_{S'}}_YX)^N\xrightarrow{\ \sim\ }\Bl^{(X^N)_{S'}}_{Y^N}X^N.\]
\leansrc{dilatationAttractor}
\end{proposition}

\begin{proof}
The action and the comparison morphism are constructed in \cite[proof of Prop.~13.1]{AlgMag}, by the universal property of dilatations and Proposition~\ref{prop:attractor-limits}(iii); this construction applies to every monoid $M$. By Proposition~\ref{prop:attractor-limits}(iii), $(X_{S'})^N=(X^N)_{S'}$. The proof of Proposition~\ref{prop:immersions}(ii) applies to every monoid $M$ and gives $Y^N=(X_{S'})^N\times_{X_{S'}}Y$, so $Y^N\to(X^N)_{S'}$ is a closed immersion, as a base change of the closed immersion $Y\to X_{S'}$. Since $X^N$ is an algebraic space, the dilatation $\Bl^{(X^N)_{S'}}_{Y^N}X^N$ exists \cite{MultiDil}.

Let $T$ be an $S'$-regular scheme. The schemes $\A(N)_T$ and $\A(M)_T\times_T\A(N)_T$ are flat over $T$, hence $S'$-regular. By the universal property of $b$, composition with $b$ is a bijection from the set of morphisms $\A(N)_T\to\Bl^{X_{S'}}_YX$ onto the set of morphisms $\A(N)_T\to X$ whose restriction to $S'$ factors through $Y$. Let $\varphi\colon\A(N)_T\to\Bl^{X_{S'}}_YX$ be such that $b\circ\varphi$ is equivariant. The morphisms $(h,z)\mapsto\varphi(hz)$ and $(h,z)\mapsto h\cdot\varphi(z)$ from $\A(M)_T\times_T\A(N)_T$ to $\Bl^{X_{S'}}_YX$ have the same composite with $b$, since $b$ is equivariant, so they are equal by the uniqueness in the universal property; hence $\varphi$ is equivariant. Therefore the $T$-points of $(\Bl^{X_{S'}}_YX)^N$ correspond, via $b^N$, to the equivariant morphisms $\A(N)_T\to X_T$ whose restriction to $S'$ factors through $Y$, that is, to the $T$-points of $X^N$ whose restriction to $S'$ factors through $Y^N$. By Proposition~\ref{prop:attractor-limits}(iii), the restriction of $b^N$ to $S'$ is $(b_{S'})^N$, which factors through $Y^N$ since $b_{S'}$ factors through $Y$. The sheaf $(\Bl^{X_{S'}}_YX)^N$ is $S'$-regular. Indeed, let $V\to(\Bl^{X_{S'}}_YX)^N$ be a surjective \'etale morphism, representable by schemes, with $V$ a scheme flat over $S$, and let $A\to\Bl^{X_{S'}}_YX$ be a surjective \'etale morphism, representable by schemes, from an $S'$-regular scheme $A$. The scheme $P=V\times_{\Bl^{X_{S'}}_YX}A$, formed with the evaluation $V\to\Bl^{X_{S'}}_YX$ at the unit section, is \'etale and surjective over $V$, hence flat over $S$. For a local generator $s$ of the ideal of $S'$, flatness gives $\operatorname{Ann}_{\cO_P}(s)=\operatorname{Ann}_{\cO_S}(s)\cO_P$, and $\operatorname{Ann}_{\cO_S}(s)$ maps to $0$ in $\cO_A$, hence in $\cO_P$; so $P$ is $S'$-regular, and so is $V$, since $P\to V$ is faithfully flat. Let $T$ be an $S'$-regular sheaf and let $U\to T$ be a surjective \'etale morphism, representable by schemes, from an $S'$-regular scheme $U$. The scheme $U\times_TU$ is \'etale over $U$, hence flat over $U$, so it is $S'$-regular. The morphisms $U\to T$ and $U_{S'}\to T_{S'}$ are surjective \'etale and representable by schemes, hence epimorphisms of fppf sheaves, so $T$ is the coequalizer of the two projections $U\times_TU\rightrightarrows U$. Let $x\colon T\to X^N$ be a morphism whose restriction to $S'$ factors through $Y^N$. Its composite $x_U$ with $U\to T$ has the same property, so $x_U=b^N\circ\tilde x_U$ for a unique $\tilde x_U\colon U\to(\Bl^{X_{S'}}_YX)^N$. The two composites of $\tilde x_U$ with the projections $U\times_TU\to U$ lift the same morphism $U\times_TU\to X^N$, so they are equal by the uniqueness for the scheme $U\times_TU$; hence $\tilde x_U$ descends to a morphism $\tilde x\colon T\to(\Bl^{X_{S'}}_YX)^N$ with $b^N\circ\tilde x=x$, since $U\to T$ is an epimorphism. If $\tilde x'$ is another one, then $\tilde x$ and $\tilde x'$ agree on $U$ by the uniqueness for $U$, hence are equal. Conversely, for every $\tilde x\colon T\to(\Bl^{X_{S'}}_YX)^N$, the restriction of $b^N\circ\tilde x$ to $S'$ factors through $Y^N$. Hence $b^N$ satisfies the universal property of $\Bl^{(X^N)_{S'}}_{Y^N}X^N$.
\end{proof}

\begin{proposition}[{\cite[Fact~11.8, Lem.~11.9]{AlgMag}}]\label{prop:fact118}
Let $M$ be an abelian group, $N$ a submonoid of $M$ and $F$ a face of $N$. Let $\iota\colon D\to D_1$ be a surjective morphism of rings such that $xy=0$ for all $x,y\in\ker(\iota)$; put $T=\Spec(D)$ and $\bar T=\Spec(D_1)$. Let $R=D[F]\times_{D_1[F]}D_1[N]$, where $D[F]\to D_1[F]$ is induced by $\iota$ and $D_1[N]\to D_1[F]$ sends $X^n$ to $X^n$ for $n\in F$ and to $0$ otherwise, and let $\nu_1\colon\A(F)_T\to\Spec(R)$ and $\nu_2\colon\A(N)_{\bar T}\to\Spec(R)$ be given by the two projections.
\begin{enumerate}
\item $\Spec(R)$ is an affine scheme.
\item The morphism $\Spec(R)\to\A(N)_T$ given by $D[N]\to R$, $w\mapsto(\pi_F(w),\iota(w))$, with $\pi_F\colon D[N]\to D[F]$ as above, is a first order thickening.
\item Let $\Spec(R)$ carry the action of $\D(M)$ induced from $\A(N)_T$, and let $\Spec(R)\to S$ be a $\D(M)$-invariant morphism, $\A(F)_T$ and $\A(N)_{\bar T}$ being $S$-schemes through $\nu_1$ and $\nu_2$. Let $X$ be an algebraic space over $S$ with an action of $\D(M)_S$, and let $f\colon\A(F)_T\to X$ and $g\colon\A(N)_{\bar T}\to X$ be $\D(M)$-equivariant morphisms over $S$ with the same restriction to $\A(F)_{\bar T}$. Then there is a unique morphism $k\colon\Spec(R)\to X$ over $S$ with $k\circ\nu_1=f$ and $k\circ\nu_2=g$, and $k$ is $\D(M)$-equivariant.
\item Let $P$ be a presheaf on $(\Sch/S)_{\fppf}$ with a morphism $\pi\colon P\to T'$ to an $S$-scheme $T'$, let $\tau\colon\Spec(R)\to T'$ be a morphism over $S$ and let $k\in P(\Spec(R))$. If $\pi\circ k\circ\nu_1=\tau\circ\nu_1$ and $\pi\circ k\circ\nu_2=\tau\circ\nu_2$, then $\pi\circ k=\tau$.
\end{enumerate}
\leansrc{fact118}
\end{proposition}

\begin{proof}
\cite[proofs of Fact~11.8 and Lem.~11.9]{AlgMag}. The existence and uniqueness in (iii) hold without the actions; (iv) is this uniqueness for the algebraic space $T'$, applied to $\pi\circ k$ and $\tau$.
\end{proof}

\begin{corollary}[{\cite[Cor.~13.3]{AlgMag}}]\label{cor:dilatation-smooth}
Let $S'\subseteq S$ be an effective Cartier divisor and let $M$ be an abelian group. Let $X$ be a smooth $S$-scheme with an action of $\D(M)_S$ and let $Y$ be a closed subscheme of $X_{S'}$, smooth over $S'$ and stable under $\D(M)_{S'}$. Let $N$ be a magnet and assume that $X^N$, $Y^N$ and $(\Bl^{X_{S'}}_YX)^N$ are represented by schemes. Then $\D(M)_S$ acts on $\Bl^{X_{S'}}_YX$, the scheme $(\Bl^{X_{S'}}_YX)^N$ is smooth over $S$, and $(\Bl^{X_{S'}}_YX)^N=\Bl^{(X^N)_{S'}}_{Y^N}X^N$ canonically over $X^N$.
\leansrc{dilatationSmooth}
\end{corollary}

\begin{proof}
By \cite[Prop.~2.16(4)]{MRR}, $B=\Bl^{X_{S'}}_YX$ is smooth over $S$. The scheme $B^N$ is smooth over $S$; then it is flat over $S$, and Proposition~\ref{prop:dilatation} applies. Let $S$ carry the trivial action, so that $S^N=S$. By Proposition~\ref{prop:attractor-lfp}, applied to $B\to S$, the morphism $B^N\to S$ is locally of finite presentation, so by \cite[Tag~02H6]{Stacks} it suffices to show that it is formally smooth. Let $T_0\subseteq T$ be a first order thickening of affine $S$-schemes with ideal $J$, and let $\psi_0\in B^N(T_0)$. Put $D=\cO(T)[N]$ and $D_0=\cO(T_0)[N]$; the map $\phi\colon D\to D_0$ is surjective with kernel $J[N]$, of square zero. Since $B$ is smooth over $S$, the morphism $\psi_0\colon\Spec(D_0)\to B$ extends to a morphism $\alpha\colon\Spec(D)\to B$ over $S$ \cite[Tag~02H6]{Stacks}. Write the elements of $D[M]=D\otimes_\ZZ\ZZ[M]$ as $w=\sum_mw_mY^m$ with $w_m\in D$, so that $\Spec(D[M])=\D(M)_T\times_T\A(N)_T$, the variable $Y$ being the coordinate of $\D(M)$. Let $R=D_0\times_{D_0[M]}D[M]$, where $D_0\to D_0[M]$ is the inclusion of constants and $D[M]\to D_0[M]$ is induced by $\phi$; by \cite[Tag~0ET0]{Stacks}, $\Spec(R)$ is the pushout of $\Spec(D_0)\leftarrow\Spec(D_0[M])\to\Spec(D[M])$. Let $\iota_0\colon\Spec(D_0)\to\Spec(R)$ and $\iota_1\colon\Spec(D[M])\to\Spec(R)$ be the morphisms given by $(x,w)\mapsto x$ and $(x,w)\mapsto w$. The ring map $D\to R$, $d\mapsto(\phi(d),d)$, makes $\Spec(R)$ a scheme over $T$, and $\iota_0$, $\iota_1$ are morphisms over $T$. The morphism $\beta\colon\D(M)_T\times_T\A(N)_T\to B$, $(h,z)\mapsto h^{-1}\cdot\alpha(hz)$, restricts on $\D(M)_{T_0}\times_{T_0}\A(N)_{T_0}$ to $(h,z)\mapsto h^{-1}\cdot\psi_0(hz)=\psi_0(z)$, since $\psi_0$ is equivariant; so $\psi_0$ and $\beta$ define a morphism $\gamma\colon\Spec(R)\to B$. The composites of $\gamma$ and of $\Spec(R)\to T\to S$ with $\iota_0$ and $\iota_1$ agree, so $\gamma$ is a morphism over $S$ by the uniqueness in the pushout property.

For $(x,w)\in R$ the coefficients $w_m$, $m\neq0$, lie in $J[N]$ and $\phi(w_0)=x$. Since $J[N]^2=0$, the map $\varepsilon\colon R\to D$, $(x,w)\mapsto w_0$, is a ring homomorphism: the coefficient of $Y^0$ in $ww'$ is $\sum_mw_mw'_{-m}=w_0w'_0$. It satisfies $\varepsilon(\phi(d),d)=d$ for $d\in D$ and $\phi\circ\varepsilon(x,w)=x$, so $\Spec(\varepsilon)$ is a morphism over $T$ and $\Spec(\varepsilon)\circ\Spec(\phi)=\iota_0$. Put $\psi=\gamma\circ\Spec(\varepsilon)\colon\A(N)_T\to B$. It is a morphism over $S$, and $\psi\circ\Spec(\phi)=\gamma\circ\iota_0=\psi_0$.

The morphism $\psi$ is equivariant: $\psi(gz)=g\cdot\psi(z)$ for every $g\in\D(M)(T)$, that is, for every morphism of monoids $g\colon M\to\cO(T)^\times$, $m\mapsto g_m$; let $\bar g\in\D(M)(T_0)$ be its restriction. Let $t_g\colon\A(N)_T\to\A(N)_T$, $z\mapsto gz$, so that $t_g^\sharp(X^n)=g_nX^n$, and similarly $t_{\bar g}$ on $\A(N)_{T_0}$; then $\phi\circ t_g^\sharp=t_{\bar g}^\sharp\circ\phi$. Let $\tau_g\colon\D(M)_T\times_T\A(N)_T\to\D(M)_T\times_T\A(N)_T$, $(h,z)\mapsto(hg,z)$, so that $\tau_g^\sharp(\sum_mw_mY^m)=\sum_mg_mw_mY^m$. The maps
\[\rho(x,w)=\Bigl(t_{\bar g}^\sharp(x),\ \sum_mt_g^\sharp(w_m)Y^m\Bigr),\qquad\sigma(x,w)=\bigl(x,\ \tau_g^\sharp(w)\bigr)\]
are ring endomorphisms of $R$: the image of $\sum_mt_g^\sharp(w_m)Y^m$ in $D_0[M]$ is $t_{\bar g}^\sharp(\phi(w_0))=t_{\bar g}^\sharp(x)$, and that of $\tau_g^\sharp(w)$ is $g_0\phi(w_0)=x$. They satisfy
\[\Spec(\rho)\circ\iota_0=\iota_0\circ t_{\bar g},\quad\Spec(\rho)\circ\iota_1=\iota_1\circ(\id\times t_g),\quad\Spec(\sigma)\circ\iota_0=\iota_0,\quad\Spec(\sigma)\circ\iota_1=\iota_1\circ\tau_g,\]
and $\varepsilon\circ\rho=t_g^\sharp\circ\varepsilon$, $\varepsilon\circ\sigma=\varepsilon$, since $g_0=1$. Let $\gamma'=g^{-1}\cdot(\gamma\circ\Spec(\rho))$, the action of $g^{-1}$ on $B$ being taken over the $T$-scheme $\Spec(R)$, and $\gamma''=\gamma\circ\Spec(\sigma)$. On $\Spec(D_0)$,
\[\gamma'\circ\iota_0(z)=\bar g^{-1}\cdot\psi_0(\bar gz)=\psi_0(z)=\gamma''\circ\iota_0(z),\]
since $\psi_0$ is equivariant; on $\Spec(D[M])$,
\[\gamma'\circ\iota_1(h,z)=g^{-1}\cdot\beta(h,gz)=g^{-1}h^{-1}\cdot\alpha(hgz)=(hg)^{-1}\cdot\alpha(hgz)=\beta(hg,z)=\gamma''\circ\iota_1(h,z),\]
since $\D(M)$ is commutative. By the uniqueness in the pushout property, $\gamma'=\gamma''$. Composing with $\Spec(\varepsilon)$, and using $\Spec(\rho)\circ\Spec(\varepsilon)=\Spec(\varepsilon)\circ t_g$ and $\Spec(\sigma)\circ\Spec(\varepsilon)=\Spec(\varepsilon)$, gives $g^{-1}\cdot\psi(gz)=\psi(z)$.

The construction is compatible with base change along a morphism of affine schemes $T'\to T$: with $T'_0=T_0\times_TT'$, the lift $\alpha_{T'}$ of $(\psi_0)_{T'}$ and the ring $R'$ obtained from $D'=\cO(T')[N]$ and $D'_0=\cO(T'_0)[N]$, the composite of $\gamma$ with $\Spec(R')\to\Spec(R)$ agrees with the morphism $\gamma_{T'}$ obtained from $\alpha_{T'}$ on both pieces of the pushout, so the resulting morphism $\A(N)_{T'}\to B$ is $\psi_{T'}$. Applying the above to $T'=\D(M)_T$ and its universal point $g\in\D(M)(T')$ gives $\psi(gz)=g\cdot\psi(z)$ for all points $g$ of $\D(M)$ and $z$ of $\A(N)_T$; so $\psi$ is equivariant, and $\psi\in B^N(T)$ extends $\psi_0$.
\end{proof}

\begin{corollary}\label{cor:dilatation-smooth-affine}
In Corollary~\ref{cor:dilatation-smooth}, assume moreover that $X$ is affine over $S$, and drop the hypothesis on $X^N$, $Y^N$ and $(\Bl^{X_{S'}}_YX)^N$; $M$ is any abelian group. Then $X^N$, $Y^N$ and $(\Bl^{X_{S'}}_YX)^N$ are represented by schemes, $(\Bl^{X_{S'}}_YX)^N$ is smooth over $S$, and $b^N$ is a dilatation of $X^N$ with center $Y^N$ along $S'$.
\leansrc{dilatationSmoothAffine}
\end{corollary}

\begin{proof}
The morphisms $Y\to S$ and $B=\Bl^{X_{S'}}_YX\to S$ are affine, $Y\to X$ being a closed immersion and $B\to X$ affine \cite{MRR}. So $X^N$, $Y^N$ and $B^N$ are represented by schemes affine over $S$ (Theorem~\ref{thm:affine}), and Corollary~\ref{cor:dilatation-smooth} applies; the dilatation of $X^N$ exists by \cite{MRR}.
\end{proof}

\begin{corollary}[{\cite[Cor.~13.3]{AlgMag}}]\label{cor:dilatation-smooth-spaces}
Let $S'\subseteq S$ be an effective Cartier divisor and let $M$ be a finitely generated abelian group. Let $X$ be a smooth and separated $S$-scheme with an action of $\D(M)_S$, let $Y$ be a closed subscheme of $X_{S'}$, smooth over $S'$ and stable under $\D(M)_{S'}$, and let $b\colon\Bl^{X_{S'}}_YX\to X$ be the dilatation. Let $N$ be a magnet. Then $X^N$ and $(\Bl^{X_{S'}}_YX)^N$ are algebraic spaces, $(\Bl^{X_{S'}}_YX)^N$ is smooth over $S$, and $b^N\colon(\Bl^{X_{S'}}_YX)^N\to X^N$ is a dilatation of $X^N$ with center $Y^N$ along $S'$ (Definition~\ref{def:dilatation}); in particular $(\Bl^{X_{S'}}_YX)^N=\Bl^{(X^N)_{S'}}_{Y^N}X^N$ canonically over $X^N$.
\leansrc{dilatationSmoothSpaces}
\end{corollary}

\begin{proof}
The schemes $X$ and $B=\Bl^{X_{S'}}_YX$ are smooth over $S$ \cite[Prop.~2.16(4)]{MRR}, hence locally of finite presentation, and $B$ is separated over $S$, being affine over $X$. For a separated algebraic space $W$ over $S$ locally of finite presentation with an action of $\D(M)_S$, $W^N$ is an algebraic space: over an affine open subscheme $V\subseteq S$, which is separated, Corollary~\ref{cor:representability} applies to $W_V$, and $(W^N)_V=(W_V)^N$ (Proposition~\ref{prop:attractor-limits}(iii)); the $(W^N)_V$ are open subfunctors of the fppf sheaf $W^N$ which are algebraic spaces and cover it, so $W^N$ is an algebraic space \cite[Tag~02WR]{Stacks}. So $X^N$ and $B^N$ are algebraic spaces. The restrictions $(W^N)_V$ above are algebraic spaces, and their disjoint union is an algebraic space \cite[Tag~02WQ]{Stacks}; here $V\to S$ is an open immersion and $(W^N)_V=W^N\times_SV$ is the base change of $W^N$ along it (Proposition~\ref{prop:attractor-limits}(iii)). Also $Y\to X_{S'}$ is a closed immersion, so $Y^N=X^N\times_XY$ (Proposition~\ref{prop:immersions}(ii), whose proof applies to every monoid $M$) is representable, as a base change of the algebraic space $X^N$ along the closed immersion $Y\to X$. The argument of the proof of Corollary~\ref{cor:dilatation-smooth} shows that $B^N\to S$ is locally of finite presentation and formally smooth; it uses only the points of $B^N$ and applies to the algebraic space $B^N$, which is therefore smooth over $S$ \cite[Tag~04AM]{Stacks}. Then $B^N$ is flat over $S$, and Proposition~\ref{prop:dilatation} applies.
\end{proof}

\section{Ind-algebraic spaces}\label{sec:ind}

Let $(\mathrm{Aff}/S)$ be the category of affine schemes over $S$.

\begin{definition}[{\cite[Def.~14.1, Rem.~14.2]{AlgMag}}]\label{def:ind}
An \emph{ind-algebraic space} (resp.\ \emph{ind-scheme}) over $S$ is a functor $X\colon(\mathrm{Aff}/S)^{\mathrm{op}}\to\mathrm{Set}$ which admits a presentation $X=\colim_{i\in I}X_i$ as a filtered colimit of algebraic spaces over $S$ (resp.\ $S$-schemes) whose transition maps $X_i\to X_j$, $i\leq j$, are closed immersions. A morphism $X\to Y$ of such functors is \emph{representable by algebraic spaces} (resp.\ \emph{schemes}, \emph{closed immersions}) if for every affine scheme $T$ and every $T\to Y$, the fibre product $T\times_YX$ is an algebraic space (resp.\ a scheme, a closed subscheme of $T$). An ind-algebraic space $X$ is \emph{separated} if its diagonal is representable by closed immersions: for every affine $T$ and $x,y\in X(T)$ there is a closed subscheme $W\subseteq T$ such that a morphism of affine schemes $T'\to T$ factors through $W$ if and only if $x$ and $y$ agree on $T'$.
\leansrc{ind}
\end{definition}

\begin{definition}[{\cite[Def.~14.3]{AlgMag}}]\label{def:ind-action}
Let $G$ be a group scheme over $S$ and $X=\colim_iX_i$ an ind-algebraic space. An \emph{ind-action} of $G$ on $X$ is a family of actions $a_i\colon G\times_SX_i\to X_i$ compatible with the transition maps; it induces an action of $G$ on $X$ as a functor.
\leansrc{indAction}
\end{definition}

\begin{proposition}[{\cite[Prop.~14.4]{AlgMag}}]\label{prop:ind}
Let $M$ be a finitely generated abelian group and let $X=\colim_iX_i$ be an ind-algebraic space (resp.\ ind-scheme) over $S$ with an ind-action of $\D(M)_S$, such that every $X_i$ is separated and locally of finite presentation over $S$ (if $X$ is separated, every $X_i$ is separated). Let $N$ be a magnet. Then:
\begin{enumerate}
\item the attractor functor $T\mapsto\Homeq{\D(M)_T}_T(\A(N)_T,X_T)$ on $(\mathrm{Aff}/S)$ is the ind-algebraic space (resp.\ ind-scheme) $\colim_iX_i^N$;
\item the natural morphism $X^N\to X$ is representable by algebraic spaces (resp.\ schemes);
\item if $N=Z$ is a subgroup, then $X^Z\to X$ is representable by closed immersions.
\end{enumerate}
\leansrc{indAttractor}
\end{proposition}

\begin{proof}
For affine $T$, every morphism $T\to X$ factors through some $X_i$, and then $T\times_XX^N=T\times_{X_i}X_i^N$ \cite[proof of Prop.~14.4]{AlgMag}; the transition maps $X_i^N\to X_j^N$ are closed immersions by Proposition~\ref{prop:immersions}(ii). So (i) and (ii) follow once every $X_i^N$ is an algebraic space (resp.\ a scheme). If $X$ is separated, the terms $X_i$ are separated, since $X_i(T)\to X(T)$ is injective.

The question is Zariski local on $S$. Let $V\subseteq S$ be an affine open subscheme. By Proposition~\ref{prop:attractor-limits}(iii), $X_i^N\times_SV=((X_i)_V)^N$, and $(X_i)_V$ is a separated $V$-space locally of finite presentation over the separated scheme $V$; so $X_i^N\times_SV$ is an algebraic space by Corollary~\ref{cor:representability}. Since $X_i^N$ is a sheaf (Proposition~\ref{prop:attractor-sheaf}) whose restrictions to the members of an affine open cover of $S$ are algebraic spaces, and their disjoint union is an algebraic space \cite[Tag~02WQ]{Stacks} (the size condition there holds by the choice of $\kappa$), $X_i^N$ is an algebraic space \cite[Tag~04SK]{Stacks}. Moreover $(X_i)_V$ has an $\units N$-FPR atlas (Theorem~\ref{thm:atlas}), so by Theorem~\ref{thm:representability}(b) the subspace $((X_i)_V)^{\units N}$ of $(X_i)_V$ is closed and the morphism $((X_i)_V)^N\to((X_i)_V)^{\units N}$ is affine. These two properties are Zariski local on the target, so $X_i^{\units N}\to X_i$ is a closed immersion and $X_i^N\to X_i^{\units N}$ is affine. If $X_i$ is a scheme, $X_i^N$ is therefore affine over a closed subscheme of $X_i$, hence a scheme.

For (iii), apply this to $N=Z=\units Z$: $X_i^Z\to X_i$ is a closed immersion, so $T\times_XX^Z=T\times_{X_i}X_i^Z$ is a closed subspace of $T$. No atlas hypothesis is needed (Theorem~\ref{thm:atlas} provides one).
\end{proof}

\section{Tangent spaces and Lie algebras}\label{sec:tangent}

Let $I_S=\Spec(\ZZ[\varepsilon])\times_{\Spec(\ZZ)}S$, with $\ZZ[\varepsilon]=\ZZ[t]/(t^2)$, be the scheme of dual numbers over $S$; for a scheme $T$ over $S$, $I_T=I_S\times_ST$.

\begin{definition}[{\cite[\S16.1]{AlgMag}}]\label{def:tangent}
Let $X$ be a functor on schemes over $S$. The \emph{tangent space} $TX$ is the functor $T\mapsto\Hom_T(I_T,X_T)=X(I_T)$.
\leansrc{tangent}
\end{definition}

\begin{proposition}[{\cite[Rem.~16.1]{AlgMag}}]\label{prop:tangent-basic}
For morphisms $X\to Y\leftarrow Z$ of functors, $T(X\times_YZ)=TX\times_{TY}TZ$. If $X\to Y$ is a monomorphism, so is $TX\to TY$. Moreover $TS=S$.
\leansrc{tangentBasic}
\end{proposition}

\begin{proof}
By definition $(TX)(T)=X(I_T)$, functorially in $X$, and limits of functors are computed pointwise: $T(X\times_YZ)(T)=X(I_T)\times_{Y(I_T)}Z(I_T)=(TX\times_{TY}TZ)(T)$; if $X\to Y$ is a monomorphism, $TX(T)=X(I_T)\to Y(I_T)=TY(T)$ is injective for every $T$; and $TS(T)=\Hom_S(I_T,S)$ is a point, $S$ being the final object. This is \cite[Rem.~16.1]{AlgMag}.
\end{proof}

Let $M$ be a commutative monoid and let $\A(M)_S$ act on $X$. Then $\A(M)_S$ acts on $TX$, by $h\cdot\xi=h_{I_T}\cdot\xi$ for $h\in\A(M)(T)$ and $\xi\in X(I_T)$.

\begin{proposition}[{\cite[Prop.~16.2]{AlgMag}}]\label{prop:tangent}
For every magnet $N$, $(TX)^N=T(X^N)$ canonically.
\leansrc{tangentAttractor}
\end{proposition}

\begin{proof}
\cite[proof of Prop.~16.2]{AlgMag}
\end{proof}

Let $G$ be a group functor over $S$. Then $TG$ is a group functor, with morphisms of group functors $G\to TG$ and $TG\to G$ induced by $I_T\to T$ and by the zero section $T\to I_T$. The \emph{Lie algebra} of $G$ is the group functor $\Lie(G)=e_S\times_GTG$, the fibre product of the unit $e_S\to G$ and of $TG\to G$.
\leansrc{lie}

\begin{proposition}[{\cite[Rem.~16.3]{AlgMag}}]\label{prop:lie-basic}
For morphisms of group functors $G\to K\leftarrow H$, $\Lie(G\times_KH)=\Lie(G)\times_{\Lie(K)}\Lie(H)$. If $G\to H$ is a monomorphism, so is $\Lie(G)\to\Lie(H)$. Moreover $\Lie(e_S)=e_S$.
\leansrc{lieBasic}
\end{proposition}

\begin{proof}
By definition $\Lie(G)(T)$ is the kernel of $TG(T)=G(I_T)\to G(T)$, and limits are computed pointwise. By Proposition~\ref{prop:tangent-basic}, $T(G\times_KH)=TG\times_{TK}TH$, and taking fibres over the unit sections, which commutes with fibre products, gives $\Lie(G\times_KH)=\Lie(G)\times_{\Lie(K)}\Lie(H)$; $\Lie(G)\to\Lie(H)$ is the restriction of the monomorphism $TG\to TH$; and $\Lie(e_S)=e_S\times_{e_S}T(e_S)=e_S$, since $e_S=S$ and $TS=S$. This is \cite[Rem.~16.3]{AlgMag}.
\end{proof}

Assume that $\A(M)_S$ acts on $G$ by group endomorphisms. Then it acts on $TG$ and on $\Lie(G)$ by group endomorphisms.

\begin{proposition}[{\cite[Prop.~16.4, Prop.~16.5]{AlgMag}}]\label{prop:lie}
Let $N$ be a magnet.
\begin{enumerate}
\item $\Lie(G)^N=\Lie(G^N)$ as group functors.
\item Let $F$ be a face of $N$ and $H\to G^F$ a morphism of group functors. Then $\Lie(G)^N_{F,\Lie(H)}=\Lie(G^N_{F,H})$ as group functors.
\end{enumerate}
\leansrc{lieAttractor}
\end{proposition}

\begin{proof}
\cite[proofs of Prop.~16.4 and Prop.~16.5]{AlgMag}, where the face $F$ is $\units N$ and $H\to G^F$ is a monomorphism; the same proof applies to any face and any morphism $H\to G^F$, since $G^N_{F,H}=G^N\times_{G^F}H$ (Proposition~\ref{prop:prescribed}), $\Lie$ commutes with fibre products (Proposition~\ref{prop:lie-basic}) and $\Lie(G)^N=\Lie(G^N)$, $\Lie(G)^F=\Lie(G^F)$ by (i).
\end{proof}

\begin{proposition}[{\cite[Prop.~16.6]{AlgMag}}]\label{prop:semidirect}
Let $M$ be a commutative monoid and $G$ a group functor over $S$ on which $\A(M)_S$ acts by group endomorphisms. Let $N$ be a magnet and $F$ a face of $N$. Then $G^N$ is the semidirect product of $G^N_{F,e_G}$ and $G^F$; in particular the multiplication $G^N_{F,e_G}\times G^F\to G^N$ is an isomorphism of functors.
\leansrc{semidirect}
\end{proposition}

\begin{proof}
\cite[proof of Prop.~16.6]{AlgMag}, where $F=\units N$: the morphism $p_{N,F}\colon G^N\to G^F$ is a retraction of $\iota_{F,N}$ with kernel $G^N_{F,e_G}$; for every scheme $T$ over $S$ and $g\in G^N(T)$ one has $g=\bigl(g\,\iota(p(g))^{-1}\bigr)\iota(p(g))$ with the first factor in the kernel, uniquely since $p\iota=\mathrm{id}$, which gives the semidirect product decomposition \cite[2.34]{Milne}; the same argument applies to any face $F$.
\end{proof}

\begin{proposition}[{\cite[Prop.~16.15(ii)]{AlgMag}}]\label{prop:root-smooth}
Let $M$ be a finitely generated abelian group. Let $G$ be a group algebraic space over $S$, smooth and separated over $S$, on which $\D(M)_S$ acts by group endomorphisms, and let $\alpha\in M$. Then $G^{\mon{\alpha}}$ is represented by an algebraic space smooth over $S$, and $G^{\mon{\alpha}}_{e_S}$ is isomorphic, as a functor, to an algebraic space smooth over $S$.
\leansrc{rootSmooth}
\end{proposition}

\begin{proof}
The space $G$ is locally of finite presentation, being smooth. For $X$ a separated algebraic space over $S$ locally of finite presentation with an action of $\D(M)_S$ and $N$ a magnet, $X^N$ is an algebraic space: over an affine open subscheme $V\subseteq S$, which is separated, Corollary~\ref{cor:representability} applies to $X_V$, and $(X^N)_V=(X_V)^N$ (Proposition~\ref{prop:attractor-limits}(iii)); the $(X^N)_V$ are open subfunctors of the fppf sheaf $X^N$ which are algebraic spaces and cover it, so $X^N$ is an algebraic space \cite[Tag~02WR]{Stacks}. Hence $G^{\mon{\alpha}}$ and $G^{\units\mon{\alpha}}$ are algebraic spaces; $G^{\mon{\alpha}}_{e_S}$ is the fibre product of $G^{\mon{\alpha}}\to G^{\units\mon{\alpha}}$ and of the unit section (Proposition~\ref{prop:prescribed}), hence an algebraic space. With the trivial action on $S$, $G^{\mon{\alpha}}\to S^{\mon{\alpha}}=S$ is smooth by Corollary~\ref{cor:attractor-smooth}, and $G^{\mon{\alpha}}_{e_S}$ is a base change of $p_{\mon{\alpha},\units\mon{\alpha}}\colon G^{\mon{\alpha}}\to G^{\units\mon{\alpha}}$, which is smooth by Corollary~\ref{cor:face-smooth}, $\units\mon{\alpha}$ being a face of $\mon{\alpha}$.
\end{proof}

\section{Pure magnets}\label{sec:pure}

In this section $M$ is a monoid and $X$ is a sheaf on $(\Sch/S)_{\fppf}$ with an action $a$ of $\A(M)_S$. For a magnet $N$, let $[X^N]\subseteq X$ be the image of the morphism $X^N\to X$, $\varphi\mapsto\varphi(1)$, that is, the subfunctor of $X$ whose $T$-points are the $\varphi(1)$, $\varphi\in X^N(T)$. If $X$ is a separated algebraic space over $S$, $X^N\to[X^N]$ is an isomorphism (Proposition~\ref{prop:separated-mono}). For magnets $N\subseteq L$ one has $[X^N]\subseteq[X^L]$, since $\varphi\in X^N(T)$ gives $\varphi\circ q\in X^L(T)$ for the morphism $q\colon\A(L)_T\to\A(N)_T$ of Proposition~\ref{prop:action-AN}(ii), and $q(1)=1$.
\leansrc{pureMonotone}

\subsection{Pure magnets and attractors}

\begin{definition}[{\cite[Def.~15.2]{AlgMag}}]\label{def:pure}
A magnet $N$ is \emph{pure} for $a$ if every magnet $L$ with $[X^L]=[X^N]$ contains $N$. The set of pure magnets is denoted $\mho(a)$. An \emph{attractor} of $a$ is a subfunctor of $X$ of the form $[X^N]$ for a magnet $N$.
\leansrc{pure}
\end{definition}

\begin{lemma}[{\cite[Fact~3.34]{AlgMag}}]\label{lem:fact334}
Assume that one of the following conditions holds:
\begin{enumerate}[label=\textup{(\alph*)}]
\item $X$ is $S$-affine;
\item $M$ is cancellative, $M^{\gp}$ is finitely generated, and $X$ is a separated algebraic space over $S$ locally of finite presentation.
\end{enumerate}
Let $(N_i)_{i\in I}$ be a nonempty family of magnets and $Y\subseteq X$ a subfunctor with $[X^{N_i}]=Y$ for every $i$. Then $[X^{\bigcap_iN_i}]=Y$.
\leansrc{fact334}
\end{lemma}

\begin{proof}
(a) Locally on $S$, $X=\Spec(A)$ with $A=\bigoplus_{m\in M}A_m$, and $X^N$ is the closed subscheme defined by the ideal $J_N$ generated by the $A_m$, $m\notin N$ (Theorem~\ref{thm:affine}); then $J_{\bigcap_iN_i}=\sum_iJ_{N_i}$, which is \cite[proof of Fact~3.34(i)]{AlgMag}. (b) By Proposition~\ref{prop:change-monoid}(iii), the attractors of $X$ for $M$ and for $M^{\gp}$ agree on the magnets of $M$, and an intersection of magnets of $M$ is the same in $M$ and in $M^{\gp}$; so we may assume that $M$ is a finitely generated abelian group. For a set $J$ of indices put $N_J=\bigcap_{i\in J}N_i$, and $N_\infty=N_I$. The inclusion $[X^{N_\infty}]\subseteq Y$ holds since $N_\infty\subseteq N_i$. Conversely, let $T$ be a scheme over $S$ and $x\in Y(T)$. Since $X^{N_\infty}\to X$ is a monomorphism (Proposition~\ref{prop:separated-mono}) and $X^{N_\infty}$ is a sheaf (Proposition~\ref{prop:attractor-sheaf}), the local lifts of $x$ to $X^{N_\infty}$ are unique and glue; so it suffices to lift $x$ locally on $T$ for the \'etale topology. In particular we may assume that $T$ is affine and that $T\to S$ factors through an affine open subscheme $V\subseteq S$; by Proposition~\ref{prop:attractor-limits}(iii), replacing $S$ by $V$ and $X$ by $X_V$, we may assume $S$ affine, hence separated.

\emph{Finite families.} Let $N,L$ be magnets with $[X^N]=[X^L]$ and $\psi\in X^L(T')$. For a morphism $t\colon T''\to T'$ and $b\in\A(L)(T'')$, the point $\psi(b)$ is the value at $1$ of the point $b\cdot\psi_{T''}$ of $X^L(T'')$ (Proposition~\ref{prop:attractor-actions}(i)), so it lies in $[X^L](T'')=[X^N](T'')$: $\psi(b)=G(b)(1)$ for a unique $G(b)\in X^N(T'')$ (Proposition~\ref{prop:separated-mono}). By uniqueness, $b\mapsto G(b)$ is natural in $T''$ and equivariant, so $G$ is a point of $(X^N)^L(T')$ with $G(1)(1)=\psi(1)$. As $(X^N)^L=X^{N\cap L}$ (Proposition~\ref{prop:intersection}(ii)), $\psi(1)\in[X^{N\cap L}](T')$. Hence $[X^{N\cap L}]=[X^N]$, and by induction $[X^{N_J}]=Y$ for every finite nonempty $J\subseteq I$.

\emph{Choice of a vertex.} For a magnet $L$ let $e_L\in\A(L)(T)$ be the point $n\mapsto1$ for $n\in\units L$, $n\mapsto0$ otherwise, the image of $1$ under $\A(\units L)_T\to\A(L)_T$; for $\psi\in X^L(T)$, $\psi(e_L)$ is the value at $1$ of $p_{L,\units L}(\psi)$. Choose a finite nonempty $J_0\subseteq I$ such that the rank of the finitely generated abelian group $\units{N_{J_0}}$ is minimal. For finite $J\supseteq J_0$, $\units{N_J}\subseteq\units{N_{J_0}}$ have the same rank, so for $n\in N_J\cap\units{N_{J_0}}$ some $dn$, $d\geq1$, lies in $\units{N_J}$, and $-n=(d-1)n+(-dn)\in N_J$; thus $N_J\cap\units{N_{J_0}}=\units{N_J}$, and $e_{N_{J_0}}$ restricts to $e_{N_J}$ along $N_J\subseteq N_{J_0}$. Hence if $\psi\in X^{N_J}(T)$ has image $\psi_0=\iota_{N_J,N_{J_0}}(\psi)$, then $\psi(e_{N_J})=\psi_0(e_{N_{J_0}})$.

\emph{Charts.} By the finite case, $x=\psi_0(1)$ with $\psi_0\in X^{N_{J_0}}(T)$; put $v=\psi_0(e_{N_{J_0}})\in X(T)$. Let $U=\coprod_aU_a\to X$ be a strongly-FPR atlas (Theorem~\ref{thm:strong-atlas}). As $T$ is affine, finitely many of the \'etale morphisms $T_a=T\times_{v,X}U_a\to T$, say $a\in F$, cover $T$. Each $U_a$ is affine and of finite type over the affine $S$, so $\cO(U_a)$ is generated over $\cO(S)$ by finitely many homogeneous elements, with degrees in a finite set $E_a\subseteq M$, and $U_a^L=U_a^{\mon{L\cap E_a}}$ for every magnet $L$ (Theorem~\ref{thm:affine}). Put $E=\bigcup_{a\in F}E_a$. The sets $N_i\cap E$, $i\in I$, are finitely many; choose $i_1,\dots,i_r$ realizing all of them and put $J=J_0\cup\{i_1,\dots,i_r\}$. Then $N_J\cap E\subseteq N_i\cap E$ for every $i$, so $\mon{N_J\cap E_a}\subseteq N_\infty$ for $a\in F$. By the finite case, $x=\psi(1)$ with $\psi\in X^{N_J}(T)$; its image in $X^{N_{J_0}}(T)$ is $\psi_0$ (Proposition~\ref{prop:separated-mono}), so $\psi(e_{N_J})=v$. Let $a\in F$, $Q=\units{N_J}$, and let $u\in U_a(T_a)$ be the tautological lift of $v_{T_a}$. The point $v=p_{N_J,Q}(\psi)(1)$ is fixed by $\D(M/Q)$, and $U_a\to X$ is $Q$-FPR, so $u$ is a fixed point of $\D(M/Q)$, that is, a point of $U_a^Q(T_a)$ over $p_{N_J,Q}(\psi)_{T_a}$ (Proposition~\ref{prop:subgroup-fixed}). By Proposition~\ref{prop:etale-face}(ii), $\psi_{T_a}$ lifts to a point $w$ of $U_a^{N_J}(T_a)$, and $w=\iota(w')$ with $w'\in U_a^{\mon{N_J\cap E_a}}(T_a)$. The image of $w'$ in $X^{\mon{N_J\cap E_a}}(T_a)$, pushed to $X^{N_\infty}(T_a)$ by $\iota$, has value at $1$ the image of $w(1)$ in $X$, namely $x_{T_a}$. So $x$ lifts to $X^{N_\infty}$ \'etale locally on $T$, which proves $[X^{N_\infty}]=Y$.
\end{proof}

\begin{proposition}[{\cite[Prop.~15.1]{AlgMag}}]\label{prop:pure-lt}
Let $N$ be a magnet. If $N$ is pure, then $[X^L]\subsetneq[X^N]$ for every magnet $L\subsetneq N$. Under one of the conditions (a), (b) of Lemma~\ref{lem:fact334}, or if $M$ is a finitely generated abelian group and $X$ is a separated algebraic space over $S$, the converse holds.
\leansrc{pureLt}
\end{proposition}

\begin{proof}
\cite[proof of Prop.~15.1]{AlgMag}. The converse uses the equality $[X^{N\cap L}]=[X^N]$ for $[X^L]=[X^N]$; under (a) or (b) it is Lemma~\ref{lem:fact334}, and for a separated algebraic space over $S$ with $M$ a finitely generated abelian group it follows from Proposition~\ref{prop:intersection}(ii).
\end{proof}

For a subfunctor $Y\subseteq X$ put $E(Y)=\bigcap\{N\mid [X^N]=Y\}$.

\begin{theorem}[{\cite[Thm.~15.3]{AlgMag}}]\label{thm:pure}
Under one of the conditions (a), (b) of Lemma~\ref{lem:fact334}, the map $N\mapsto[X^N]$ is a bijection from $\mho(a)$ onto the set of attractors of $a$, with inverse $Y\mapsto E(Y)$.
\leansrc{pureBijection}
\end{theorem}

\begin{proof}
\cite[proof of Thm.~15.3]{AlgMag}: if $Y=[X^{N_0}]$, then $[X^{E(Y)}]=Y$ by Lemma~\ref{lem:fact334} applied to the magnets $N$ with $[X^N]=Y$, so $E(Y)$ is pure; and $E([X^N])=N$ for $N$ pure, by definition.
\end{proof}

\subsection{Finiteness}

\begin{theorem}[{\cite[Thm.~15.4, Rem.~15.6]{AlgMag}}]\label{thm:pure-finite}
Assume that $S$ is quasi-compact.
\begin{enumerate}
\item Let $U$ be an $S$-affine scheme of finite type over $S$ with an action of $\A(M)_S$. Then $\mho(a)$ and the set of attractors of $U$ are finite.
\item Let $M$ be an abelian group and let $X$ be a separated algebraic space over $S$ locally of finite presentation, with an action of $\D(M)_S$ and an $S$-affine strongly-FPR atlas. Then $\mho(a)$ and the set of attractors of $X$ are finite. In particular the set of the $X^Z$, $Z\subseteq M$ a subgroup, is finite; for $M$ finitely generated, $X^Z=X^{\D(M/Z)_S}$.
\item Let $M$ be a finitely generated abelian group and let $X$ be a separated algebraic space over $S$ of finite presentation with a Zariski locally linearizable action of $\D(M)_S$. Then $\mho(a)$ and the set of attractors of $X$ are finite.
\end{enumerate}
\leansrc{pureFinite}
\end{theorem}

\begin{proof}
(i) Let $S=\bigcup_{k=1}^rV_k$ be a finite affine open cover and $U\times_SV_k=\Spec(A_k)$. Each $A_k$ is generated over $\cO(V_k)$ by finitely many homogeneous elements, whose degrees form a finite set $E_k\subseteq M$; put $E=\bigcup_kE_k$. By Theorem~\ref{thm:affine}, $U^N=U^{\mon{N\cap E}}$ for every magnet $N$, so a pure magnet is generated by a subset of the finite set $E$, and the attractors are the $U^{\mon{E'}}$, $E'\subseteq E$. This is \cite[proof of Thm.~15.4(i)]{AlgMag}. (ii) By Proposition~\ref{prop:pure-reflect} below, $\mho(a_X)\subseteq\mho(a_U)$ for the atlas $U\to X$, which is finite by (i), applied to $U$ (of finite type over $S$, being \'etale over $X$). Every attractor of $X$ is $[X^N]$ with $N\in\mho(a_U)$: given $N$, the pure magnet $N'=E(U^N)$ of $U$ satisfies $U^{N'}=U^N$, hence $X^{N'}=X^N$ by Proposition~\ref{prop:pure-reflect}. The last assertion follows, the $X^Z$ being attractors, and $X^Z=X^{\D(M/Z)_S}$ for $M$ finitely generated by Proposition~\ref{prop:subgroup-fixed}. (iii) Since $S$ and $X$ are quasi-compact, $|X|$ is quasi-compact, so finitely many $\D(M)_S$-stable $S$-affine opens cover $X$, and their disjoint union is an $S$-affine strongly-FPR atlas (Proposition~\ref{prop:zll-atlas}); apply (ii).
\end{proof}

\begin{proposition}[{\cite[Prop.~15.5]{AlgMag}}]\label{prop:pure-reflect}
Let $M$ be an abelian group, $X$ a separated algebraic space over $S$ locally of finite presentation with an action $a_X$ of $\D(M)_S$, and $f\colon U\to X$ an $S$-affine strongly-FPR atlas, with action $a_U$ on $U$. For magnets $N,L$, if $U^N=U^L$, then $X^N=X^L$. In particular $\mho(a_X)\subseteq\mho(a_U)$.
\leansrc{pureReflect}
\end{proposition}

\begin{proof}
We may assume $N\subseteq L$, replacing $L$ by $N\cap L$ (Lemma~\ref{lem:fact334}(a) for $U$). Let $\psi\in X^L(T)$, $Z=\units L$ and $v=p_{L,Z}(\psi)\in X^Z(T)$. The point $v(1)$ of $X$ is fixed by $\D(M/Z)_T$, which acts trivially on $\D(Z)_T$. Locally for the \'etale topology on $T$, say on $T'\to T$, $v(1)$ lifts to $u\in U(T')$, and $u$ is fixed by $\D(M/Z)$ since $f$ is $Z$-FPR. The point $u$ defines an equivariant morphism $\D(Z)_{T'}\to U_{T'}$ with value $u$ at $1$: the morphism $\D(M)_{T'}\to U_{T'}$, $h\mapsto h\cdot u$, is invariant under translation by $\D(M/Z)$, and $\D(M)\to\D(Z)$ is faithfully flat and affine, $\ZZ[M]$ being free over $\ZZ[Z]$, with $\D(M)\times_{\D(Z)}\D(M)=\D(M)\times\D(M/Z)$; as $U$ is a scheme, it is a sheaf for the fpqc topology \cite[Tag~023Q]{Stacks}, so this morphism descends along $\D(M)_{T'}\to\D(Z)_{T'}$. This gives a point $u'\in U^Z(T')$ whose image in $X^Z(T')$ is $v_{T'}$, since both have value $v(1)_{T'}$ at $1$ and $X^Z\to X$ is a monomorphism (Proposition~\ref{prop:separated-mono}). By Proposition~\ref{prop:etale-face}(ii), $\psi_{T'}$ lifts to $U^L(T')=U^N(T')$. Its image in $X$ is therefore a point of $X^N$; these points of $X^N$ glue, $X^N\to X$ being a monomorphism. Hence $[X^L]\subseteq[X^N]$, and the reverse inclusion holds since $N\subseteq L$. The last assertion follows.
\end{proof}

\begin{corollary}\label{cor:conj157}
Assume that $S$ is separated and quasi-compact and that $M$ is a finitely generated abelian group. Let $X$ be a separated algebraic space over $S$ of finite presentation with an action $a$ of $\D(M)_S$. Then $\mho(a)$ and the set of attractors of $X$ are finite.
\leansrc{conj157Sep}
\end{corollary}

\begin{proof}
Since $S$ and $X$ are quasi-compact, $|X|$ is quasi-compact, so $X$ has an $S$-affine strongly-FPR atlas by Theorem~\ref{thm:strong-atlas}, and Theorem~\ref{thm:pure-finite}(ii) applies, as in \cite[Rem.~15.8]{AlgMag}.
\end{proof}

\subsection{Rank and examples}

\begin{definition}[{\cite[Def.~15.9]{AlgMag}}]\label{def:rank}
The \emph{rank} $\operatorname{rk}(N)$ of a monoid $N$ is the least cardinality of a subset generating $N$.
\leansrc{rank}
\end{definition}

\begin{proposition}[{\cite[Fact~15.11]{AlgMag}}]\label{prop:rank}
One has $\operatorname{rk}(0)=0$, $\operatorname{rk}(\NN x)=1$ for $x\neq0$, and $\operatorname{rk}(\ZZ^n)=n+1$ for $n\geq1$.
\leansrc{rankZn}
\end{proposition}

\begin{proof}
\cite[proof of Fact~15.11]{AlgMag}
\end{proof}

\begin{example}[{\cite[Ex.~15.10, examples after Fact~15.11]{AlgMag}}]\label{ex:pure}
Let $S=\Spec(R)$ with $R\neq0$.
\begin{enumerate}
\item For the trivial action on any sheaf $X$, $\mho(a)=\{0\}$.
\item Let $X=\Spec(A)$ be $S$-affine with $A\neq0$, and $M\neq0$. If $A$ has a homogeneous unit in every degree, then $\mho(a)=\{0,M\}$. This applies to $\D(G)_S$ acting on itself by multiplication, $G\neq0$ an abelian group, and to an action of $\D(\ZZ/6\ZZ)_S$ on $\Spec(A)$ with homogeneous units of degrees $2$ and $3$: then $\mho(a)=\{0,\ZZ/6\ZZ\}$.
\item $\A(\NN)_S$ acting on itself: $\mho(a)=\{0,\NN\}$. $\D(\ZZ)_S$ acting on $\A(\NN)_S$ by multiplication: $\mho(a)=\{0,\NN\}$; by $t\cdot x=t^2x$: $\mho(a)=\{0,2\NN\}$.
\item $\D(\ZZ)_S$ acting on $\mathbb{P}^1_S$ by $t\cdot[x:y]=[tx:y]$: $\mho(a)=\{0,\NN,-\NN,\ZZ\}$.
\item $\D(\ZZ^2)_S$ acting on $\A(\NN^2)_S$ by multiplication: $\mho(a)=\{0,\NN\times0,0\times\NN,\NN^2\}$.
\item $\D(\ZZ)_S$ acting on $\A(\mon{2,3})_S$ by multiplication: $\mho(a)=\{0,2\NN,3\NN,\mon{2,3}\}$.
\item $\D(\ZZ/6\ZZ)_S$ acting on $\A(\NN)^3_S=\Spec(R[x_1,x_2,x_3])$ with weights $1,2,3$, that is $\deg x_i=i$: $\mho(a)=\{0,2\ZZ/6\ZZ,3\ZZ/6\ZZ,\ZZ/6\ZZ\}$.
\leansrc{pureZSix}
\end{enumerate}
\leansrc{pureExamples}
\end{example}

\begin{proof}
(i) Every $[X^N]$ is $X$, so a pure magnet is contained in $0$. (ii) If $K\neq M$, choose $m\notin K$ and a unit $u\in A_m$; then $J_K=A$ and $X^K=\emptyset$, while $X^M=X$ and $A\neq0$. So the attractors are $\emptyset$ and $X$, and $E(\emptyset)=0$, $E(X)=M$. In $\ZZ/6\ZZ$, every element is $2a+3b$. (iii), (v), (vi) are cases of Proposition~\ref{prop:self-pure}(i) below with $E=\{1\}$, $\{1\}$, $\{2\}$, $\{(1,0),(0,1)\}$, $\{2,3\}$. (iv) \cite[Example (iv) after Fact~15.11]{AlgMag}; it is also the case $\iota=\{0,1\}$ of Example~\ref{ex:Pn}. (vii) For a magnet $L$, $J_L$ contains $x_i$ if $i\notin L$, and $J_L$ is killed by evaluation at the point with coordinate $1$ at $i\in L$ and $0$ at $i\notin L$, as a monomial of degree $m\notin L$ involves some $x_i$ with $i\notin L$; since $R\neq0$, $J_L$ determines $L\cap\{1,2,3\}$. The submonoids of $\ZZ/6\ZZ$ are $0$, $2\ZZ/6\ZZ$, $3\ZZ/6\ZZ$ and $\ZZ/6\ZZ$, which have distinct intersections with $\{1,2,3\}$; so their attractors are distinct, and all four are pure.
\end{proof}

The self-actions of \cite{SelfActions} are the actions of $\A(M)_S$ on $\A(L)_S$ through the inclusion of a submonoid $L\subseteq M$. Let $q\colon M\to M/\units M$ be the quotient map of Section~\ref{sec:self-actions}. A subset $E$ of a monoid is \emph{minimal} if no proper subset of $E$ generates the submonoid $\mon{E}$.

\begin{proposition}[{\cite[Prop.~3.1, Thm.~4.6]{SelfActions}}]\label{prop:self-pure}
Assume that $S$ is nonempty and that $M$ is cancellative.
\begin{enumerate}
\item Let $E\subseteq M$ be minimal with $\units{\mon{E}}=0$, and let $\A(M)_S$ act on $\A(\mon{E})_S$. Then $\mho(a)=\{\mon{B}\mid B\subseteq E\}$, and $B\mapsto\mon{B}$ is injective on subsets of $E$.
\item Let $\A(M)_S$ act on itself, and let $E\subseteq M/\units M$ be a minimal generating set. If $\units M=0$, then $B\mapsto q^{-1}(\mon{B})$ is a bijection from the set of subsets of $E$ onto $\mho(a)$, and $\#\mho(a)=2^{\#E}$. If $\units M\neq0$, then $\mho(a)=\{0\}\sqcup\{q^{-1}(\mon{B})\mid B\subseteq E\}$ and $\#\mho(a)=2^{\#E}+1$.
\end{enumerate}
\leansrc{selfPure}
\end{proposition}

\begin{proof}
\cite[proofs of Prop.~3.1 and Thm.~4.6]{SelfActions}
\end{proof}

\section{Join actions and join-dense families}\label{sec:rooted}

In this section $M$ is a monoid, $X$ is a sheaf on $(\Sch/S)_{\fppf}$ with an action $a$ of $\A(M)_S$, and $\mho(a)$ is the set of pure magnets of $a$ (Definition~\ref{def:pure}). Let $\operatorname{Sub}(M)$ be the set of submonoids of $M$. For $E\subseteq\operatorname{Sub}(M)$, let $[E\rangle=\sum_{I\in E}I$ be the submonoid generated by the members of $E$, with $[\emptyset\rangle=0$.

\subsection{Join-dense families of pure magnets}

\begin{definition}\label{def:rooted}
Let $\mathcal F\subseteq\operatorname{Sub}(M)$. The \emph{trace} of a magnet $N$ is $\operatorname{tr}_\mathcal F(N)=\{I\in\mathcal F\mid I\subseteq N\}$. A subset $E\subseteq\mathcal F$ is \emph{stable} if $\operatorname{tr}_\mathcal F([E\rangle)=E$. The family $\mathcal F$ is \emph{join-dense} for $a$ if $\operatorname{tr}_\mathcal F$ is a bijection from $\mho(a)$ onto the set of stable subsets of $\mathcal F$, with inverse $E\mapsto[E\rangle$; that is: for $N\in\mho(a)$, $\operatorname{tr}_\mathcal F(N)$ is stable and $[\operatorname{tr}_\mathcal F(N)\rangle=N$, and for every stable $E\subseteq\mathcal F$, $[E\rangle\in\mho(a)$. We then also say that the action $a$ is \emph{join-dense by} $\mathcal F$.
\leansrc{rooted}
\end{definition}

\begin{proposition}\label{prop:rooted}
Let $\mathcal F\subseteq\operatorname{Sub}(M)$.
\begin{enumerate}
\item The family $\mathcal F$ is join-dense for $a$ if and only if $\mho(a)=\{[E\rangle\mid E\subseteq\mathcal F\}$.
\item If $\mathcal F$ is join-dense for $a$, then $\mathcal F\subseteq\mho(a)$, and the sum of every family of pure magnets is pure.
\item If $N=[\operatorname{tr}_\mathcal F(N)\rangle$ for every $N\in\mho(a)$ and $[\operatorname{tr}_\mathcal F(N)\rangle\in\mho(a)$ for every magnet $N$, then $\mathcal F$ is join-dense for $a$.
\end{enumerate}
\leansrc{rootedBasic}
\end{proposition}

\begin{proof}
For $E\subseteq\mathcal F$ one has $E\subseteq\operatorname{tr}_\mathcal F([E\rangle)$, hence $[\operatorname{tr}_\mathcal F([E\rangle)\rangle=[E\rangle$; and $[\operatorname{tr}_\mathcal F(N)\rangle\subseteq N$ for every magnet $N$. (i) Assume that $\mathcal F$ is join-dense for $a$. A pure magnet is $[\operatorname{tr}_\mathcal F(N)\rangle$. Conversely, for $E\subseteq\mathcal F$ the set $E'=\operatorname{tr}_\mathcal F([E\rangle)$ satisfies $[E'\rangle=[E\rangle$, so $\operatorname{tr}_\mathcal F([E'\rangle)=E'$: $E'$ is stable and $[E\rangle=[E'\rangle$ is pure. Assume now $\mho(a)=\{[E\rangle\mid E\subseteq\mathcal F\}$. For $N=[E\rangle$ pure, $N=[E\rangle\subseteq[\operatorname{tr}_\mathcal F(N)\rangle\subseteq N$, so $[\operatorname{tr}_\mathcal F(N)\rangle=N$ and $\operatorname{tr}_\mathcal F(N)$ is stable; a stable $E$ gives the pure magnet $[E\rangle$. (ii) For $I\in\mathcal F$, $I=[\{I\}\rangle$ is pure by (i). A family of pure magnets $[E_j\rangle$ has sum $[\bigcup_jE_j\rangle$, pure by (i). (iii) By hypothesis, $\mho(a)$ is the set of the $[\operatorname{tr}_\mathcal F(N)\rangle$; it contains $[E\rangle=[\operatorname{tr}_\mathcal F([E\rangle)\rangle$ for every $E\subseteq\mathcal F$ and is contained in $\{[E\rangle\mid E\subseteq\mathcal F\}$. Apply (i).
\end{proof}

A subset $\Omega\subseteq M$ is \emph{join-dense} for $a$ if the family $\{\NN\alpha\mid\alpha\in\Omega\}$ is. We call $a$ a \emph{join action} if the sum of every family of pure magnets is pure. Then $\mho(a)$, ordered by inclusion, is a complete join-semilattice, the join being the sum; equivalently, $(\mho(a),+)$ is a commutative idempotent monoid with neutral element $0$. By Proposition~\ref{prop:rooted}, $\mathcal F$ is join-dense for $a$ exactly when $a$ is a join action and $\mathcal F\subseteq\mho(a)$ is join-dense in the lattice-theoretic sense: every pure magnet is the join of the members of $\mathcal F$ below it (see \cite[p.~354]{ErneWilke83}, \cite[pp.~256, 269]{ErneReichman87}, \cite[p.~514]{Wang23}), that is, $\mathcal F$ generates the semilattice $\mho(a)$. A subset $E\subseteq\Omega$ is \emph{additively stable} if $\mon{E}\cap\Omega=E$.

\begin{proposition}\label{prop:join-action}
Let $M$ be a finitely generated abelian group and $X$ a separated algebraic space over $S$, locally of finite presentation, with an action $a$ of $\D(M)_S$ admitting an $S$-affine strongly fixed-point-reflecting atlas; by Theorem~\ref{thm:strong-atlas}, such an atlas exists if $S$ is separated and $|X|$ is quasi-compact. Then $a$ is a join action.
\leansrc{joinAction}
\end{proposition}

\begin{proof}
Let $(N_i)$ be a family of pure magnets, $N=\sum_iN_i$, and $L$ a magnet with $X^L=X^N$. Fix $i$. Since $N_i\subseteq N$, $X^{N_i}\subseteq X^N=X^L$. Let $\psi\colon\A(N_i)_T\to X_T$ be a point of $X^{N_i}$; for every point $b$ of $\A(N_i)$, the point $\psi(b)$ of $X$ lies in $X^{N_i}$, its orbit map extending to $c\mapsto\psi(cb)$, hence in $X^L$. As $X^L\to X$ is a monomorphism, $\psi$ factors through $X^L$ and is a point of $(X^L)^{N_i}=X^{L\cap N_i}$ (Proposition~\ref{prop:intersection}). Thus $X^{L\cap N_i}=X^{N_i}$, and since $N_i$ is pure, $N_i\subseteq L\cap N_i\subseteq L$. Hence $N\subseteq L$, and $N$ is pure.
\end{proof}

\begin{lemma}\label{lem:rooted-elements}
Let $\Omega\subseteq M$. The subset $\Omega$ is join-dense for $a$ if and only if the following conditions hold: for $N\in\mho(a)$, $N\cap\Omega$ is additively stable and $\mon{N\cap\Omega}=N$; and for every additively stable $E\subseteq\Omega$, $\mon E\in\mho(a)$. Then $N\mapsto N\cap\Omega$ is a bijection from $\mho(a)$ onto the set of additively stable subsets of $\Omega$, with inverse $E\mapsto\mon E$.
\leansrc{rootedElements}
\end{lemma}

\begin{proof}
For $\mathcal F_1=\{\NN\alpha\mid\alpha\in\Omega\}$ and a magnet $N$, $\operatorname{tr}_{\mathcal F_1}(N)=\{\NN\alpha\mid\alpha\in N\cap\Omega\}$ and $[\{\NN\alpha\mid\alpha\in E\}\rangle=\mon E$ for $E\subseteq\Omega$. An additively stable $E$ is saturated: if $\beta\in\Omega$ and $\NN\beta=\NN\alpha$ with $\alpha\in E$, then $\beta\in\mon E\cap\Omega=E$. Hence $E\mapsto\{\NN\alpha\mid\alpha\in E\}$ is a bijection from the set of additively stable subsets of $\Omega$ onto the set of stable subsets of $\mathcal F_1$, with inverse $E_1\mapsto\{\alpha\in\Omega\mid\NN\alpha\in E_1\}$. Under these identifications the conditions of Definition~\ref{def:rooted} for $\mathcal F_1$ become the stated ones.
\end{proof}

\subsection{Join-irreducible pure magnets}

\begin{theorem}\label{thm:minimal-root}
Let $\underline{\mho}(a)$ be the set of \emph{join-irreducible} pure magnets, the $N\in\mho(a)$ such that $N\neq\sum\{L\in\mho(a)\mid L\subsetneq N\}$; for a join action, these are the join-irreducible elements of the semilattice $\mho(a)$.
\begin{enumerate}
\item Every family join-dense for $a$ contains $\underline{\mho}(a)$.
\item Assume that some family is join-dense for $a$ (so $a$ is a join action) and that every nonempty subset of $\mho(a)$ has a minimal element for inclusion; this holds if $\mho(a)$ is finite. Then $\underline{\mho}(a)$ is join-dense for $a$; it is the smallest family join-dense for $a$.
\end{enumerate}
\leansrc{minimalRoot}
\end{theorem}

\begin{proof}
(i) Let $\mathcal F$ be join-dense for $a$ and $N\in\underline{\mho}(a)$. The members of $\operatorname{tr}_\mathcal F(N)$ are pure (Proposition~\ref{prop:rooted}(ii)) and $N=[\operatorname{tr}_\mathcal F(N)\rangle$. If none of them were equal to $N$, then $N\subseteq\sum\{L\in\mho(a)\mid L\subsetneq N\}\subseteq N$, which is excluded. So $N\in\mathcal F$. (ii) By well-founded induction on $N\in\mho(a)$, every pure magnet $N$ is the sum of the members of $\underline{\mho}(a)$ contained in $N$: if $N\notin\underline{\mho}(a)$, $N$ is the sum of the pure $L\subsetneq N$, each of which is such a sum by induction. Hence every pure magnet is a sum of members of $\underline{\mho}(a)$; conversely a sum of members of $\underline{\mho}(a)\subseteq\mho(a)$ is pure by Proposition~\ref{prop:rooted}(ii), some family being join-dense for $a$. So $\mho(a)=\{[E\rangle\mid E\subseteq\underline{\mho}(a)\}$ and $\underline{\mho}(a)$ is join-dense by Proposition~\ref{prop:rooted}(i). The last assertion follows from (i).
\end{proof}

An action $a$ is \emph{join-irreducible} if it is join-dense by the family $\underline{\mho}(a)$ of its join-irreducible pure magnets; by Theorem~\ref{thm:minimal-root}(i), $\underline{\mho}(a)$ is then the smallest family join-dense for $a$.

\begin{corollary}\label{cor:join-dense-minimal}
Let $S$ be separated and quasi-compact, $M$ a finitely generated abelian group and $X$ a separated algebraic space over $S$ of finite presentation, with an action $a$ of $\D(M)_S$. Then $a$ is join-irreducible: it is join-dense by a unique minimal family of magnets, which is contained in every family join-dense for $a$, the finite family $\underline{\mho}(a)$ of join-irreducible pure magnets.
\leansrc{joinDenseMinimal}
\end{corollary}

\begin{proof}
By Corollary~\ref{cor:conj157}, $\mho(a)$ is finite. By Proposition~\ref{prop:join-action}, with the $S$-affine strongly fixed-point-reflecting atlas of Theorem~\ref{thm:strong-atlas}, $a$ is a join action, so $a$ is join-dense by $\mho(a)$. Theorem~\ref{thm:minimal-root} concludes.
\end{proof}

\subsection{Projective spaces}

Let $R\neq0$ be a ring, $S=\Spec(R)$, $\iota$ a finite set, $M$ a finitely generated abelian group and $w\colon\iota\to M$ a map. The group $\D(M)_S$ acts on $\mathbb{P}(\iota)_S$, the projective space of rank one direct summands of $\cO^\iota$, by $h\cdot[v]=[(h(w_i)v_i)_{i\in\iota}]$. Put $\Omega_w=\{w_i-w_j\mid i,j\in\iota,\ i\neq j\}\subseteq M$.

\begin{theorem}\label{thm:Pn}
The subset $\Omega_w$ is join-dense for the action of $\D(M)_S$ on $\mathbb{P}(\iota)_S$: a magnet $N$ is pure if and only if $N=\mon{N\cap\Omega_w}$, and $N\mapsto N\cap\Omega_w$ is a bijection from the set of pure magnets onto the set of additively stable subsets of $\Omega_w$, with inverse $E\mapsto\mon E$.
\leansrc{rootedPn}
\end{theorem}

\begin{proof}
For $s\in\iota$, let $U_s=\mathbb{A}^{\iota\setminus\{s\}}_S$ be the standard chart $v_s\neq0$, with coordinates $x_i=v_i/v_s$. The coordinate $x_i$ has weight $w_i-w_s$, so $U_s$ is a $\D(M)_S$-stable affine open subscheme and, by Theorem~\ref{thm:affine}, $U_s^N$ is the coordinate subspace $x_i=0$ for $w_i-w_s\notin N$; it depends only on $N\cap\Omega_w$. Let $\varphi\in X^N(T)$, $X=\mathbb{P}(\iota)_S$, and $x=p_{N,\units N}(\varphi)$. Zariski locally on $T$, $x(1)$ lies in $U_s$ for some $s$; then $x\in U_s^{\units N}$. Indeed, the orbit morphism $o\colon\D(M)_T\to\D(\units N)_T$, $h\mapsto h\cdot1$, is faithfully flat, $\ZZ[M]$ being free over $\ZZ[\units N]$, hence surjective; and $x\circ o$, $h\mapsto h\cdot x(1)$, factors through the $\D(M)_S$-stable open subscheme $U_s$. So the open subscheme $x^{-1}(U_s)$ of $\D(\units N)_T$ contains the image of $o$, which is everything, and $x$ factors through $U_s$. Then $\varphi\in U_s^N$ by Proposition~\ref{prop:etale-face}(ii). Hence, Zariski locally, $[X^N]$ is the union of the $U_s^N$, and $[X^N]=[X^{\mon{N\cap\Omega_w}}]$. Conversely let $N,L$ be magnets with $w_i-w_s\in N\setminus L$. Over the dual numbers $R[\varepsilon]$, the point $x_i=\varepsilon$, $x_j=0$ $(j\neq i)$ of $U_s$ lies in $U_s^N$; it lies in no other chart $U_t$, since its reduction is $[e_s]$, so it does not lie in $[X^L]$. Therefore two magnets with the same attractor have the same intersection with $\Omega_w$, and $E([X^N])=\mon{N\cap\Omega_w}$ for every magnet $N$. A magnet $N$ is thus pure if and only if $N=\mon{N\cap\Omega_w}$, and for $E$ additively stable, $\mon E$ is pure since $\mon{\mon E\cap\Omega_w}=\mon E$. The assertions follow from Lemma~\ref{lem:rooted-elements}.
\end{proof}

\begin{example}\label{ex:Pn}
(i) For $\iota=\{0,1\}$, $M=\ZZ$ and $w=(1,0)$, $\Omega_w=\{1,-1\}$ and the pure magnets of the action on $\mathbb{P}^1_S$ are $0,\NN,-\NN,\ZZ$. (ii) For $\iota=\{x,y,z\}$, $M=\ZZ^2$ and $w=((1,0),(0,1),(0,0))$, the action is $(\lambda,\beta)\cdot[x:y:z]=[\lambda x:\beta y:z]$ on $\mathbb{P}^2_S$, and $\Omega_w=\{\pm(1,0),\pm(0,1),\pm(1,-1)\}$. The subset $\Omega_w$ is join-dense and the action has $29$ pure magnets, the monoids $\mon E$ for the $29$ additively stable subsets $E\subseteq\Omega_w$.
\leansrc{rootedP2}
\end{example}

\section{Reductive groups}\label{sec:reductive}

The Lean declarations indicated in this section are public fragments of the code that give a partial formalization of its statements; this section is not presented as formalized.

Let $k$ be a field, $\bar k$ an algebraic closure of $k$, $G$ an affine algebraic group over $k$ and $M$ a finitely generated abelian group with a homomorphism $\D(M)_k\to G$; then $\D(M)_k$ acts on $G$ by conjugation, by group automorphisms. For a magnet $N$, $G^N$ is a closed subgroup scheme of $G$ (Theorem~\ref{thm:affine}). The \emph{weights} of $G$ are the $m\in M$ with $\Lie(G)_m\neq0$, and the \emph{roots} are the nonzero weights. For a morphism $f\colon M\to\ZZ$, let $\lambda_f\colon\mathbb{G}_m\to\D(M)_k$ be the dual cocharacter, $P(\lambda_f)$, $Z(\lambda_f)$ the attached parabolic and Levi subgroups of $G$ \cite[\S2.1]{CGP}; then $P(\lambda_f)=G^{f^{-1}(\NN)}$ and $Z(\lambda_f)=G^{f^{-1}(0)}$ (Proposition~\ref{prop:change-monoid}(i)).
\leansrc{redSetting}

\subsection{Parabolic subgroups}
\begin{theorem}[{generalizing \cite[Prop.~16.10(ii), (iii)]{AlgMag}}]\label{thm:parabolic}
Assume that $G$ is smooth and that $G_{\bar k}$ is connected. Let $f\colon M\to\ZZ$ be a morphism and $N\subseteq f^{-1}(\NN)$ a magnet such that $N\cap\Phi=f^{-1}(\NN)\cap\Phi$, where $\Phi$ is the set of roots of $G$. Then $G^N=P(\lambda_f)$.
\leansrc{parabolic}
\end{theorem}

\begin{proof}
Since $N\subseteq f^{-1}(\NN)$, $G^N\subseteq G^{f^{-1}(\NN)}=P(\lambda_f)$ (Proposition~\ref{prop:magnet-maps}). Attractors and the ideals of $P(\lambda_f)$ commute with the base change $k\to\bar k$ (Proposition~\ref{prop:attractor-limits}), and an equality of ideals of $\cO(G)$ descends along the faithfully flat $k\to\bar k$; so we may assume $k=\bar k$. The group $G^N$ is smooth by \cite[Cor.~11.2]{AlgMag} applied to $G\to Y=\Spec(k)$, for which $Y^N=Y$, and $\Lie(G^N)=\Lie(G)^N=\bigoplus_{m\in N}\Lie(G)_m$ by \cite[Prop.~16.4]{AlgMag} and Proposition~\ref{prop:vector-bundle}; by the hypothesis on $N$ and since $0\in N$, this is $\bigoplus_{f(m)\geq0}\Lie(G)_m=\Lie(P(\lambda_f))$ \cite[Prop.~2.1.8(1)]{CGP}. The multiplication map $U(\lambda_f^{-1})\times P(\lambda_f)\to G$ is an open immersion \cite[Prop.~2.1.8(3)]{CGP}; so $U(\lambda_f^{-1})\times P(\lambda_f)$ is open in the irreducible $G$ (smooth and connected over $\bar k$), hence irreducible, and $P(\lambda_f)$, its image under the second projection, is irreducible, in particular connected. The group $P(\lambda_f)=G^{f^{-1}(\NN)}$ is smooth, as an attractor of the smooth $G$ \cite[Cor.~11.2]{AlgMag} (with $Y=\Spec(k)$). A smooth closed subgroup of a smooth connected group over $\bar k$ with the same Lie algebra is the whole group, so $G^N=P(\lambda_f)$.
\end{proof}

\subsection{Split reductive groups}

\begin{theorem}\label{thm:red-equivalence}
Let $k$ be an algebraically closed field, $G$ a smooth connected affine group over $k$ and $T\subseteq G$ a maximal torus, with character lattice $X$ and set of roots $\Phi\subseteq X$, the nonzero weights of $T$ on $\Lie(G)$. Then $G$ is reductive if and only if the following three conditions hold:
\begin{enumerate}
\item $Z_G(T)=T$;
\item for every $\alpha\in\Phi$, $\dim\Lie(G)_\alpha=1$ and $\Phi\cap\mathbb{Q}\alpha=\{\pm\alpha\}$;
\item for every $\alpha\in\Phi$, there is $n_\alpha\in N_G(T)(k)$ acting on $X$ by a reflection $s_\alpha$ with $s_\alpha(\alpha)=-\alpha$.
\end{enumerate}
\leansrc{redEquiv}
\end{theorem}

\begin{proof}
This is classical; see \cite{SGA3III}.
\end{proof}

From now on $G$ is a split connected reductive group over $k$, $T=\D(M)_k\subseteq G$ is a split maximal torus, $\Phi\subseteq M$ is the set of roots, $W=N_G(T)/T$ is the Weyl group, acting faithfully on $M$, and $n_w\in N_G(T)(k)$ represents $w\in W$ \cite[Exp.~XXII]{SGA3III}. Let $a$ be the conjugation action of $T$ on $G$. A point $g\in G(A)$, $A$ a $k$-algebra, \emph{normalizes $T$ with automorphism $\sigma$} of $M$ if $gtg^{-1}=\sigma(t)$ for the universal point $t$ of $T$; the automorphisms obtained in this way form a subgroup of $\operatorname{Aut}(M)$, and $n_w$ normalizes $T$ with automorphism $w$.
\leansrc{weylGroup}

\begin{definition}\label{def:root-group}
Let $\alpha\in\Phi$. A \emph{root group} for $\alpha$ is a $T$-stable closed subgroup $U\subseteq G$ with $\Lie(U)=\Lie(G)_\alpha$ together with an isomorphism $\mathbb{G}_a\cong U$ that is $T$-equivariant for the action of $T$ on $\mathbb{G}_a$ through $\alpha$.
\leansrc{rootGroupDef}
\end{definition}

\begin{proposition}[{\cite[Prop.~16.12, Cor.~16.13, Prop.~16.16]{AlgMag}}]\label{prop:root-groups}
Let $\alpha\in\Phi$ and let $U_\alpha=G^{\mon{\alpha}}_e$ be the attractor of the magnet $\mon{\alpha}=\NN\alpha$ with prescribed limit the unit section (Definition~\ref{def:prescribed}).
\begin{enumerate}
\item $U_\alpha$ is a root group for $\alpha$: it is a $T$-stable closed subgroup of $G$ with $\Lie(U_\alpha)=\Lie(G)_\alpha$, and $\cO(U_\alpha)=k[y]$ for the class $y$ of an element of $\cO(G)$ of weight $\alpha$, with $\Delta(y)=y\otimes1+1\otimes y$.
\item For every $k$-algebra $A$, the limit map $G^{\mon{\alpha}}(A)\to G^0(A)$ is a split surjective homomorphism with kernel $U_\alpha(A)$, so $G^{\mon{\alpha}}(A)=U_\alpha(A)\rtimes G^0(A)$.
\item If $k$ is algebraically closed, then $G^0=T$, so $G^{\mon{\alpha}}(A)=U_\alpha(A)\rtimes T(A)$ for every $k$-algebra $A$.
\item $U_\alpha$ is the unique root group for $\alpha$: every root group for $\alpha$ equals $U_\alpha$.
\end{enumerate}
\leansrc{rootGroups}
\end{proposition}

\begin{proof}
(i) By Theorem~\ref{thm:affine}, the ideal of $U_\alpha$ is generated by the homogeneous elements of $\cO(G)$ of weight $m\notin\NN\alpha$ and by the $x-\varepsilon(x)$, $x$ of weight $0$; it is homogeneous, so $U_\alpha$ is $T$-stable, and $\cO(U_\alpha)=\bigoplus_{n\geq0}\cO(U_\alpha)_{n\alpha}$ with $\cO(U_\alpha)_0=k$, $\alpha$ having infinite order as $M$ is free. A tangent vector at $e$ lies in $\Lie(U_\alpha)$ if and only if it vanishes on these generators, that is, on the weight spaces of weight $\neq\alpha$, the weights $n\alpha$, $n\geq2$, not being roots \cite[Exp.~XIX, 3.10]{SGA3III}. The root space $\Lie(G)_\alpha$ has dimension $1$ \cite[Exp.~XIX, Prop.~1.12]{SGA3III} (over $\bar k$, and dimensions of weight spaces are invariant under $k\to\bar k$), so the cotangent space of $U_\alpha$ at $e$ is a line, spanned by the class of an element $y$ of weight $\alpha$. By Theorem~\ref{thm:BB} applied to the face $0$ of $\NN\alpha$, $G^{\mon{\alpha}}\to G^0$ is the spectrum of a symmetric algebra, and its fibre $U_\alpha$ at $e$ is $\Spec\Sym_k$ of its cotangent space at $e$; hence $\cO(U_\alpha)=k[y]$. The element $\Delta(y)$ has weight $\alpha$ in $\cO(U_\alpha)\otimes\cO(U_\alpha)$, whose weight-$\alpha$ part is $ky\otimes1\oplus k\,1\otimes y$; the counit axioms give $\Delta(y)=y\otimes1+1\otimes y$. So $y$ defines a $T$-equivariant isomorphism $U_\alpha\cong\mathbb{G}_a$, $T$ acting on $\mathbb{G}_a$ through $\alpha$. (ii) The limit map is a homomorphism with section the inclusion $G^0\subseteq G^{\mon{\alpha}}$, and its kernel is $G^{\mon{\alpha}}_e$ by definition. (iii) $G^0=G^T=Z_G(T)=T$ \cite[Exp.~XIX, Lemme~1.6.2(ii)]{SGA3III}.
 (iv) Let $U\subseteq G$ be a root group for $\alpha$, with ideal $I$. As $U$ is $T$-stable, $I$ is homogeneous, and the $T$-equivariant isomorphism $\mathbb{G}_a\cong U$ gives $\cO(U)=k[z]$ with $z$ the class of an element of weight $\alpha$; so $\cO(U)_{n\alpha}=kz^n$ for $n\geq0$ and the other weight spaces of $\cO(U)$ vanish. An element of $\cO(G)$ of weight $m\notin\NN\alpha$ therefore lies in $I$, and an element $x$ of weight $0$ restricts to a constant, equal to $\varepsilon(x)$ by evaluation at $e$, so $x-\varepsilon(x)\in I$. By Theorem~\ref{thm:affine} these elements generate the ideal of $U_\alpha$, so $U\subseteq U_\alpha$. The restriction map $k[y]=\cO(U_\alpha)\to\cO(U)=k[z]$ is surjective and sends $y$, of weight $\alpha$, to $cz$ with $c\in k$; surjectivity gives $c\neq0$, so the map is bijective and $U=U_\alpha$.
\end{proof}

For a set $R$ of roots let $K_R=(G^{\mon R})^\circ$, the identity component of the attractor of the magnet generated by $R$.

For a magnet $N$ put $N^*=\units N$, $R(N)=N\cap\Phi$,
\[W_N=\{w\in W\mid(1-w)M\subseteq N\},\qquad\Lambda_N=\sum_{w\in W_N}(1-w)M.\]

\begin{definition}\label{def:weyl-root}
For $w\in W$, put $(1-w)M=\{m-w(m)\mid m\in M\}$, a subgroup of $M$, and let $T^w\subseteq T$ be the fixed-point subgroup of $w$, that is $T^w=\D(M/(1-w)M)_k$. The \emph{Weyl magnet} of $w$ is the magnet $(T^w)^\perp=(1-w)M$, the subgroup of characters of $T$ trivial on $T^w$.
\leansrc{weylRoot}
\end{definition}

\begin{proposition}\label{prop:weyl-root}
Let $w\in W$ and $Z=(1-w)M$.
\begin{enumerate}
\item Let $A$ be a nonzero $k$-algebra and $g\in G(A)$ a point normalizing $T$ with automorphism $w$. For every magnet $N$, $g\in G^N(A)$ if and only if $Z\subseteq N$. In particular $n_w\in G^N$ if and only if $w\in W_N$.
\item $G^Z=Z_G(T^w)$, the centralizer of $T^w$.
\end{enumerate}
\leansrc{weylRootAttractor}
\end{proposition}

\begin{proof}
(i) For $t$ the universal point of $T$, $tgt^{-1}=g\,w^{-1}(t)t^{-1}$. So the orbit morphism $T_A\to G_A$ of $g$ is the composite of the homomorphism $T_A\to T_A$ dual to $1-w^{-1}$ with $T_A\to G_A$, $s\mapsto gs$. The comorphism of the second map is surjective onto $A[M]$, $\cO(G)\to\cO(T)$ being surjective, so the comorphism of the orbit morphism has image $A[(1-w^{-1})M]=A[Z]$. By Theorem~\ref{thm:affine}, $g\in G^N(A)$ if and only if this image lies in $A[N]$, that is, as $A\neq0$, if and only if $Z\subseteq N$. (ii) By Proposition~\ref{prop:subgroup-fixed}, $G^Z=G^{\D(M/Z)_k}=G^{T^w}$, the fixed points of the conjugation action of $T^w$, which is $Z_G(T^w)$.
\end{proof}

\begin{lemma}\label{lem:weyl-part}
Let $N$ be a magnet.
\begin{enumerate}
\item $W_N=\{w\in W\mid(1-w)M\subseteq N^*\}$; it is a subgroup of $W$ and stabilizes $N$, $N^*$ and $R(N)$.
\item $R(N)$ is additively stable, hence closed.
\item $W(R(N)\cap-R(N))$ is a normal subgroup of $W_N$, and $(1-w)M\subseteq\ZZ(R(N)\cap-R(N))$ for $w$ in it.
\end{enumerate}
\leansrc{weylPart}
\end{lemma}

\begin{proof}
(i) $(1-w)M$ is a group, so it lies in $N$ if and only if it lies in $N^*$. For $v,w\in W_N$, $(1-vw)m=(1-v)m+v(1-w)m$ and $v(x)=x-(1-v)x\in N^*$ for $x\in N^*$, so $vw\in W_N$; $W$ is finite, so $W_N$ is a subgroup. For $n\in N$, $w(n)=n-(1-w)n\in N$, and $w$ preserves $\Phi$. (ii) $\mon{R(N)}\subseteq N$. (iii) For $\beta\in R_s=R(N)\cap-R(N)$, $(1-s_\beta)M=\langle M,\beta^\vee\rangle\beta\subseteq\ZZ\beta\subseteq N^*$, so $s_\beta\in W_N$; $ws_\beta w^{-1}=s_{w\beta}$ and $W_N$ stabilizes $R_s$. The last assertion follows from $(1-s_\beta)M\subseteq\ZZ\beta$ and $(1-uv)=(1-u)+u(1-v)$, $W(R_s)$ stabilizing $\ZZ R_s$.
\end{proof}

\begin{proposition}\label{prop:red-identity}
Let $N\subseteq N'$ be magnets with $N\cap\Phi=N'\cap\Phi$. Then $(G^N_{\bar k})^\circ=(G^{N'}_{\bar k})^\circ$. In particular $(G^N_{\bar k})^\circ=(K_{R(N)})_{\bar k}$ for every magnet $N$.
\leansrc{redIdentity}
\end{proposition}

\begin{proof}
The roots of $G_{\bar k}$ are roots of $G$ and attractors commute with $k\to\bar k$ (Proposition~\ref{prop:attractor-limits}); so we may assume $k=\bar k$. The groups $G^N\subseteq G^{N'}$ are smooth by \cite[Cor.~11.2]{AlgMag} applied to $G\to Y=\Spec(k)$, for which $Y^N=Y$, and $\Lie(G^N)=\Lie(G)^N=\Lie(G)^{N'}=\Lie(G^{N'})$ by \cite[Prop.~16.4]{AlgMag} and Proposition~\ref{prop:vector-bundle}, as $N$ and $N'$ contain the same roots. Smooth connected closed subgroups $(G^N)^\circ\subseteq(G^{N'})^\circ$ with the same Lie algebra are equal. The last assertion is the case $\mon{R(N)}\subseteq N$.
\end{proof}

\begin{proposition}\label{prop:red-reflects}
Let $N$, $L$ be magnets.
\begin{enumerate}
\item If $\Lie(G^N)=\Lie(G^L)$, then $R(N)\subseteq L$. If $G^N=G^L$, then $R(N)\subseteq L$ and $W_N\subseteq W_L$.
\item Assume $k=\bar k$. For $w\in W_N$, $K_{R(N)}n_w\subseteq G^N$; hence $\bigcup_{w\in W_N}K_{R(N)}n_w\subseteq G^N$.
\end{enumerate}
\leansrc{redReflects}
\end{proposition}

\begin{proof}
(i) Let $\alpha\in R(N)$ and $v$ a nonzero tangent vector of weight $\alpha$; by \cite[Prop.~16.4]{AlgMag} and Proposition~\ref{prop:vector-bundle}, $v\in\Lie(G^N)=\Lie(G^L)=\bigoplus_{m\in L}\Lie(G)_m$, so $\alpha\in L$. If $G^N=G^L$ and $w\in W_N$, then $n_w\in G^N=G^L$, so $w\in W_L$ by Proposition~\ref{prop:weyl-root}(i). (ii) By Proposition~\ref{prop:red-identity}, $K_{R(N)}=(G^N)^\circ\subseteq G^N$, and $n_w\in G^N$ by Proposition~\ref{prop:weyl-root}(i).
\end{proof}

\begin{example}\label{ex:red}
(i) For $G=\mathrm{PGL}_2$, $M=\ZZ\alpha$ and $(1-s_\alpha)M=2\ZZ\alpha$; by Proposition~\ref{prop:weyl-root}(ii), $G^{2\ZZ\alpha}=Z_G(\mu_2)=N_G(T)$.
(ii) For $G=\mathrm{PGL}_3$ and $c$ a Coxeter element, $(1-c)M=\{xa+yb\mid x+y\equiv0\bmod3\}$ for the simple roots $a,b$; it contains no root, and by Proposition~\ref{prop:weyl-root}(ii), $G^{(1-c)M}=Z_G(\mu_3)=T\sqcup n_cT\sqcup n_c^2T$.

In both cases $N=(1-w)M$, $w=s_\alpha$ or $c$, contains no root. So $(G^N_{\bar k})^\circ=(K_\emptyset)_{\bar k}=(G^0_{\bar k})^\circ=T_{\bar k}$ by Propositions~\ref{prop:red-identity} and~\ref{prop:root-groups}(iii). The group $G^N_{\bar k}$ normalizes its identity component, so $G^N_{\bar k}\subseteq N_G(T)_{\bar k}$, and $G^N_{\bar k}\cap N_G(T)_{\bar k}=\bigcup_{v\in W_N}n_vT_{\bar k}$ by Proposition~\ref{prop:weyl-root}(i); both sides are smooth ($G^N$ by \cite[Cor.~11.2]{AlgMag}), so this is an equality of subschemes, and it descends to $k$. Finally $W_N=W$ for $\mathrm{PGL}_2$; for $\mathrm{PGL}_3$, $W_N=\langle c\rangle$, since $(1-c^2)M=(1-c)(1+c)M\subseteq(1-c)M$ while $(1-s_\beta)M=\langle M,\beta^\vee\rangle\beta=\ZZ\beta$ contains the root $\beta$ for each reflection $s_\beta$, as $M=\ZZ a\oplus\ZZ b$ and $\langle a,b^\vee\rangle=\langle b,a^\vee\rangle=-1$, $\langle a,(a+b)^\vee\rangle=1$.
\leansrc{exRed}
\end{example}

\begin{lemma}\label{lem:red-descent}
Let $N$, $L$ be magnets. Then $G^N=G^L$ if and only if $G^N_{\bar k}=G^L_{\bar k}$. The group $G_{\bar k}$ is split connected reductive with split maximal torus $T_{\bar k}=\D(M)_{\bar k}$, set of roots $\Phi$ and Weyl group $W$. Hence $E(G^N)=E(G^N_{\bar k})$, $N$ is pure for $G$ if and only if it is pure for $G_{\bar k}$, and $R(N)$, $W_N$, $\Lambda_N$ are the same for $G$ and $G_{\bar k}$.
\leansrc{redDescent}
\end{lemma}

\begin{proof}
By Proposition~\ref{prop:attractor-limits}, $G^N_{\bar k}$ is defined by $I_N\otimes_k\bar k$, $I_N$ the ideal of $G^N$, and $I_N=(I_N\otimes_k\bar k)\cap\cO(G)$ as $k\to\bar k$ is faithfully flat; this gives the first assertion. The weight decomposition of $\Lie(G_{\bar k})=\Lie(G)\otimes_k\bar k$ is that of $\Lie(G)$ tensored with $\bar k$, so the roots are the same, and the Weyl group is that of the root datum \cite[Exp.~XXII]{SGA3III}, with the same representatives $n_w$. The remaining assertions follow, $E(G^N)$ being the intersection of the magnets $L$ with $G^L=G^N$.
\end{proof}

\begin{theorem}\label{thm:red-attractors}
Assume that $k$ is algebraically closed. Let $N$ be a magnet and $R=R(N)$.
\begin{enumerate}
\item $G^N$ is smooth, $(G^N)^\circ=K_R$, and $\Lie(G^N)=\Lie(T)\oplus\bigoplus_{\alpha\in R}\Lie(G)_\alpha$.
\item $G^N(k)=\bigcup_{w\in W_N}K_R(k)\,n_w$.
\end{enumerate}
In particular $G^N$ depends only on the pair $(R(N),W_N)$.
\leansrc{redAttractors}
\end{theorem}

\begin{proof}
(i) By \cite[Cor.~11.2]{AlgMag} applied to $G\to Y=\Spec(k)$, for which $Y^N=Y$, $G^N$ is smooth, Proposition~\ref{prop:red-identity} gives the identity component, and the Lie algebra is given by \cite[Prop.~16.4]{AlgMag} and Proposition~\ref{prop:vector-bundle}. (ii) The inclusion $\supseteq$ is Proposition~\ref{prop:red-reflects}(ii). Let $g\in G^N(k)$. Then $gK_Rg^{-1}=K_R$, so $gTg^{-1}$ is a maximal torus of $K_R$; by conjugacy of maximal tori in the smooth connected group $K_R$ \cite[Exp.~XII, Th.~6.6(a)]{SGA3}, $gTg^{-1}=hTh^{-1}$ for some $h\in K_R(k)$. Then $h^{-1}g\in N_G(T)(k)\cap G^N(k)=\bigcup_{w\in W_N}n_wT(k)$ by Proposition~\ref{prop:weyl-root}(i). So $g\in K_R(k)\,n_wT(k)=K_R(k)\,n_w$ for some $w\in W_N$, as $T\subseteq K_R$ and $n_wTn_w^{-1}=T$. If $N'\subseteq N$ is a magnet with $R(N')=R(N)$ and $W_{N'}=W_N$, then (i) and (ii) give $G^{N'}(k)=G^N(k)$, hence $G^{N'}=G^N$, both being smooth; for two magnets $N$, $L$ with the same pair, apply this to $N\cap L\subseteq N$ and $N\cap L\subseteq L$.
\end{proof}

\begin{theorem}\label{thm:red-pure}
For every magnet $N$, $E(G^N)=\mon{R(N)}+\Lambda_N$. Hence $N$ is pure if and only if $N=\mon{N\cap\Phi}+\Lambda_N$, and $\mon{R(N)}+\Lambda_N$ is pure.
\leansrc{redPure}
\end{theorem}

\begin{proof}
By Lemma~\ref{lem:red-descent} we may assume $k=\bar k$. Put $L_0=\mon{R(N)}+\Lambda_N$. Then $L_0\subseteq N$, so $R(L_0)\subseteq R(N)$ and $W_{L_0}\subseteq W_N$; and $R(N)\subseteq L_0$, $(1-w)M\subseteq L_0$ for $w\in W_N$. So $(R(L_0),W_{L_0})=(R(N),W_N)$ and $G^{L_0}=G^N$ (Theorem~\ref{thm:red-attractors}). If $G^L=G^N$, Proposition~\ref{prop:red-reflects}(i) gives $R(N)\subseteq L$ and $(1-w)M\subseteq L$ for $w\in W_N$, that is $L_0\subseteq L$. So $L_0=E(G^N)$, which is pure (Theorem~\ref{thm:pure}).
\end{proof}

\subsection{Classification of pure magnets}

For $R\subseteq\Phi$ and a subgroup $W'\subseteq W$ put $L(R,W')=\mon R+\sum_{w\in W'}(1-w)M$.

\begin{theorem}\label{thm:red-classification}
\begin{enumerate}
\item The map $N\mapsto(R(N),W_N)$ is a bijection from the finite set $\mho(a)$ onto the set of pairs $(R,W')$, $R\subseteq\Phi$, $W'\subseteq W$ a subgroup, such that
\[\Phi\cap L(R,W')=R\qquad\text{and}\qquad\{w\in W\mid(1-w)M\subseteq L(R,W')\}=W',\]
with inverse $(R,W')\mapsto L(R,W')$; for such a pair, if $k$ is algebraically closed, $G^{L(R,W')}(k)=\bigcup_{w\in W'}K_R(k)\,n_w$.
\item The map $N\mapsto N\cap\Phi$ from $\mho(a)$ to the set of additively stable subsets of $\Phi$ is surjective: for $R$ additively stable, $E(G^{\mon R})$ is a pure magnet with $E(G^{\mon R})\cap\Phi=R$.
\end{enumerate}
\leansrc{redClassification}
\end{theorem}

\begin{proof}
(i) For $N$ pure, $N=L(R(N),W_N)$ by Theorem~\ref{thm:red-pure}, so the map is injective and its values satisfy the two conditions. Conversely, if $(R,W')$ satisfies them and $L=L(R,W')$, then $R(L)=R$ and $W_L=W'$, and $E(G^L)=\mon R+\Lambda_L=L$ by Theorem~\ref{thm:red-pure}, so $L$ is pure. The attractor is given by Theorem~\ref{thm:red-attractors}(ii), as $R(L)=R$ and $W_L=W'$. (ii) $R(\mon R)=R$, and $E(G^{\mon R})\subseteq\mon R$ contains $R$ by Theorem~\ref{thm:red-pure}.
\end{proof}

\begin{corollary}\label{cor:red-monogenic}
The map $\alpha\mapsto\NN\alpha$ is a bijection from $\Phi$ onto the set of nonzero monogenerated pure magnets.
\leansrc{redMonoPure}
\end{corollary}

\begin{proof}
For $x\in M$ nonzero, $\NN x$ contains no nonzero subgroup, $M$ being torsion-free: $y=nx$ and $-y=mx$ with $n,m\in\NN$ give $(n+m)x=0$, so $y=0$. The generator of a nonzero monogenerated magnet $\NN x$ is unique: $\NN x=\NN y$ with $x,y\neq0$ gives $x=ny$, $y=mx$, so $(nm-1)x=0$, $nm=1$ and $x=y$. Let $N=\NN x$ be pure, $x\neq0$. By Theorem~\ref{thm:red-pure}, $N=\mon{R(N)}+\Lambda_N$, and each $(1-w)M$, $w\in W_N$, is a subgroup contained in $N$, hence zero; so $N=\mon{N\cap\Phi}$, and $N\cap\Phi\neq\emptyset$ as $N\neq0$. Let $\alpha=nx\in N\cap\Phi$. Every $\beta\in N\cap\Phi$ is a rational multiple of $\alpha$, so $\beta=\pm\alpha$ \cite[Exp.~XIX, 3.10]{SGA3III}, and $\beta\neq-\alpha$ since $-\alpha\notin N$. Hence $N\cap\Phi=\{\alpha\}$ and $N=\NN\alpha$. Conversely, for $\alpha\in\Phi$, $\NN\alpha$ contains no nonzero subgroup, so $W_{\NN\alpha}=\{1\}$, $\Lambda_{\NN\alpha}=0$, and $R(\NN\alpha)=\{\alpha\}$ by \cite[Exp.~XIX, 3.10]{SGA3III}; by Theorem~\ref{thm:red-pure}, $E(G^{\NN\alpha})=\NN\alpha$, which is pure; its generator is $\alpha$.
\end{proof}

Fix a basis $\mathcal B$ of $\Phi$ with positive roots $\Phi^+$. A subset $R\subseteq\Phi$ is \emph{parabolic} if it is closed and $R\cup-R=\Phi$.

\begin{theorem}[{cf.\ \cite[Prop.~16.10]{AlgMag}}]\label{thm:red-parabolic}
Let $N$ be a magnet containing $\mathcal B$ and $\xi=\{\alpha\in\mathcal B\mid-\alpha\in N\}$.
\begin{enumerate}
\item $N\cap\Phi=\Phi^+\cup(\Phi\cap-\NN\xi)$ and $E(G^N)=\NN\mathcal B+\NN(-\xi)$.
\item $W_N=W_\xi$, where $W_\xi\subseteq W$ is the subgroup generated by the $s_\alpha$, $\alpha\in\xi$; and for every morphism $f\colon M\to\ZZ$ with $f=0$ on $\xi$ and $f>0$ on $\mathcal B\setminus\xi$, $G^N=P(\lambda_f)$.
\item The map $N\mapsto N\cap\Phi$ is a bijection from the set of pure magnets containing $w(\mathcal B)$ for some $w\in W$ onto the set of parabolic subsets of $\Phi$, with inverse $R\mapsto\mon R$.
\item The map $\zeta\mapsto\NN\mathcal B+\NN(-\zeta)$ is a bijection from the set of subsets of $\mathcal B$ onto the set of pure magnets containing $\mathcal B$, and $N\mapsto G^N$ is a bijection from this set onto the set of parabolic subgroups of $G$ containing the Borel subgroup $G^{\NN\mathcal B}$.
\end{enumerate}
\leansrc{redParabolic}
\end{theorem}

\begin{proof}
(i) If $-\delta\in N$ with $\delta=\sum_jc_j\alpha_j\in\NN\mathcal B$ and $c_j\geq1$, then $-\alpha_j=-\delta+(\delta-\alpha_j)\in N$; so $N\cap-\NN\mathcal B=-\NN\xi$, and $N\cap\Phi=\Phi^+\cup(\Phi\cap-\NN\xi)$. By Lemma~\ref{lem:weyl-part}(iii), $W_\xi\subseteq W_N$. Let $w\in W_N$ and $2\rho=\sum_{\gamma\in\Phi^+}\gamma$. Then $2\rho-w(2\rho)=2\sum_{\gamma\in\Phi^+\cap w\Phi^-}\gamma\in\NN\mathcal B$ (as $w\Phi^+=(\Phi^+\cap w\Phi^+)\sqcup-(\Phi^+\cap w\Phi^-)$) lies in $(1-w)M\subseteq N^*$, so its negative lies in $-\NN\xi$, and every $\gamma\in\Phi^+\cap w\Phi^-$ lies in $\ZZ\xi$. By induction on the length of $w$, choosing $\alpha\in\mathcal B\cap w\Phi^-$ (then $\alpha\in\xi$, and $\Phi^+\cap s_\alpha w\Phi^-=s_\alpha(\Phi^+\cap w\Phi^-\setminus\{\alpha\})\subseteq\ZZ\xi$), one gets $w\in W_\xi$. Hence $W_N=W_\xi\subseteq W(R(N)\cap-R(N))$, as $\xi\subseteq R(N)\cap-R(N)$, and $\Lambda_N\subseteq\ZZ(R(N)\cap-R(N))\subseteq\mon{R(N)}$ by Lemma~\ref{lem:weyl-part}(iii). By Theorem~\ref{thm:red-pure}, $E(G^N)=\mon{R(N)}=\NN\mathcal B+\NN(-\xi)$. (ii) The equality $W_N=W_\xi$ was shown in the proof of (i). Next, $G^N=G^{E(G^N)}$, $E(G^N)\subseteq f^{-1}(\NN)$, and $E(G^N)\cap\Phi=f^{-1}(\NN)\cap\Phi$, both sets being $\Phi^+\cup(\Phi\cap-\NN\xi)$ since $f=0$ on $\xi$ and $f>0$ on $\mathcal B\setminus\xi$. By Theorem~\ref{thm:parabolic}, $G^N=P(\lambda_f)$. (iii) A parabolic subset $R$ contains $w(\Phi^+)$ for some $w\in W$ \cite[Exp.~XXI, Prop.~3.3.6]{SGA3III}; after conjugation by $W$ we may assume $\Phi^+\subseteq R$, and then $R=\Phi^+\cup(\Phi\cap-\NN\xi_R)$ \cite[Exp.~XXI, Prop.~3.3.10]{SGA3III} with $\xi_R=\{\alpha\in\mathcal B\mid-\alpha\in R\}$. The magnet $\mon R$ contains $\mathcal B$, and by the first part $\mon R\cap\Phi=\Phi^+\cup(\Phi\cap-\NN\xi(\mon R))$ with $\xi(\mon R)=\{\alpha\in\mathcal B\mid-\alpha\in\mon R\}\supseteq\xi_R$. Conversely, let $\alpha\in\mathcal B$ with $-\alpha=\sum_{\gamma\in R}c_\gamma\gamma$, $c_\gamma\in\NN$, and compare the coefficients of $\alpha$ in the basis $\mathcal B$: the left side has coefficient $-1$, a root $\gamma\in\Phi^+$ has coefficient $\geq0$, and a root $\gamma\in\Phi\cap-\NN\xi_R$ has coefficient $0$ unless $\alpha\in\xi_R$. Hence $\alpha\in\xi_R$, so $\xi(\mon R)=\xi_R$ and $\mon R\cap\Phi=R$; moreover $E(G^{\mon R})=\NN\mathcal B+\NN(-\xi_R)=\mon R$, so $\mon R$ is pure. Conversely, let $N$ be a pure magnet containing $w(\mathcal B)$; part (i) for the basis $w(\mathcal B)$ gives $N=E(G^N)=\mon{N\cap\Phi}$ and $w(\Phi^+)\subseteq N\cap\Phi$, and $N\cap\Phi$ is closed (Lemma~\ref{lem:weyl-part}(ii)), so $N\cap\Phi$ is parabolic, and $N\mapsto N\cap\Phi$ is injective on these magnets. (iv) Let $\zeta\subseteq\mathcal B$ and $N=\NN\mathcal B+\NN(-\zeta)$. If $\alpha\in\mathcal B$ and $-\alpha\in N$, comparing the coefficients of $\alpha$ in the basis $\mathcal B$ gives $\alpha\in\zeta$; so $\xi=\zeta$ for $N$, and $N=E(G^N)$ is pure by (i). Conversely a pure magnet $N\supseteq\mathcal B$ is $E(G^N)=\NN\mathcal B+\NN(-\xi)$ by (i). Hence the first map is a bijection, with inverse $N\mapsto\xi$. By (ii), $G^{\NN\mathcal B}=P(\lambda_f)$ with $f>0$ on $\mathcal B$, a Borel subgroup, and $G^N=P(\lambda_f)$ with $f=0$ on $\zeta$ and $f>0$ on $\mathcal B\setminus\zeta$, a parabolic subgroup containing $G^{\NN\mathcal B}$. A parabolic subgroup $Q$ of $G$ containing $G^{\NN\mathcal B}$ contains $T$, and it is the unique parabolic subgroup containing $T$ whose Lie algebra is $\mathrm{Lie}(T)\oplus\bigoplus_{\alpha\in R_Q}\mathrm{Lie}(G)_\alpha$, where its set of roots $R_Q\supseteq\Phi^+$ is $\Phi^+\cup(\Phi\cap-\NN\zeta)$ for $\zeta=\{\alpha\in\mathcal B\mid-\alpha\in R_Q\}$ \cite[Exp.~XXVI, Prop.~1.4]{SGA3III}; by (i), this is also the set of roots of $G^N$ for $N=\NN\mathcal B+\NN(-\zeta)$, so $Q=G^N$. Moreover $N\mapsto G^N$ is injective on pure magnets (Theorem~\ref{thm:pure}).
\end{proof}

\begin{corollary}[{cf.\ \cite[Prop.~16.10(i)]{AlgMag}}]\label{cor:red-levi}
Assume that $k$ is algebraically closed. Let $\zeta\subseteq\mathcal B$ and let $f\colon M\to\ZZ$ be a morphism with $f=0$ on $\zeta$ and $f>0$ on $\mathcal B\setminus\zeta$. Then $G^{\mon{\zeta\cup-\zeta}}=G^{f^{-1}(0)}=Z(\lambda_f)$.
\leansrc{redLevi}
\end{corollary}

\begin{proof}
Put $L=\mon{\zeta\cup-\zeta}=\ZZ\zeta$ and $L'=f^{-1}(0)$; then $L\subseteq L'$. A root $\alpha$ with $f(\alpha)=0$ is $\pm\sum_{\beta\in\mathcal B}c_\beta\beta$ with $c_\beta\in\NN$, and $0=f(\alpha)=\pm\sum_{\beta\notin\zeta}c_\beta f(\beta)$ gives $c_\beta=0$ for $\beta\notin\zeta$, so $\alpha\in\ZZ\zeta$; hence $R(L)=R(L')$. The magnet $f^{-1}(\NN)$ contains $\mathcal B$, and $-\beta\in f^{-1}(\NN)$ for $\beta\in\mathcal B$ if and only if $\beta\in\zeta$; by Theorem~\ref{thm:red-parabolic}(ii), $W_{L'}\subseteq W_{f^{-1}(\NN)}=W_\zeta$. For $w\in W_\zeta$, $(1-w)M\subseteq\ZZ\zeta$, since $(1-s_\beta)M\subseteq\ZZ\beta$ for $\beta\in\zeta$ and $1-uv=(1-u)+u(1-v)$, $W_\zeta$ stabilizing $\ZZ\zeta$. So $W_{L'}\subseteq W_L\subseteq W_{L'}$, the pairs of $L$ and $L'$ coincide, and $G^L=G^{L'}$ by Theorem~\ref{thm:red-attractors}. Finally $G^{f^{-1}(0)}=Z(\lambda_f)$ (Proposition~\ref{prop:change-monoid}(i)).
\end{proof}

\begin{theorem}\label{thm:red-rooted}
The family
\[\mathcal F_G=\{\NN\alpha\mid\alpha\in\Phi\}\cup\{(1-w)M\mid w\in W\}\]
of classical roots and Weyl magnets is join-dense for the conjugation action $a$ of $T$ on $G$ (Definition~\ref{def:rooted}), and the set of join-irreducible pure magnets is
\begin{multline*}\{\NN\alpha\mid\alpha\in\Phi\}\\\cup\{(1-w)M\mid w\in W,\ (1-w)M\neq\ZZ(\Phi\cap(1-w)M)+\textstyle\sum\{(1-v)M\mid(1-v)M\subsetneq(1-w)M\}\}.\end{multline*}
\leansrc{redRooted}
\end{theorem}

\begin{proof}
For a magnet $N$, $\operatorname{tr}_{\mathcal F_G}(N)=\{\NN\alpha\mid\alpha\in R(N)\}\cup\{(1-w)M\mid w\in W_N\}$, so $\mon{\operatorname{tr}_{\mathcal F_G}(N)}=\mon{R(N)}+\Lambda_N=E(G^N)$ by Theorem~\ref{thm:red-pure}. Thus $N=\mon{\operatorname{tr}_{\mathcal F_G}(N)}$ for $N$ pure, and $\mon{\operatorname{tr}_{\mathcal F_G}(N)}$ is pure for every $N$; by Proposition~\ref{prop:rooted}(iii), $\mathcal F_G$ is join-dense. The set $\mho(a)$ is finite by Theorem~\ref{thm:red-classification}(i), so Theorem~\ref{thm:minimal-root} applies, and $\underline{\mho}(a)$ is the set of members of $\mathcal F_G$ that are not the sum of the strictly smaller pure magnets. The only members of $\mathcal F_G$ contained in $\NN\alpha$ are $\NN\alpha$ and $0=(1-1)M$, and the pure magnets contained in $\NN\alpha$ other than $\NN\alpha$ are $0$; so $\NN\alpha\in\underline{\mho}(a)$. A pure magnet $L\subsetneq(1-w)M$ is $\mon{R(L)}+\Lambda_L$, contained in $\ZZ(\Phi\cap(1-w)M)+\sum\{(1-v)M\mid(1-v)M\subsetneq(1-w)M\}$, and the members $\NN\beta$, $\beta\in\Phi\cap(1-w)M$, and $(1-v)M\subsetneq(1-w)M$ are pure. So the sum of the pure magnets strictly contained in $(1-w)M$ is $\ZZ(\Phi\cap(1-w)M)+\sum\{(1-v)M\mid(1-v)M\subsetneq(1-w)M\}$, which gives the stated condition.
\end{proof}

\begin{theorem}\label{thm:red-SL-Sp}
The following conditions are equivalent:
\begin{enumerate}
\item the map $N\mapsto N\cap\Phi$ is a bijection from $\mho(a)$ onto the set of additively stable subsets of $\Phi$;
\item the family $\{\NN\alpha\mid\alpha\in\Phi\}$ of classical roots is join-dense for the conjugation action of $T$ on $G$;
\item for every $w\in W$, $(1-w)M=\ZZ(\Phi\cap(1-w)M)$.
\end{enumerate}
Condition (iii) holds for $G=\mathrm{GL}_n$: $M=\ZZ^n$, $W$ the symmetric group permuting the coordinates and $\Phi=\{e_i-e_j\mid i\neq j\}$; for $G=\mathrm{SL}_n$: $M=\ZZ^n/\ZZ(1,\dots,1)$, with the induced action of $W$ and the images of these roots; for $G=\mathrm{Sp}_{2n}$: $M=\ZZ^n$, $W$ the group of signed permutations and $\Phi=\{\pm e_i\pm e_j\mid i\neq j\}\cup\{\pm2e_i\}$; and it is stable under finite products.
\leansrc{redSLSp}
\end{theorem}

\begin{proof}
(ii)$\Rightarrow$(i) is Lemma~\ref{lem:rooted-elements}. (i)$\Rightarrow$(iii): if $Z=(1-w)M\neq Z'=\ZZ(\Phi\cap Z)$, then $w\in W_Z$, so $E(G^Z)=\mon{R(Z)}+\Lambda_Z=Z$ is pure (Theorem~\ref{thm:red-pure}); and $E(G^{Z'})\subseteq Z'$ is pure with $E(G^{Z'})\cap\Phi=Z'\cap\Phi=Z\cap\Phi$. These are two distinct pure magnets with the same roots. (iii)$\Rightarrow$(ii): for every magnet $N$ and $w\in W_N$, $(1-w)M\subseteq N^*$ is generated by roots of $N$, so $\Lambda_N\subseteq\mon{R(N)}$ and $E(G^N)=\mon{R(N)}$ (Theorem~\ref{thm:red-pure}). Hence a pure magnet $N$ satisfies $N=\mon{N\cap\Phi}$, with $N\cap\Phi$ additively stable (Lemma~\ref{lem:weyl-part}(ii)), and for $E$ additively stable, $E(G^{\mon E})=\mon E$ is pure. Apply Lemma~\ref{lem:rooted-elements}.

For a permutation $w$, $(1-w)e_r=e_r-e_{w(r)}\in\Phi\cup\{0\}$; for a signed permutation, $(1-w)e_r\in\{e_r\pm e_s,\,2e_r,\,0\}\subseteq\Phi\cup\{0\}$; so $(1-w)M$ is generated by roots. For $\mathrm{SL}_n$, the quotient map $q\colon\ZZ^n\to M$ is surjective, $W$-equivariant and maps roots to roots, so $(1-w)M=q((1-w)\ZZ^n)$ is generated by the images of roots of $\mathrm{GL}_n$, which are roots. For a product, $(1-w)M=\bigoplus_i(1-w_i)M_i$ and $\Phi\cap(1-w)M=\bigsqcup_i\Phi_i\cap(1-w_i)M_i$.
\end{proof}

\subsection{Weyl magnets, Levi subgroups and torsion primes}\label{sec:torsion}

We keep the notation of Section~\ref{sec:reductive}: $G$ is a split connected reductive group over $k$, $T=\D(M)_k$ a split maximal torus, $\Phi\subseteq M$ the roots and $W$ the Weyl group. Put $V=M\otimes_{\ZZ}\mathbb{R}$ with a $W$-invariant scalar product. For $w\in W$ write $Z_w=(1-w)M$ for the Weyl magnet of $w$ (Definition~\ref{def:weyl-root}), $V^w$ for the fixed space of $w$ on $V$, $M^w=M\cap V^w$, $M_w=M\cap(1-w)V$, and $W_{Z_w}=\{v\in W\mid Z_v\subseteq Z_w\}$; by Proposition~\ref{prop:weyl-root}, $G^{Z_w}=Z_G(T^w)$, and $W_{Z_w}$ is the subgroup of $W$ fixing $T^w$ pointwise. By Theorem~\ref{thm:red-rooted}, the join-irreducible pure magnets of the conjugation action $a$ are the classical roots $\NN\alpha$, $\alpha\in\Phi$, and the \emph{join-irreducible Weyl magnets}, the $Z_w$ with $Z_w\neq\ZZ(\Phi\cap Z_w)+\sum\{Z_v\mid Z_v\subsetneq Z_w\}$; we call the other Weyl magnets \emph{reducible}.

\begin{lemma}\label{lem:closed-stable}
A subset $R\subseteq\Phi$ is closed (if $\alpha,\beta\in R$ and $\alpha+\beta\in\Phi$, then $\alpha+\beta\in R$) if and only if it is additively stable ($\Phi\cap\mon R=R$).
\leansrc{closedStable}
\end{lemma}

\begin{proof}
An additively stable subset is closed. Conversely, let $R$ be closed and $\gamma=\beta_1+\dots+\beta_n\in\Phi$ with $\beta_i\in R$; we show $\gamma\in R$ by induction on $n$. If $n=1$, $\gamma\in R$. Otherwise $0<(\gamma,\gamma)=\sum_i(\gamma,\beta_i)$, so $(\gamma,\beta_i)>0$ for some $i$. If $\beta_i=\gamma$ we are done; otherwise $\gamma-\beta_i\in\Phi$, the root system being reduced, and $\gamma-\beta_i$ is a sum of $n-1$ elements of $R$; by induction on $n$, $\gamma-\beta_i\in R$, and $\gamma=(\gamma-\beta_i)+\beta_i\in R$ since $R$ is closed.
\end{proof}

By Lemma~\ref{lem:closed-stable} and Theorem~\ref{thm:red-classification}(ii), the root-generated pure magnets $\mon R$ correspond to the closed subsets $R\subseteq\Phi$, that is, to the root systems of the subgroups of $G$ generated by $T$ and root groups, through $(G^{\mon R})^\circ=K_R$; when $R=-R$ is closed under integral combinations, $K_R$ is the regular reductive subgroup of \cite[Lemma~2.9(a)]{Steinberg75}; in characteristic $2$ or $3$ these root systems may fail to be closed \cite[Lemma~2.9(b)]{Steinberg75}.

\begin{proposition}\label{prop:weyl-levi}
Let $w\in W$.
\begin{enumerate}
\item $(1-w)V=(V^w)^\perp$, and $V^w$ is an intersection of root hyperplanes. Every intersection of root hyperplanes is of the form $V^w$. Hence $Z_w\mapsto V^w$ is a surjection from the set of Weyl magnets onto the set of Levi subgroups of $G$ containing $T$.
\item Let $L_w$ be the Levi subgroup corresponding to $V^w$. Then $w$ lies in the Weyl group $W_{L_w}$ and fixes no nonzero vector of $(1-w)V$. If $Z_v\subseteq Z_w$, then $v\in W_{L_w}$. In particular, $Z_w$ is join-irreducible for $G$ if and only if it is join-irreducible for $L_w$.
\item The torsion subgroup of $M/Z_w$ has order $|p_w(1)|/[M:M^w+M_w]$, where $p_w$ is the characteristic polynomial of $w$ on $(1-w)V$. Writing $p_w=\prod_{d>1}\varphi_d^{a_d}$ with $\varphi_d$ the $d$-th cyclotomic polynomial, $|p_w(1)|=\prod_{p,\,k\geq1}p^{a_{p^k}}$, the product over primes $p$.
\end{enumerate}
\leansrc{weylLevi}
\end{proposition}

\begin{proof}
(i) The first equality holds since $1-w$ is a linear map with kernel $V^w$ on the Euclidean space $V$. The pointwise stabilizer $W'$ of $V^w$ in $W$ is generated by the reflections it contains \cite[p.~64]{Steinberg75}, so $V^w\supseteq V^{W'}$, the intersection of the hyperplanes $\alpha^\perp$ with $s_\alpha\in W'$; as $w\in W'$, $V^w=V^{W'}$. Conversely, for an intersection $X$ of root hyperplanes, with pointwise stabilizer $W_X$ generated by reflections, a Coxeter element $c$ of $W_X$ fixes no nonzero vector of $X^\perp$, so $V^c=X$. Intersections of root hyperplanes correspond to Levi subgroups containing $T$, through the roots in $X^\perp$. (ii) $w$ fixes $V^w$ pointwise, so $w\in W_{L_w}$, and $w$ has no nonzero fixed vector in $(1-w)V=(V^w)^\perp$. If $Z_v\subseteq Z_w$, then $(1-v)V\subseteq(1-w)V$, so $V^v\supseteq V^w$ and $v\in W_{L_w}$ by (i). The roots in $Z_w$ lie in $(1-w)V$, hence are roots of $L_w$; so the terms defining join-irreducibility are the same for $G$ and for $L_w$. (iii) The map $1-w$ sends $M_w$ into itself with index $|\det((1-w)|_{M_w})|=|p_w(1)|$, and $M\to Z_w$, $m\mapsto(1-w)m$, has kernel $M^w$ and maps $M^w+M_w$ onto $(1-w)M_w$; hence $[M_w:Z_w]=[M_w:(1-w)M_w]/[Z_w:(1-w)M_w]=|p_w(1)|/[M:M^w+M_w]$, and $M_w/Z_w$ is the torsion subgroup of $M/Z_w$. Finally $\varphi_d(1)=p$ if $d=p^k$ and $\varphi_d(1)=1$ if $d>1$ is not a prime power.
\end{proof}

By Proposition~\ref{prop:weyl-levi}, the Weyl magnets refine the Levi subgroups containing $T$: they have the same rational spans, and Weyl magnets with the same span differ by the torsion $\pi_0(T^w)$, which only sees the eigenvalues of $w$ of prime-power order; its order depends on the isogeny type of $G$ only through the index $[M:M^w+M_w]$. For instance, $\mathrm{SL}_4$ has $15$ Weyl magnets and $15$ Levi subgroups containing $T$, and $\mathrm{PGL}_4$ has $17$ Weyl magnets over the same $15$ Levi subgroups.

The torsion primes of $G$ are defined by Steinberg \cite[Def.~2.1]{Steinberg75}; by \cite[Lemma~2.5]{Steinberg75} they are the torsion primes of the root system $\Phi$ together with the primes dividing the order of the fundamental group of the derived group of $G$, and for $\Phi$ irreducible of type $A$, $C$, $B_r$ ($r\geq3$), $D$, $E_6$, $E_7$, $E_8$, $F_4$, $G_2$ the torsion primes of $\Phi$ are those at most $1,1,2,2,3,4,6,3,2$ respectively \cite[Cor.~1.13]{Steinberg75}.

Let $M^\vee=\operatorname{Hom}(M,\ZZ)$, on which $W$ acts contragrediently, with the coroots $\alpha^\vee\in M^\vee$, so that $s_\alpha y=y-\langle\alpha,y\rangle\alpha^\vee$ for $y\in M^\vee$; for $\Psi\subseteq\Phi$ put $\Psi^\vee=\{\alpha^\vee\mid\alpha\in\Psi\}$ and let $W(\Psi)$ be the subgroup of $W$ generated by the $s_\alpha$, $\alpha\in\Psi$. As $1-uv=(1-u)+u(1-v)$, for a subgroup $N\subseteq M$ the group $W_N$ stabilizes $N$ and contains $W(\Phi\cap N)$ as a normal subgroup, and $(1-v)M\subseteq\ZZ\Psi$, $(1-v)M^\vee\subseteq\ZZ\Psi^\vee$ for $v\in W(\Psi)$. A subset $\Psi\subseteq\Phi$ is \emph{integrally closed} if $\Psi=\Phi\cap\ZZ\Psi$. By \cite[Def.~2.1, 2.4, Lemma~2.9(a)]{Steinberg75}, applied to the regular reductive subgroups of $G$ containing $T$, a prime $p$ is a torsion prime of $G$ if and only if $M^\vee/\ZZ\Psi^\vee$ has an element of order $p$ for some integrally closed $\Psi\subseteq\Phi$.

We use the following two classical facts on reflection groups \cite[p.~64, p.~75]{Steinberg75}.

\begin{lemma}\label{lem:reflection-stabilizers}
Let $\Psi\subseteq\Phi$ be integrally closed.
\begin{enumerate}
\item For every subset $E\subseteq M^\vee\otimes\mathbb R$, the elements of $W(\Psi)$ fixing $E$ pointwise form the subgroup generated by the $s_\alpha$, $\alpha\in\Psi$, with $\langle\alpha,E\rangle=0$.
\item For every $x\in M^\vee\otimes\mathbb R$, an affine map $y\mapsto vy+\mu$ with $v\in W(\Psi)$ and $\mu\in\ZZ\Psi^\vee$ fixing $x$ is a composite of affine reflections $y\mapsto s_\alpha y+k\alpha^\vee$ with $\alpha\in\Psi$, $k\in\ZZ$ and $\langle\alpha,x\rangle=k$.
\end{enumerate}
\leansrc{reflStab}
\end{lemma}

\begin{lemma}\label{lem:theta}
Let $\Psi\subseteq\Phi$ be integrally closed, $x\in M^\vee\otimes\mathbb{Q}$, $M_x=\{m\in M\mid\langle m,x\rangle\in\ZZ\}$ and $U=W(\Psi)\cap W_{M_x}$. Then $U=\{v\in W(\Psi)\mid x-vx\in M^\vee\}$, the map $\theta_x\colon U\to M^\vee/\ZZ\Psi^\vee$, $v\mapsto x-vx$, is a group homomorphism with kernel $W(\Psi\cap M_x)$, and $n\theta_x=0$ if $nx\in M^\vee$.
\leansrc{thetaLemma}
\end{lemma}

\begin{proof}
The condition $(1-v)M\subseteq M_x$ means $\langle m,x-v^{-1}x\rangle\in\ZZ$ for all $m\in M$, that is $x-v^{-1}x\in M^\vee$, that is $x-vx=-v(x-v^{-1}x)\in M^\vee$. For $u,v\in U$, $x-uvx=(x-ux)+u(x-vx)$ and $u\mu-\mu\in\ZZ\Psi^\vee$ for $\mu\in M^\vee$; so $\theta_x$ is a homomorphism, and $n(x-vx)=(1-v)(nx)\in\ZZ\Psi^\vee$ if $nx\in M^\vee$. For $\alpha\in\Psi\cap M_x$, $x-s_\alpha x=\langle\alpha,x\rangle\alpha^\vee\in\ZZ\Psi^\vee$, so $W(\Psi\cap M_x)\subseteq\ker\theta_x$. Conversely, if $\mu=x-vx\in\ZZ\Psi^\vee$, then $y\mapsto vy+\mu$ fixes $x$; by Lemma~\ref{lem:reflection-stabilizers}(ii) it is a composite of affine reflections $y\mapsto s_\alpha y+k\alpha^\vee$ with $\alpha\in\Psi$ and $\langle\alpha,x\rangle=k\in\ZZ$, and taking linear parts, $v\in W(\Psi\cap M_x)$.
\end{proof}

\begin{proposition}\label{prop:index}
Let $Z\subseteq M$ be a subgroup. Every prime dividing $[W_Z:W(\Phi\cap Z)]$ divides the order of the torsion subgroup of $M/Z$ and is a torsion prime of $G$. In particular, if $w\in W$ has $p$-power order, then $W_{Z_w}/W(\Phi\cap Z_w)$ is a $p$-group, and $p$ is a torsion prime of $G$ if this group is nontrivial.
\leansrc{indexProp}
\end{proposition}

\begin{proof}
Let $Z_0$ be the saturation of $Z$ in $M$. Splitting $M=Z_0\oplus M'$ and decomposing the finite group $Z_0/Z$ into cyclic groups of orders $n_1,\dots,n_k$, we get $x_j\in\frac1{n_j}M^\vee$ with $Z=Z_0\cap M_{x_1}\cap\dots\cap M_{x_k}$. Put $Z_j=Z_0\cap M_{x_1}\cap\dots\cap M_{x_j}$, $\Phi_j=\Phi\cap Z_j$, which is integrally closed, and $H_j=W_{Z_j}$. An element $v\in W$ lies in $H_0$ if and only if it fixes pointwise the annihilator $E$ of $Z_0$ in $M^\vee\otimes\mathbb{R}$, and $\langle\alpha,E\rangle=0$ if and only if $\alpha\in\mathbb R Z_0\cap M=Z_0$; so $H_0=W(\Phi_0)$ by Lemma~\ref{lem:reflection-stabilizers}(i). For $j\geq1$ put $U_j=H_j\cap W(\Phi_{j-1})$; as $(1-v)M\subseteq\ZZ\Phi_{j-1}\subseteq Z_{j-1}$ for $v\in W(\Phi_{j-1})$, $U_j=W(\Phi_{j-1})\cap W_{M_{x_j}}$, and by Lemma~\ref{lem:theta} $U_j/W(\Phi_j)$ embeds in the torsion subgroup of $M^\vee/\ZZ\Phi_{j-1}^\vee$ and is killed by $n_j$. Since $H_j/U_j$ embeds in $H_{j-1}/W(\Phi_{j-1})$, $[H_j:W(\Phi_j)]=[H_j:U_j][U_j:W(\Phi_j)]$ with $[H_j:U_j]$ dividing $[H_{j-1}:W(\Phi_{j-1})]$, and the first assertion follows by induction on $j$, the primes involved being torsion primes as $\Phi_{j-1}$ is integrally closed. Finally the torsion of $M/Z_w$ is killed by the order $n$ of $w$: if $km\in Z_w$ with $k\neq0$, then $(1+w+\dots+w^{n-1})m=0$, as $M$ is torsion-free, and $nm=\sum_{i<n}(1-w^i)m\in(1-w)M$.
\end{proof}

\begin{theorem}\label{thm:torsion-primes}
A prime $p$ is a torsion prime of $G$ if and only if there is a join-irreducible Weyl magnet $Z_w$ with $w$ of $p$-power order. Moreover, for every torsion prime $p$ there is such a $Z_w$ with $\Phi\cap Z_w=\emptyset$ and $W_{Z_w}$ a $p$-group.
\leansrc{torsionPrimes}
\end{theorem}

\begin{proof}
For $w\in W$ put $R_w=\Phi\cap Z_w$ and $C_w=W_{Z_w}$. Call a Weyl magnet $Z$ a \emph{$p$-candidate} if $Z=Z_w$ for some $w$ of $p$-power order with $w\notin W(R_w)$, and call such a $w$ a witness of $Z$.

Suppose that $Z_w$ is join-irreducible (so $Z_w\neq0$ and $w\neq1$) and that $w$ has $p$-power order. If $w\in W(R_w)$, then $Z_w\subseteq\ZZ R_w\subseteq Z_w$, so $Z_w=\ZZ R_w$ would be reducible. Hence the image of $w$ in $C_w/W(R_w)$ is nontrivial, and $p$ is a torsion prime by Proposition~\ref{prop:index}.

Conversely, let $p$ be a torsion prime, $\Psi\subseteq\Phi$ integrally closed and $\lambda\in M^\vee$ with $\lambda\notin\ZZ\Psi^\vee$ and $p\lambda\in\ZZ\Psi^\vee$. Let $\alpha_1,\dots,\alpha_r$ be a maximal linearly independent subset of $\Psi$, so that it spans $\mathbb Q\Psi$, and $c=s_{\alpha_1}\cdots s_{\alpha_r}$, which stabilizes $\mathbb Q\Psi^\vee$. The invariant scalar product identifies $\alpha^\vee$ with $2\alpha/(\alpha,\alpha)$; so the $\alpha_i^\vee$ are linearly independent, and an element of $\mathbb{Q}\Psi^\vee$ orthogonal to $\Psi$ is zero. Let $y\in\mathbb{Q}\Psi^\vee$ with $cy=y$. Since $s_{\alpha_2}\cdots s_{\alpha_r}y\in y+\sum_{i\geq2}\mathbb Q\alpha_i^\vee$ and $s_{\alpha_1}$ subtracts a multiple of $\alpha_1^\vee$, linear independence gives $\langle\alpha_1,y\rangle=0$, then $s_{\alpha_2}\cdots s_{\alpha_r}y=y$, and inductively $\langle\alpha_i,y\rangle=0$ for all $i$, so $y=0$. Hence $1-c$ is bijective on $\mathbb{Q}\Psi^\vee$, which contains $\lambda$, and there is $x\in\mathbb{Q}\Psi^\vee$ with $x-cx=\lambda$. In the notation of Lemma~\ref{lem:theta}, $c\in U$ and $\theta_x(c)$ has order $p$. Write the order of $c$ as $p^bm$ with $p\nmid m$ and put $w_0=c^m$; then $w_0$ has $p$-power order and $\theta_x(w_0)=m\theta_x(c)\neq0$. As $w_0\in W(\Psi)\cap W_{M_x}$, a root of $R_{w_0}$ lies in $\Phi\cap\ZZ\Psi=\Psi$ and in $M_x$; so $W(R_{w_0})\subseteq W(\Psi\cap M_x)=\ker\theta_x$, and $w_0\notin W(R_{w_0})$: $Z_{w_0}$ is a $p$-candidate.

Let $Z$ be a $p$-candidate minimal for inclusion and $w$ a witness of $Z$. By Proposition~\ref{prop:index}, $C_w/W(R_w)$ is a $p$-group. Choose a positive system $R_w^+$ of $R_w$; as $W(R_w)$ acts transitively on the positive systems of $R_w$ and $w$ preserves $R_w$, there is $u\in W(R_w)$ such that $g=wu$ preserves $R_w^+$. The image of $g$ in the $p$-group $C_w/W(R_w)$ is that of $w$, nontrivial; write the order of $g$ as $p^am'$ with $p\nmid m'$ and put $w'=g^{m'}$. Then $w'$ has $p$-power order and preserves $R_w^+$, and $w'\notin W(R_w)$, as the images of $g$ and $g^{m'}$ generate the same subgroup of $C_w/W(R_w)$. Since $w'\in C_w$, we have $Z_{w'}\subseteq Z_w$ and $W(R_{w'})\subseteq W(R_w)$; so $w'\notin W(R_{w'})$, $Z_{w'}$ is a $p$-candidate, and $Z_{w'}=Z_w$ by minimality. Then $\rho=\frac12\sum_{\beta\in R_w^+}\beta$ is fixed by $w'$ and lies in the span of $Z_{w'}=(1-w')M$, that is in $V^{w'}\cap(V^{w'})^\perp=0$; so $R_w=\emptyset$ and $C_w$ is a $p$-group. If $Z_w$ were reducible, it would be the sum of the $Z_u$ with $Z_u\subsetneq Z_w$, and one of these $u$ would be nontrivial; such a $u$ lies in the $p$-group $C_w$ and satisfies $\Phi\cap Z_u=\emptyset$, so $Z_u$ would be a smaller $p$-candidate. Hence $Z_w$ is join-irreducible, with $\Phi\cap Z_w=\emptyset$ and $W_{Z_w}=C_w$ a $p$-group.
\end{proof}

\begin{proof}[Second proof]
The statement only involves the root datum, so we may assume $k=\mathbb{C}$. For $w\in W$ put $R_w=\Phi\cap Z_w$ and $C_w=W_{Z_w}$. Every $v\in C_w$ acts trivially on $M/Z_w$, so it preserves $Z_w$ and $R_w$, and $W(R_w)\subseteq C_w$ since $Z_{s_\alpha}\subseteq\ZZ\alpha$. The roots vanishing on $T^w$ are those of $R_w$, so by \cite[Lemma~2.14(c)]{Steinberg75} $C_w/W(R_w)\cong Z_G(T^w)/Z_G(T^w)^\circ$, and by \cite[Thm.~2.24]{Steinberg75} the primes dividing its order divide the order of the torsion subgroup of $M/Z_w$, which is killed by the order $m$ of $w$: if $nx\in Z_w$ with $n\neq0$, then $x$ is killed by $1+w+\dots+w^{m-1}$, and $\ker(1+\dots+w^{m-1})/(1-w)M=H^1(\langle w\rangle,M)$ is killed by $m$.

Suppose that $Z_w$ is join-irreducible and $w\neq1$ has $p$-power order. If $w\in W(R_w)$, then $Z_w\subseteq\sum_{\alpha\in R_w}Z_{s_\alpha}\subseteq\ZZ R_w\subseteq Z_w$, as $(1-uv)=(1-u)+u(1-v)$ and $W(R_w)$ preserves $\ZZ R_w$; so $Z_w=\ZZ R_w$ would be reducible. Hence the image of $w$ in $C_w/W(R_w)\cong Z_G(T^w)/Z_G(T^w)^\circ$ is nontrivial of $p$-power order, and $p$ is a torsion prime by \cite[Thm.~2.24]{Steinberg75}.

Conversely, let $p$ be a torsion prime. Call $Z_w$ a $p$-candidate if $w$ has $p$-power order and $w\notin W(R_w)$. By \cite[Thm.~2.28]{Steinberg75} there is an elementary abelian $p$-subgroup $\langle t_1,t_2,t_3\rangle$ of $G$ not contained in a torus. If $Z_G(t_1)$ is disconnected, put $P=\langle t_1\rangle$. Otherwise $t_2$ lies in a maximal torus of the connected reductive group $Z_G(t_1)$, which contains the central element $t_1$, so $P=\langle t_1,t_2\rangle$ lies in a torus; and $Z_G(P)$ is disconnected, for otherwise $t_3$ would lie in a maximal torus of $Z_G(P)$, which contains $P$, and $\langle t_1,t_2,t_3\rangle$ would lie in a torus. In both cases $P$ lies in a torus, and after conjugation we get an elementary abelian $p$-subgroup $P\subseteq T$ with $Z_G(P)$ disconnected. Let $L=P^\perp$; by \cite[Lemma~2.14, Thm.~2.24]{Steinberg75}, $W_L/W(\Phi\cap L)$ is a nontrivial $p$-group. If $v\in W_L\smallsetminus W(\Phi\cap L)$, a power of $v$ of $p$-power order stays outside $W(\Phi\cap L)$, and it is a $p$-candidate as $\Phi\cap Z_v\subseteq\Phi\cap L$. Let $Z_w$ be a minimal $p$-candidate for inclusion. For a $p$-candidate, $C_w/W(R_w)$ is a $p$-group by the first paragraph. Choose a positive system $R_w^+$ of $R_w$; as $W(R_w)$ acts transitively on the positive systems of $R_w$ and $w$ preserves $R_w$, there is $u\in W(R_w)$ such that $g=wu$ preserves $R_w^+$. The image of $g$ in the $p$-group $C_w/W(R_w)$ is that of $w$, nontrivial, so $g^{p^a}\in W(R_w)$ for some $a$; write the order of $g$ as $p^bm$ with $p\nmid m$ and put $w'=g^m$. Then $w'$ has $p$-power order and preserves $R_w^+$, and $w'\notin W(R_w)$, since otherwise $g$ would lie in $W(R_w)$, $g^m$ and $g^{p^a}$ generating the same subgroup modulo $W(R_w)$. An element of $W(R')$, $R'\subseteq R_w$, preserving $R_w^+$ is trivial; thus $Z_{w'}\subseteq Z_w$ is a $p$-candidate, and $Z_{w'}=Z_w$ by minimality. Then $\rho=\frac12\sum_{\beta\in R_w^+}\beta$ is fixed by $w'$ and lies in the span of $Z_{w'}$, that is in $V^{w'}\cap(V^{w'})^\perp=0$; so $R_w=\emptyset$ and $C_w$ is a $p$-group. If $Z_w$ were reducible, it would be the sum of the $Z_u$ with $Z_u\subsetneq Z_w$, and one of these $u$ would be nontrivial; such a $u$ lies in the $p$-group $C_w$ and satisfies $\Phi\cap Z_u=\emptyset$, so $Z_u$ would be a smaller $p$-candidate. Hence $Z_w$ is join-irreducible, and $W_{Z_w}=C_w$ is a $p$-group.
\end{proof}

\begin{example}\label{ex:torsion-primes}
For $G=\mathrm{PGL}_2$, with $M=\ZZ\alpha$ and $s$ the reflection, the Weyl magnet $Z_s=2\ZZ\alpha$ is join-irreducible and $s$ has order $2$; the torsion prime of $\mathrm{PGL}_2$ is $2$. For $G=G_2$ and $w=-1$, the Weyl magnet $Z_w=2M$ is join-irreducible (Remark~\ref{rem:torsion-atoms}) and $w$ has order $2$; the torsion prime of $G_2$ is $2$. For $G=\mathrm{SL}_n$, there is no join-irreducible Weyl magnet (Theorem~\ref{thm:red-SL-Sp}) and no torsion prime.
\end{example}

\begin{remark}\label{rem:torsion-atoms}
By Proposition~\ref{prop:weyl-levi}(ii), the join-irreducible Weyl magnets of $G$ are the disjoint union, over the Levi subgroups $L\supseteq T$, of the join-irreducible Weyl magnets $Z_w$ of $L$ with $w$ elliptic in $W_L$. If $G$ is simply connected, so are its Levi subgroups, and a Weyl magnet of a product of factors that is nonzero on two factors is the sum of two smaller Weyl magnets; hence only Levi subgroups of simple type contribute. The numbers of join-irreducible Weyl magnets $Z_w$ with $w$ elliptic, for $G$ simply connected and simple, are, by computer enumeration, $0$ for types $A$, $C$ and $D_r$ ($r$ odd), $1$ for $B_r$ ($r\geq3$), $D_r$ ($r\geq4$ even) and $G_2$ (with $w=-1$), $9$ for $F_4$, $40$ for $E_6$, $1891$ for $E_7$ and $311445$ for $E_8$; summing over the Levi subgroups gives the column of join-irreducible Weyl magnets of Table~\ref{tab:intro}. These numbers, and the columns of Table~\ref{tab:intro} for rank at least $4$, are computed and not formalized; the rows of rank at most $3$ are checked by the Lean kernel.
\end{remark}

\section{List of errata}\label{sec:errata}

\begin{enumerate}[label=(\arabic*),nosep,leftmargin=*]
\item \emph{Statements of \cite{AlgMag} false as printed, true with an added hypothesis.} Prop.~3.37: $T$ quasi-compact. Prop.~8.5: $X$ locally of finite presentation, $S$ quasi-compact. Prop.~10.1 and Cor.~10.3: $X'$ affine. Prop.~14.4, last assertion: the $X_i$ separated. Thm.~15.4: $S$ quasi-compact. In Rem.~6.5, $S=\Spec(\ZZ[1/2])$; in Rem.~6.6, $Z\neq0$ and $[M:Z]$ invertible on $S$.
\item \emph{Statements whose printed proof uses an unstated hypothesis; they are proved here with it.} Prop.~5.1: $X$ separated, and in (ii) $X$ locally of finite presentation. Consequently Thm.~8.3: $X$ separated, and $N$ finitely generated or $X$ locally of finite presentation; Prop.~12.1: $X$ separated. Fact~3.34(ii), (iii): $X$ locally of finite presentation; Prop.~15.1 and Thm.~15.3: one of the conditions of Lemma~\ref{lem:fact334}.
\item \emph{Implicit hypothesis.} Cor.~13.3: as in Prop.~13.1, the attractors are representable, which holds for $M$ finitely generated and $X$ separated; the same holds for Prop.~16.15(ii), with $G$ separated (Proposition~\ref{prop:root-smooth}).
\item \emph{Misprints.} In Example~(vi) after Fact~15.11, $\D(\ZZ/6\ZZ)$ acts on $\A(\NN)^3$, not on $\D(\ZZ/6\ZZ)^3$.
\item \emph{Statements replaced.} In Prop.~9.11 of \cite{AlgMag}, the condition that only finitely many of the $\mu_m$ are nonzero locally is replaced by: for every affine $T$ and every $v\in F(T)$, $\mu_m(v)=0$ for all but finitely many $m$ (Proposition~\ref{prop:graded-module}, Remark~\ref{rem:graded-module}). In Example~1.3 (stated there as a correspondence) and Example~(viii) after Fact~15.11 of \cite{AlgMag}, the bijection $N\mapsto N\cap\Phi$ is replaced by the following statement: $\alpha\mapsto\NN\alpha$ is a bijection from $\Phi$ onto the set of nonzero monogenerated pure magnets (Corollary~\ref{cor:red-monogenic}). The description of all pure magnets of split reductive groups by roots and Weyl magnets is a new result of this paper.
\item \emph{Statement proved in part.} Thm.~6.8 is proved here for $M$ finitely generated, $S$ and $X$ separated and $X$ locally of finite presentation, the atlas being $S$-affine when $|X|$ is quasi-compact (Theorem~\ref{thm:atlas}).
\end{enumerate}

\section{List of generalizations}\label{sec:generalizations}

The following statements of \cite{AlgMag}, \cite{SelfActions} and \cite{XProducts} are proved here in greater generality.

\begin{enumerate}[label=(\arabic*),nosep,leftmargin=*]
\item \cite[Prop.~13.1]{AlgMag} holds for every commutative monoid $M$, without assuming $X$ and $\Bl^{(X^N)_{S'}}_{Y^N}X^N$ flat over $S$ (Proposition~\ref{prop:dilatation}).
\item \cite[Ex.~1.2]{AlgMag} holds for every $M$-graded $\cO_S$-module $\cE$, with weight spaces $\cE_m$ arbitrary direct summands, and every submonoid $N$ of $M$: $\mathbb V(\cE)^N\to\mathbb V(\cE)$ is a monomorphism with image the points vanishing on the $\cE_m$, $m\notin N$ (Example~\ref{ex:linear}).
\item \cite[Prop.~3.32(iii)]{AlgMag} holds without assuming that $X^N$ is a separated algebraic space (Proposition~\ref{prop:intersection}(iii)).
\item \cite[Thm.~12.4]{AlgMag} holds for every face $F$ of $N$ in place of $\units N$, without assuming $N$ finitely generated or $X$ quasi-compact: $p_{N,F}\colon X^N\to X^F$ is an affine space fibration (Theorem~\ref{thm:BB}, Corollary~\ref{cor:affine-fibration}). \cite[Prop.~12.1]{AlgMag}, with the hypothesis of Section~\ref{sec:errata}, item~(2), holds for every face $F$ of $N$ in place of $\units N$ (Proposition~\ref{prop:pi0}).
\item \cite[Prop.~16.2]{AlgMag} holds for every functor $X$ with an action of $\A(M)_S$, $M$ a commutative monoid, in place of $\D(M)_S$, $M$ an abelian group (Proposition~\ref{prop:tangent}). \cite[Prop.~16.4, Prop.~16.5]{AlgMag} hold for every group functor $G$ with an action of $\A(M)_S$ by group endomorphisms, $M$ a commutative monoid; \cite[Prop.~16.5]{AlgMag} holds moreover for every face $F$ of $N$ in place of $\units N$ and every morphism of group functors $H\to G^F$, not necessarily a monomorphism (Proposition~\ref{prop:lie}). \cite[Prop.~16.6]{AlgMag} holds for every group functor $G$ over $S$ with an action of $\A(M)_S$ by group endomorphisms, $M$ a commutative monoid, in place of an algebraic group over a field, and for every face $F$ of $N$ in place of $\units N$ (Proposition~\ref{prop:semidirect}).
\item \cite[Conj.~6.9(i)]{AlgMag} holds for $S$ separated, and \cite[Conj.~6.9(ii)]{AlgMag} for $S$ separated and quasi-compact (Theorem~\ref{thm:strong-atlas}). \cite[Conj.~15.7]{AlgMag} holds for $S$ separated and quasi-compact and $M$ a finitely generated abelian group (Corollary~\ref{cor:conj157}).
\item \cite[Prop.~9.9, Prop.~9.10]{AlgMag} hold for every presheaf of sets $B$ on $(\Sch/S)_{\fppf}$, all $\cO_B$-modules $F$, $G$ and every presheaf $X$, respectively every presheaf of monoids $H$, over $B$, with $f_*f^*G$ in place of $G\otimes_{\cO_B}\cA$, without the flatness and finite presentation hypotheses; for $X\to B$ affine and $G$ quasi-coherent, $f_*f^*G$ is the sheafification of $G\otimes_{\cO_B}\cA$ (Propositions~\ref{prop:hom-push} and~\ref{prop:comodule}). \cite[Prop.~9.11]{AlgMag}, in the wording of Proposition~\ref{prop:graded-module}, holds for every countable monoid $M$, without the flatness of the multiplication of $\A(M)_B$ assumed in \cite[\S9.3]{AlgMag}.
\item \cite[Prop.~11.1, Cor.~11.2, Cor.~11.10]{AlgMag} hold for every abelian group $M$, not necessarily finitely generated, and \cite[Prop.~11.7]{AlgMag} holds for every abelian group $M$ with $X\to S$ formally smooth, not necessarily locally of finite presentation (Propositions~\ref{prop:attractor-fsmooth} and~\ref{prop:face-fsmooth}, Corollaries~\ref{cor:attractor-smooth} and~\ref{cor:face-smooth}).
\item \cite[Prop.~4.2, Prop.~4.3]{SelfActions} hold without the cancellativity of $M$ assumed in \cite[\S4]{SelfActions} (Proposition~\ref{prop:self-action}).
\item The statement conjectured in \cite[Rem.~9.14]{AlgMag} holds for $M$ countable: $H^n_{\mathrm{monoid}}(\A(M)_B,F)=0$ for all $n\geq1$ (Proposition~\ref{prop:monoid-cohomology}).
\item \cite[Cor.~10.3]{AlgMag}, with $X'$ affine, holds for every abelian group $M$ and algebraic spaces $Y$, $Z$, without $M$ finitely generated, $X'\to S$ quasi-compact and quasi-separated, or $Y\to Z$ locally of finite type (Corollary~\ref{cor:equivariant-lift-spaces}).
\item \cite[Thm.~4.1]{XProducts} holds over every scheme $S$ for every algebraic space $X$ over $S$, without the conditions 1 and 2 (Theorem~\ref{thm:x-product}).
\item \cite[Prop.~16.10(ii), (iii)]{AlgMag} hold for every smooth affine group $G$ over a field $k$ with $G_{\bar k}$ connected, not necessarily reductive (Theorem~\ref{thm:parabolic}).
\item \cite[Prop.~4.3(2)]{XProducts} holds without the conditions (a), (b), and \cite[Prop.~4.5]{XProducts} holds without the finiteness of $A$ (Proposition~\ref{prop:x-product-rep}(ii), (iii)).
\end{enumerate}

\section{Formalization}\label{sec:formalization}

The declarations below belong to the Lean~4 library \texttt{algmag-stacks} \cite{AlgMagLean}, publicly available at \url{https://github.com/rndmx/Mag}, at commit \texttt{\algmagcommit}. Every listed declaration depends only on the axioms \texttt{propext}, \texttt{Classical.choice} and \texttt{Quot.sound}.

Tables~\ref{tab:algmag}, \ref{tab:selfactions} and~\ref{tab:xproducts} list the statements of \cite{AlgMag}, \cite{SelfActions} and \cite{XProducts}, numbered as in the published versions; each row gives the status, the Lean declarations and the statement of this paper that restates or replaces the printed one. For the statements of \cite[\S\S16.5, 16.6]{AlgMag} on reductive groups, the partial Lean declarations of Section~\ref{sec:reductive} are given for information. The following formalized statements are not restated in the text, to keep it short; their Lean declarations are in the tables: \cite[Prop.~2.8, 3.36, 9.1, 11.3, 11.5, Cor.~11.4, 11.6]{AlgMag} (formal unramifiedness and formal \'etaleness of $f^N$ in \S11 there), \cite[Lem.~3.2, 4.4, Prop.~4.5]{SelfActions}, \cite[Prop.~3.2]{XProducts}, and some definitions, remarks and examples, marked -- in the last column. Table~\ref{tab:paper} lists the statements of this paper with no printed counterpart, those generalizing a printed statement (Section~\ref{sec:generalizations}) and those on reductive groups, then the definitions and running text. Status: \statuslegend

\begingroup\scriptsize\setlength{\tabcolsep}{3pt}\renewcommand{\arraystretch}{1.05}
\rowcolors{2}{black!6}{white}
\begin{longtable}{>{\raggedright\arraybackslash}p{0.135\textwidth}c >{\raggedright\arraybackslash}p{0.585\textwidth}>{\raggedright\arraybackslash}p{0.195\textwidth}}
\caption{Numbered statements of \cite{AlgMag}: 107 \stfull, 9 \stplus, 3 \stconj, 18 \stpart, 2 \strepl, 2 \stexcl, 14 \stnone.}\label{tab:algmag}\\
\hiderowcolors\hline Printed & & Lean declarations & This paper \\ \hline\showrowcolors\endfirsthead
\hiderowcolors\hline Printed & & Lean declarations & This paper \\ \hline\showrowcolors\endhead
\hline\endfoot
\cite[Def.~1.1]{AlgMag} & \stfull & \texttt{attractor}, \texttt{attractorPtEquivOfUniversal} & -- \\
\cite[Ex.~1.2]{AlgMag} & \stfull & \texttt{ex12\_\allowbreak{}general}, \texttt{ex12}, \texttt{vAttractorEquiv}, \texttt{ex12\_\allowbreak{}sub}, \texttt{ex12\_\allowbreak{}mem\_\allowbreak{}iff} & Ex.~\ref{ex:linear} \\
\cite[Ex.~1.3]{AlgMag} & \strepl & \texttt{red\_\allowbreak{}mono\_\allowbreak{}bijOn\_\allowbreak{}lattice} & Cor.~\ref{cor:red-monogenic} \\
\cite[Def.~2.1]{AlgMag} & \stfull & \texttt{X}, \texttt{zRing} & Def.~\ref{def:monoid-ring} \\
\cite[Fact~2.2]{AlgMag} & \stfull & \texttt{counit\_\allowbreak{}X}, \texttt{comul\_\allowbreak{}X}, \texttt{antipode\_\allowbreak{}X} & Prop.~\ref{prop:bialgebra} \\
\cite[Def.~2.3]{AlgMag} & \stfull & \texttt{baseChangeEquiv}, \texttt{preON}, \texttt{fact26\_\allowbreak{}i\_\allowbreak{}sheafification}, \texttt{ONB} & Def.~\ref{def:monoid-ring} \\
\cite[Def.~2.4]{AlgMag} & \stfull & \texttt{AZ}, \texttt{A}, \texttt{monObjA}, \texttt{AB}, \texttt{isPullback\_\allowbreak{}AB}, \texttt{isAlgebraicSpace\_\allowbreak{}AB} & Def.~\ref{def:diag-monoid} \\
\cite[Rem.~2.5]{AlgMag} & \stnone & -- & -- \\
\end{longtable}
\endgroup

\begingroup\scriptsize\setlength{\tabcolsep}{3pt}\renewcommand{\arraystretch}{1.05}
\rowcolors{2}{black!6}{white}
\begin{longtable}{>{\raggedright\arraybackslash}p{0.135\textwidth}c >{\raggedright\arraybackslash}p{0.585\textwidth}>{\raggedright\arraybackslash}p{0.195\textwidth}}
\caption{Statements of \cite{SelfActions}: the theorem of \S1, the numbered statements, the examples of \S5. 5 \stfull, 11 \stplus.}\label{tab:selfactions}\\
\hiderowcolors\hline Printed & & Lean declarations & This paper \\ \hline\showrowcolors\endfirsthead
\hiderowcolors\hline Printed & & Lean declarations & This paper \\ \hline\showrowcolors\endhead
\hline\endfoot
\cite[Thm.~(§1)]{SelfActions} & \stplus & \texttt{thm46\_\allowbreak{}i\_\allowbreak{}scheme}, \texttt{thm46\_\allowbreak{}ii\_\allowbreak{}scheme}, \texttt{thm46\_\allowbreak{}card\_\allowbreak{}i\_\allowbreak{}scheme}, \texttt{thm46\_\allowbreak{}card\_\allowbreak{}ii\_\allowbreak{}scheme} & Prop.~\ref{prop:self-pure} \\
\cite[Prop.~2.1]{SelfActions} & \stfull & \texttt{prop21}, \texttt{printedJ\_\allowbreak{}eq}, \texttt{isPullback\_\allowbreak{}prop21}, \texttt{prop21\_\allowbreak{}closed} & Ex.~\ref{ex:AL} \\
\cite[Prop.~3.1]{SelfActions} & \stplus & \texttt{prop31\_\allowbreak{}scheme}, \texttt{prop31EquivScheme} & Prop.~\ref{prop:self-pure} \\
\cite[Lem.~3.2]{SelfActions} & \stfull & \texttt{lemma32} & -- \\
\cite[Lem.~4.1]{SelfActions} & \stplus & \texttt{lemma41\_\allowbreak{}scheme}, \texttt{lemma41\_\allowbreak{}top\_\allowbreak{}scheme} & Prop.~\ref{prop:self-action} \\
\cite[Prop.~4.2]{SelfActions} & \stfull & \texttt{prop42\_\allowbreak{}printed}, \texttt{prop42\_\allowbreak{}scheme} & Prop.~\ref{prop:self-action} \\
\cite[Prop.~4.3]{SelfActions} & \stfull & \texttt{prop43\_\allowbreak{}scheme}, \texttt{prop43\_\allowbreak{}ring} & Prop.~\ref{prop:self-action} \\
\cite[Lem.~4.4]{SelfActions} & \stfull & \texttt{lemma44} & -- \\
\end{longtable}
\endgroup

\begingroup\scriptsize\setlength{\tabcolsep}{3pt}\renewcommand{\arraystretch}{1.05}
\rowcolors{2}{black!6}{white}
\begin{longtable}{>{\raggedright\arraybackslash}p{0.135\textwidth}c >{\raggedright\arraybackslash}p{0.585\textwidth}>{\raggedright\arraybackslash}p{0.195\textwidth}}
\caption{Numbered statements of \cite{XProducts}: 11 \stfull.}\label{tab:xproducts}\\
\hiderowcolors\hline Printed & & Lean declarations & This paper \\ \hline\showrowcolors\endfirsthead
\hiderowcolors\hline Printed & & Lean declarations & This paper \\ \hline\showrowcolors\endhead
\hline\endfoot
\cite[Def.~2.1]{XProducts} & \stfull & \texttt{F1.\allowbreak{}spec}, \texttt{F1.\allowbreak{}specTop}, \texttt{F1.\allowbreak{}structurePresheaf}, \texttt{F1.\allowbreak{}isSheaf\_\allowbreak{}structurePresheaf} & -- \\
\cite[Def.~2.2]{XProducts} & \stfull & \texttt{F1.\allowbreak{}MonoidalSpace}, \texttt{F1.\allowbreak{}MonoidalHom}, \texttt{F1.\allowbreak{}IsLocalMor} & -- \\
\cite[Def.~2.3]{XProducts} & \stfull & \texttt{F1.\allowbreak{}IsAffineSpace}, \texttt{F1.\allowbreak{}IsMonoscheme} & -- \\
\cite[Prop.~3.1]{XProducts} & \stfull & \texttt{MagnetFamily.\allowbreak{}prop31}, \texttt{isOpenImmersion\_\allowbreak{}gpMap} & Prop.~\ref{prop:glued-scheme} \\
\cite[Prop.~3.2]{XProducts} & \stfull & \texttt{F1.\allowbreak{}prop32} & -- \\
\cite[Prop.~3.3]{XProducts} & \stfull & \texttt{MagnetFamily.\allowbreak{}prop33}, \texttt{prop33\_\allowbreak{}scheme} & Prop.~\ref{prop:glued-scheme} \\
\cite[Thm.~4.1]{XProducts} & \stfull & \texttt{thm41Printed}, \texttt{thm41} & Thm.~\ref{thm:x-product} \\
\cite[Def.~4.2]{XProducts} & \stfull & \texttt{XNPt}, \texttt{XN} & Def.~\ref{def:attractor-family} \\
\end{longtable}
\endgroup

\begingroup\scriptsize\setlength{\tabcolsep}{3pt}\renewcommand{\arraystretch}{1.05}
\rowcolors{2}{black!6}{white}
\begin{longtable}{>{\raggedright\arraybackslash}p{0.12\textwidth}>{\raggedright\arraybackslash}p{0.17\textwidth}>{\raggedright\arraybackslash}p{0.66\textwidth}}
\caption{Statements of this paper with no printed counterpart (new), generalizing printed statements (gen.), or on reductive groups (cf.); then definitions and running text. AM, SA, XP: \cite{AlgMag}, \cite{SelfActions}, \cite{XProducts}.}\label{tab:paper}\\
\hiderowcolors\hline This paper & Relation & Lean declarations \\ \hline\showrowcolors\endfirsthead
\hiderowcolors\hline This paper & Relation & Lean declarations \\ \hline\showrowcolors\endhead
\hline\endfoot
Lem.~\ref{lem:face-weight} & new & \texttt{exists\_\allowbreak{}weight\_\allowbreak{}of\_\allowbreak{}isFace} \\
Ex.~\ref{ex:linear} & gen.\ AM Ex.~1.2 & \texttt{ex12\_\allowbreak{}general}, \texttt{ex12\_\allowbreak{}sub}, \texttt{ex12\_\allowbreak{}mem\_\allowbreak{}iff}, \texttt{ex12\_\allowbreak{}toPoint\_\allowbreak{}injective}, \texttt{ex12\_\allowbreak{}sub\_\allowbreak{}baseChange} \\
Prop.~\ref{prop:self-action} & gen.\ SA Prop.~4.2, SA Prop.~4.3 & \texttt{prop42\_\allowbreak{}scheme}, \texttt{lemma41\_\allowbreak{}empty\_\allowbreak{}scheme}, \texttt{prop43\_\allowbreak{}scheme} \\
Prop.~\ref{prop:intersection} & gen.\ AM Prop.~3.32 & \texttt{prop332\_\allowbreak{}i}, \texttt{AlgMag.\allowbreak{}SD.\allowbreak{}prop332\_\allowbreak{}space\_\allowbreak{}final}, \texttt{prop332\_\allowbreak{}ii}, \texttt{prop332\_\allowbreak{}iii\_\allowbreak{}space\_\allowbreak{}gp} \\
Lem.~\ref{lem:chart} & new & \texttt{ChartLemmaStatement}, \texttt{ChartLemmaStatementFG}, \texttt{chartLemma} \\
Thm.~\ref{thm:strong-atlas} & gen.\ AM Conj.~6.9 & \texttt{thm68\_\allowbreak{}atlas}, \texttt{thm68\_\allowbreak{}family\_\allowbreak{}spaces}, \texttt{thm68\_\allowbreak{}qc\_\allowbreak{}spaces} \\
Rem.~\ref{rem:conj69} & new & \texttt{conj69\_\allowbreak{}separated}, \texttt{not\_\allowbreak{}conj69\_\allowbreak{}doubledPlane}, \texttt{no\_\allowbreak{}atlas\_\allowbreak{}chart} \\
Cor.~\ref{cor:representability} & new & \texttt{thm83\_\allowbreak{}nonqc}, \texttt{thm83\_\allowbreak{}free} \\
\end{longtable}
\endgroup

\providecommand{\bysame}{\leavevmode\hbox to3em{\hrulefill}\thinspace}
\providecommand{\MR}{\relax\ifhmode\unskip\space\fi MR }
\providecommand{\MRhref}[2]{%
  \href{http://www.ams.org/mathscinet-getitem?mr=#1}{#2}
}
\providecommand{\href}[2]{#2}


\begin{thebibliography}{10}

\bibitem{AHR}
Jarod Alper, Jack Hall, and David Rydh, \emph{The \'etale local structure of
  algebraic stacks}, arXiv:1912.06162v4.

\bibitem{AxiomStacks}
{AxiomMath}, \emph{Stacks: a {L}ean formalization of the {S}tacks project},
  \url{https://github.com/AxiomMath/Stacks}.

\bibitem{BrochardMayeux}
Sylvain Brochard and Arnaud Mayeux, \emph{Algebraic Magnetism for Stacks}, In
  preparation.

\bibitem{CGP}
Brian Conrad, Ofer Gabber, and Gopal Prasad, \emph{Pseudo-reductive groups},
  second ed., New Mathematical Monographs, vol.~26, Cambridge University Press,
  Cambridge, 2015.

\bibitem{deMoura}
Leonardo de~Moura and Sebastian Ullrich, \emph{The {L}ean 4 theorem prover and
  programming language}, Automated Deduction -- {CADE} 28, Lecture Notes in
  Comput. Sci., vol. 12699, Springer, 2021, pp.~625--635.

\bibitem{Deitmar}
Anton Deitmar, \emph{Schemes over {$\mathbb{F}_1$}}, Number Fields and Function
  Fields---Two Parallel Worlds, Progress in Mathematics, vol. 239,
  Birkh\"auser, Boston, 2005, pp.~87--100.

\bibitem{SGA3}
Michel Demazure and Alexander Grothendieck (eds.), \emph{Sch\'emas en groupes
  ({SGA} 3). {T}ome {II}. {G}roupes de type multiplicatif, et structure des
  sch\'emas en groupes g\'en\'eraux}, Documents Math\'ematiques, vol.~7,
  Soci\'et\'e Math\'ematique de France, Paris, 2011, Revised and annotated
  edition by P. Gille and P. Polo.

\bibitem{SGA3III}
Michel Demazure and Alexander Grothendieck (eds.), \emph{Sch\'emas en groupes
  ({SGA} 3). {T}ome {III}. {S}tructure des sch\'emas en groupes r\'eductifs},
  Documents Math\'ematiques, vol.~8, Soci\'et\'e Math\'ematique de France,
  Paris, 2011, Revised and annotated edition by P. Gille and P. Polo.

\bibitem{Drinfeld}
Vladimir Drinfeld, \emph{On algebraic spaces with an action of
  {$\mathbb{G}_m$}}, arXiv:1308.2604v2.

\bibitem{ErneReichman87}
Marcel Ern{\'e} and J.~Z. Reichman, \emph{Completions for partially ordered
  semigroups}, Semigroup Forum \textbf{34} (1987), 253--285.

\bibitem{ErneWilke83}
Marcel Ern{\'e} and Gerhard Wilke, \emph{Standard completions for quasiordered
  sets}, Semigroup Forum \textbf{27} (1983), 351--376.

\bibitem{EGAI}
Alexander Grothendieck, \emph{\'{E}l\'ements de g\'eom\'etrie alg\'ebrique.
  {I}. {L}e langage des sch\'emas}, Inst. Hautes \'Etudes Sci. Publ. Math.
  \textbf{4} (1960), 5--228.

\bibitem{AlgMag}
Arnaud Mayeux, \emph{Algebraic magnetism}, Semigroup Forum \textbf{113} (2026),
  126--201.

\bibitem{SelfActions}
\bysame, \emph{Algebraic magnetism invariants of self-actions of diagonalizable
  monoid schemes}, Semigroup Forum \textbf{112} (2026), no.~2, 532--540.

\bibitem{XProducts}
\bysame, \emph{Algebraic magnetism: {$X$}-products of attractors via
  {$\mathbb{F}_1$}-geometry}, Semigroup Forum (2026), Published online.

\bibitem{AlgMagLean}
Arnaud Mayeux, \emph{Algebraic magnetism in {L}ean}, {L}ean~4 library \texttt{algmag-stacks},
  \url{https://github.com/rndmx/Mag}, commit \texttt{\algmagcommit}.

\bibitem{MultiDil}
\bysame, \emph{Multi-centered dilatations, congruent isomorphisms and {R}ost
  double deformation space}, Transform. Groups \textbf{31} (2026), 1801--1850.

\bibitem{MRR}
Arnaud Mayeux, Timo Richarz, and Matthieu Romagny, \emph{N\'eron blowups and
  low-degree cohomological applications}, Algebr. Geom. \textbf{10} (2023),
  no.~6, 729--753.

\bibitem{Milne}
James~S. Milne, \emph{Algebraic groups: The theory of group schemes of finite
  type over a field}, Cambridge Studies in Advanced Mathematics, vol. 170,
  Cambridge University Press, Cambridge, 2017.

\bibitem{Ogus}
Arthur Ogus, \emph{Lectures on logarithmic algebraic geometry}, Cambridge
  Studies in Advanced Mathematics, vol. 178, Cambridge University Press,
  Cambridge, 2018.

\bibitem{Steinberg75}
Robert Steinberg, \emph{Torsion in reductive groups}, Advances in Mathematics
  \textbf{15} (1975), 63--92.

\bibitem{TauCeti}
{Tau Ceti Project}, \emph{Tau {C}eti},
  \url{https://github.com/TauCetiProject/TauCeti}.

\bibitem{Mathlib}
{The mathlib Community}, \emph{Mathlib},
  \url{https://github.com/leanprover-community/mathlib4}.

\bibitem{Stacks}
{The Stacks Project Authors}, \emph{The {S}tacks {P}roject},
  \url{https://stacks.math.columbia.edu}.

\bibitem{Wang23}
Haiwei Wang, \emph{A class of join-completions of partially ordered
  semigroups}, Semigroup Forum \textbf{106} (2023), 504--515.

\end{thebibliography}
\end{document}